\documentclass[11pt,reqno]{amsart}
\usepackage{amsmath}
\usepackage{amssymb}
\usepackage{color}
\usepackage{eucal}
\usepackage{color}
\usepackage{tikz}
\usepackage{tikz-cd}
\usepackage[all,arc,poly]{xy}
\usetikzlibrary{calc}
\usetikzlibrary{arrows.meta} 
\usepackage{changes}
\usepackage{calligra}
\usepackage{bm}
\DeclareMathAlphabet{\mathcalligra}{T1}{calligra}{m}{n} 
\usetikzlibrary{matrix,arrows,decorations.markings}
\usetikzlibrary{decorations.pathmorphing}
\makeatletter
\@namedef{subjclassname@2020}{
  \textup{2020} Mathematics Subject Classification}
\makeatother
\newcommand{\circt}{{{\mathop{\hbox{$\bigcirc$\kern-7.5pt\raise-0.5pt\hbox{$t$}\kern3.8pt}}}}}

\newcommand{\circu}{{{\mathop{\hbox{$\bigcirc$\kern-8.5pt\raise0.5pt\hbox{$u$}\kern2.8pt}}}}}

\newcommand{\circv}{{{\mathop{\hbox{$\bigcirc$\kern-8.5pt\raise0.5pt\hbox{$v$}\kern2.8pt}}}}}

\newcommand{\circw}{{{\mathop{\hbox{$\bigcirc$\kern-9.5pt\raise0.5pt\hbox{$w$}\kern1.5pt}}}}}

\newcommand{\circs}{{{\mathop{\hbox{$\bigcirc$\kern-8.5pt\raise0.5pt\hbox{$s$}\kern3.5pt}}}}}

\numberwithin{equation}{section}
\theoremstyle{plain}
\newtheorem{Thm}{Theorem}[section]
\newtheorem{Prop}[Thm]{Proposition}
\newtheorem{Lemma}[Thm]{Lemma}

\theoremstyle{definition}
\newtheorem{Def}[Thm]{Definition}

\newtheorem{Sum}[Thm]{Summary}

\newtheorem{Rmk}[Thm]{Remark}

\DeclareSymbolFont{rsfscript}{OMS}{rsfs}{m}{n}
\DeclareSymbolFontAlphabet{\mathrsfs}{rsfscript}
\newcommand{\co}{\operatorname{co}}
\newcommand{\inv}{^{-1}}
\newcommand{\til}[1]{\widetilde {#1}}
\newcommand{\wh}[1]{\widehat{#1}}
\newcommand{\ol}[1]{\overline{#1}}
\newcommand{\rbl}{{\xrightarrow{\hspace*{1cm}}}}
\newcommand{\lbl}{{\xleftarrow{\hspace*{1cm}}}}
\newcommand{\nc}{\newcommand}
\nc{\reldot}[3]{{#1}\mathrel{\underset{#2}{\cdot}}{#3}}
\nc{\cT}{\mathrsfs T}
\nc{\bP}{\mathbb{P}}
\nc{\bO}{\mathbb{O}}
\nc{\bS}{\mathfrak{S}}
\nc{\bX}{\mathfrak{X}}
\nc{\bY}{\mathfrak{Y}}
\nc{\bT}{\mathfrak{U}}
\nc{\bV}{\mathfrak{V}}
\nc{\bC}{\mathfrak C}
\nc{\bK}{\mathfrak K}
\nc{\bN}{\mathbb{N}}
\nc{\mC}{\mathcal{C}}
\nc{\mQ}{\mathcal{Q}}
\nc{\mE}{\mathcal{E}}
\nc{\mF}{\mathcal{F}}
\nc{\mL}{\mathcal L}
\nc{\mB}{\mathcal B}
\nc{\mK}{\mathcal K}
\nc{\mH}{\mathcal H}
\nc{\mG}{\mathcal G}
\nc{\mA}{\mathcal A}
\nc{\mM}{\mathcal M}
\nc{\mT}{\mathcal T}
\nc{\mS}{\mathcal S}
\nc{\mU}{\mathcal U}
\nc{\mX}{\mathcal X}
\nc{\mY}{\mathcal Y}
\nc{\mZ}{\mathcal Z}
\nc{\mD}{\mathcal D}
\nc{\mJ}{\mathcal J}
\nc{\mO}{\mathcal O}
\nc{\mR}{\mathcal R}
\nc{\mW}{\mathcal W}
\nc{\mP}{\mathcal{P}}
\nc{\vep}{\varepsilon}
\nc{\mcu}{\mathcalligra U}
\nc{\mcv}{\mathcalligra v}
\nc{\mcw}{\mathcalligra W}
\nc{\mcy}{\mathcalligra Y}
\nc{\mcz}{\mathcalligra Z}
\nc{\mcx}{\mathcalligra X}
\nc{\mcs}{\mathcalligra S}
\nc{\sfCE}{\mathsf{CE}}
\nc{\sfCL}{\mathsf{CL}}
\nc{\nothing}{\rule{0em}{1ex}}
\nc{\nt}{\nothing}
\nc{\highnothing}{\rule{0em}{3ex}}
\nc{\hnt}{\highnothing}
\nc{\look}{\marginpar{$\bullet$}}
\nc{\ssc}{\scriptscriptstyle}
\nc{\pss}[1]{\nt^{\ssc \subsetneq}\! #1}
\nc{\red}[1]{\textcolor{red}{#1}}
\nc{\blue}[1]{\textcolor{blue}{#1}}
\nc{\black}[1]{\textcolor{black}{#1}}
\nc{\tilA}[1]{${#1}$-compatible $\til A$}
\nc{\cutout}[1]{} 

\begin{document} 

\title[{Ash et al.~revisited}]{{Finite approximation of free groups II:\\
		The Theorems of Ash, Herwig--Lascar and Ribes--Zalesskii --- revisited and strengthened}}
	\author{K.~Auinger, J.~Bitterlich and M.~Otto}%
	\address{Fakult\"at f\"ur Mathematik, Universit\"at Wien, Oskar-Morgenstern-Platz 1, A-1090 Wien, Austria}
	\email{karl.auinger@univie.ac.at}
	\address{   }
	\email{bitt.j@protonmail.com}
	\address{Department of Mathematics,
	Technische Universit\"at Darmstadt,
		Schlossgartenstrasse 7,
		D-64289 Darmstadt,
		Germany}
	\email{otto@mathematik.tu-darmstadt.de}    
\begin{abstract} Relations and interactions between the theorems of Ash, Herwig--Lascar and Ribes--Zalesskii are discussed and it is shown that these
three theorems are equivalent in the sense that each of them can be derived
from each other one.  Some strengthenings of these theorems are obtained
with the  use of groups provided by a construction of the third author. Evidence is given that these strengthenings are substantially stronger than the classical results. Yet, it turns out that both kinds of results can be interpreted as different instances of the same common scheme, namely as \emph{finite approximation of free groups}.
   	
\end{abstract}
\dedicatory{Dedicated to Prof.~Mikhail Volkov on the occasion of his 70th birthday}

\maketitle
\tableofcontents
\section{Introduction} The three celebrated theorems mentioned in the title 
were obtained in the 1990s. The Theorem of Ash~\cite{Ash} is about 
{graphs and inverse monoids}\footnote{By ``Ash's Theorem'', in this
paper we always mean its version for inverse monoids.}, the Theorem of
Herwig--Lascar~\cite{HL} is about  relational structures and the Theorem of
Ribes--Zalesskii~\cite{RibesZalesskii} is about the profinite topology of the 
free group. Each of them has had substantial influence in its own field:%
\footnote{{The high number of citations each one of the papers~\cite{Ash,HL,RibesZalesskii} has attracted (as documented in MathSciNet) may also be seen as an indicator of this.}}
Ash's Theorem {solved several outstanding problems in finite semigroup theory~\cite{Hencketal} and} led to the introduction and investigation of the concepts of
\emph{hyperdecidability} and of \emph{tameness} of pseudovarieties of
semigroups and monoids, with many
important and interesting questions and results, initiated by the Porto school: see e.g.~\cite{almeidahyperdecidability,almeidatamenes,
ACZreducibility,delgado, almzeit,costeix, almklim, ali1,ali2};
many authors have dealt with extensions and generalisations of the Theorem of Herwig--Lascar:
see e.g.~\cite{conant, evans, konecny, sin, sol, eppa classes} and the extensive list of papers cited in~\cite{eppa classes}; and also the Theorem of Ribes--Zalesskii has been widely
discussed: new proofs have been published and the theorem has been
generalised to other topologies of the free group~\cite{RibesZalesskiipro-p,AS,ASconstructiveRZ,geometryofprofinite,auinbors,zhai} and to other groups~\cite{you, coul, minas, ros, minasmin, minspri, she, min}.

From the very beginning various connections and interactions between these
theorems have also been discussed: for Ash versus Ribes--Zalesskii this is not
surprising since both were originally motivated by giving a solution to the same
problem, namely to confirm  the so-called \textsl{Rhodes type II conjecture},
a conjecture about finite monoids that was much discussed during the 1970s
and 1980s. It is known from Steinberg's paper~\cite{steinberginverseautomata} that the Theorem
of Ribes--Zalesskii is equivalent to a special case of Ash's Theorem
(namely to the case of cycle graphs). Likewise, the fact that the 
Ribes--Zalesskii Theorem is a consequence of the Herwig--Lascar Theorem
has been one of the motivating examples in the paper by Herwig and
Lascar~\cite{HL} (the term \emph{profinite topology on free groups}
even is displayed in the title).  In addition, they presented a
version of their theorem --- let us call it the \emph{group-theoretic
version} of their theorem (Theorem~3.3 in~\cite{HL}) --- which
dealt with solvability in the free group of equations of a certain
type. 
This they announced as being equivalent to the
\emph{model-theoretic version} of their theorem: 
Theorem~3.2 in~\cite{HL}, 
which is usually referred to as
\textsl{the} Herwig--Lascar Theorem. 
A special case of the
group-theoretic version turns out to be essentially the same
as the Ribes--Zalesskii Theorem.

The overall conviction (at least in the semigroup community) is that the Theorem
of Ash and the Theorem of Herwig--Lascar are equivalent. Indeed, Almeida and
Delgado~\cite{almeidadelgado1,almeidadelgado2} proved that Ash's Theorem is
equivalent to the above-mentioned group-theoretic version of the theorem:
Theorem~3.3 of~\cite{HL}. But, strictly speaking, a proof that the Theorem of
Ash implies the model-theoretic version of the Theorem of Herwig--Lascar
(Theorem 3.2 of~\cite{HL}) has ---  to the best of the authors' knowledge
--- not been published, since a proof of the implication

\hnt 
\centerline{group-theoretic version $\Longrightarrow$ model-theoretic version} 
\hnt
\noindent is not published in~\cite{HL} nor elsewhere. So, the equivalence of Ash with Herwig--Lascar is known, but no full proof is published. Ironically, in this context we also have the opposite case: the fact that the Theorem of Ribes--Zalesskii implies the Theorem of Ash  (and hence the Theorem of Herwig--Lascar, at least in its group-theoretic version) seems to be not commonly acknowledged or recognised, yet a full  proof of this can be retrieved from the literature. In fact, it is one of the key steps in Ash's paper that the general case of his theorem follows from the special case for cycle graphs (and hence from the Ribes--Zalesskii Theorem). 

The first motivation for the present paper is to clean up this kind of mess. In Section~\ref{sec: revisited} we present full proofs for the chain of implications

\hnt
\centerline{Herwig--Lascar $\Longrightarrow$ Ribes--Zalesskii $\Longrightarrow$ Ash $\Longrightarrow$ Herwig--Lascar.}
\hnt
\noindent The most difficult step  is the second implication, the first and the third are comparatively easy. The proof of this second step is essentially taken from Ash's paper, but written up in more detail and somehow polished so that it is accessible, hopefully, also to non-specialists. The method 
of that proof --- which might be called \emph{arborisation of a labelled graph with respect to an inverse monoid} 
--- seems to be interesting in its own right. It will also be 
used in Section~\ref{sec:strengthened}.
It therefore seems justified to elaborate on 
the method of that proof in {Section~\ref{sec: revisited}}.
The proof of the first implication is inspired by  the one in~\cite{HL}. But, on the one hand, it is enriched with some pointers to the group-theoretic background: it is embedded into the context of Stallings and Schreier graphs, which seem to be \textsl{the} natural environment of these arguments. On the other hand, a slightly more general assertion is proved, which prepares the ground for a direct proof of the implication 
\hnt
\centerline{(model-theoretic) Herwig--Lascar $\Longrightarrow$ Ash } 
\hnt
which we also include (Section~\ref{subsec:HL===>Ash}). The idea for the third implication is taken from the second author's doctoral dissertation~\cite{bitterlichdiss}. 

A further intention of the discussion of the subject Ash/Herwig--Lascar/ Ribes--Zalesskii (i.e.~the content of Section~\ref{sec: revisited}) is to direct  the attention of readers interested in model-theoretic extension problems to the relevance of inverse monoids and to Ash's work. 
More importantly, we also intend to shed new
light on the context of the Herwig--Lascar result. Indeed, our treatment  in 
Section~\ref{sec: revisited} leads to a new interpretation of these model-theoretic extension problems, namely   as instances of \emph{finite approximation of free extensions} (cf.~Proposition~\ref{prop:EPPAin4.3.1}). 
 {As a by-product we get that the mathematical principle encapsulated in the results of Ash/Herwig--Lascar/Ribes--Zalesskii also appears in various other guises: there are three equivalent formulations in terms of equations over the free group (Propositions \ref{prop:HLgroupformulationABO}, \ref{prop:ashalmeida}, \ref{prop:HLgroupformulation}), an equivalent one in terms of finite groups (Proposition \ref{prop:ash-pure}) and another equivalent one in terms of unary algebras (Proposition \ref{prop:abstract-extension}).}

Our main motivation 
to recast the ``old stories'' lies in Section~\ref{sec:strengthened},
however. There we present some strengthening of the theory that can be achieved by the application of a construction  established by the third author~\cite{otto3,ABO}. The main result is Theorem~\ref{thm:otto-ash}, a strengthening of Ash's Theorem. 
 {As in the classical case, the strengthened version also appears in various equivalent formulations: most importantly we have a model-theoretic strengthening of Herwig--Lascar (Theorem \ref{thm: HL strengthened}); apart from that we have formulated two equivalent versions in terms of equations over the free group (Propositions \ref{prop:ashalmeidastrengthened} and \ref{prop:HLgroupformulationstrengthened}), one equivalent version in terms of finite groups (Proposition \ref{prop:otto-ash-pure}), and again one in terms of unary algebras (Proposition \ref{prop:strengthening-abstract-extension}).}

In Section~\ref{sec:prelims} we present an extensive aggregation of preliminaries, which should make the part of the paper concerning the old results as self-contained as possible.  For the strengthenings obtained in Section~\ref{sec:strengthened} we require non-trivial concepts and results from~\cite{ABO}. These are collected in Section~\ref{sec:prep} while Section~\ref{sec:proofs} contains the proof of Theorem~\ref{thm:otto-ash}. Finally, in Section~\ref{sec:strictly stronger} the common essence of as well as the 
differences 
between the old results (Section~\ref{sec: revisited}) and the new ones (Section~\ref{sec:strengthened}) are analysed. 
The comparison between the two kinds of results becomes most transparent in Propositions~\ref{prop:ash-pure} and~\ref{prop:otto-ash-pure}. 
The common essence lies in the fact that both results deal with \emph{finite approximation of free groups},\footnote{What  ``finite approximation of free groups'' is from 
the 
point of view of group theory is ``finite approximation of free extensions'' from a model-theoretic perspective. 
In the model-theoretic context, the comparison between the old and the new results 
becomes most evident in
Proposition~\ref{prop:EPPAin4.3.1} versus Theorem~\ref{thm: HL strengthened}.}  while the difference 
lies in 
the fact that this approximation is expressed in terms of some condition that is 
apparently stronger in Proposition~\ref{prop:otto-ash-pure}. A proof that the condition mentioned is \textsl{really} stronger (rather than just seemingly stronger) is presented in the second part of Section~\ref{sec:strictly stronger}.

\section{Preliminaries}\label{sec:prelims} We present a relatively long section of preliminaries to render the paper as self-contained as possible and  (hopefully)  accessible to non-specialists. 
\subsection{Inverse monoids} A monoid $M$ is \emph{inverse} if every
element $x\in M$
admits a \textsl{unique} element $x^{-1}$, called the \emph{inverse} of $x$, satisfying $xx^{-1}x=x$ and $x^{-1}xx^{-1}=x^{-1}$. This gives rise to a unary operation ${}^{-1}\colon M\to M$ and an inverse monoid may equivalently  be defined as an algebraic structure $(M;\cdot,{}^{-1},1)$ with $\cdot$ an associative binary operation, $1$ a neutral element with respect to $\cdot$ and a unary operation ${}^{-1}$ satisfying the identities (or laws)
\[{(x^{-1})}^{-1}=x,\ (xy)^{-1}=y^{-1}x^{-1},\ xx^{-1}x=x \mbox{ and }xx^{-1}yy^{-1}=yy^{-1}xx^{-1}.\] 
In particular, the class of all inverse monoids forms a variety of
algebraic structures (in the sense of universal algebra), the variety
of all groups $(G;\cdot,{}^{-1},1)$ being a subvariety. Both, the class of all groups as well as the class of all inverse monoids are varieties of \emph{involutory monoids},  were the variety of all involutory monoids $(M;\cdot,\inv,1)$ is defined by the identities ${(x\inv)}\inv=x$ and $(xy)\inv =y\inv x\inv$.

By the Wagner--Preston Theorem~\cite[Chapter 1, Theorem 1]{Lawson}, inverse monoids may as well be characterised as monoids of partial bijections on a set, closed under composition of partial mappings and inversion. Therefore, while groups  model symmetries of mathematical structures,  inverse monoids (or semigroups) model partial symmetries, that is, symmetries between substructures of mathematical structures. It is therefore not surprising that inverse monoids could play some r\^ole in the context of extending partial automorphisms of mathematical structures.

Every inverse monoid $M$ is equipped with a partial order $\le$, the
\emph{natural order},  defined by $x\le y$ if and only if $x=ye$ for
some idempotent $e$ of $M$ (this is equivalent to $x=fy$ for some
idempotent $f$ of $M$). 
If an inverse monoid $M$ is represented as a monoid of partial
bijections, then
the idempotents of $M$ are exactly the restrictions of the identity map. 
Hence any two idempotents commute and therefore the set of all idempotents of $M$ forms a (commutative) submonoid of $M$. In addition, for every element $x$ of $M$ and idempotent $e$, the ``conjugate'' element $x\inv ex$ is again an idempotent. For $x,y\in M$ we have $x\le y$ if and only if $x\subseteq y$, that is, $x$ is a restriction of $y$.
The order is compatible with the binary
operation, that is,  from $x\le y$ and $u\le v$ it follows that $xu\le yv$ for all $x,y,u,v$, and with inversion of $M$, that is, $x\le y$ implies $x^{-1}\le y^{-1}$. 
For further information on inverse monoids the reader is referred to the monographs by Petrich~\cite{Petrich} and Lawson~\cite{Lawson}.

\subsection{$A$-generated inverse monoids}
Throughout, for any non-empty set $X$ (of letters, of edges, etc.) we
let $X\inv:=\{x\inv \colon x\in X\}$ be a disjoint copy of $X$
consisting of formal inverses of the elements of $X$, and set 
$\til{X}:=X\cup X\inv$.
The mapping $x\mapsto x\inv$ is extended to an
involution of $\til{X}$ by setting ${(x\inv )}\inv=x$, for all $x\in
X$. We let $\til{X}^*$ be the free monoid
over $\til{X}$, which, subject to the operation $(x_1\cdots x_n)\inv=x_n\inv\cdots
x_1\inv$ (where $x_i\in \til{X}$), is the \emph{free involutory
	monoid} over $X$.
The elements of $\til{X}^*$ are called \emph{words over}
$\til{X}$, and  we identify the empty word with~$1$ in this context.
A word $w\in \til X^*$ is \emph{reduced} if it does not contain any
factor of the form $xx^{-1}$ for $x\in \til  X$. Successive deletion of such factors in a word $w$ until no more such factor is present leads to the \emph{reduced form} $\mathrm{red}(w)$ of $w$, 
which is well-defined. 
A word $w\in \til A^*$ with $\mathrm{red}(w)=1$ is called a \emph{Dyck word}.

We fix a non-empty set $A$ (called alphabet in this context).
An inverse monoid $M$ together with a (not necessarily injective)
mapping $i_M\colon A\to M$ (called \emph{assignment function}) is an
\emph{$A$-generated} inverse monoid if $M$ is generated by $i_M(A)$
as an inverse monoid, that is, generated
with respect to the operations $1,\cdot,{}^{-1}$. For every congruence
$\rho$ of an $A$-generated inverse monoid $M$, the quotient $M/\rho$
is  $A$-generated with respect to the map $i_{M/\rho}=\pi_\rho\circ
i_M$ where $\pi_\rho$ is the projection $M\to M/\rho$.
A \emph{morphism $\psi$ from the $A$-generated inverse monoid $M$ to the $A$-generated inverse monoid $N$}
is a homomorphism $M\to N$ {respecting generators from $A$}, that is,
satisfying $i_N=\psi\circ i_M$.
If it exists, such a morphism is unique and surjective and is called 
\emph{canonical morphism},
denoted $\psi\colon M\twoheadrightarrow N$. 
On a more formal level, an $A$-generated inverse monoid is an algebraic structure of the form $(M;\cdot,\inv,1,A)$ where every symbol $a\in A$ is interpreted in $M$ as a constant (that is, as a nullary operation) via the assignment function $i_M$. Canonical morphisms of $A$-generated inverse monoids then are just homomorphisms of algebraic structures in the signature $\{\cdot,\inv,1\}\cup A$.

If $M\twoheadrightarrow N$ then $M$ is an
\emph{expansion of $N$}. In our usage, the term ``expansion'' just concerns the relationship between two individual $A$-generated inverse monoids $M$ and $N$. This somehow deviates from the widespread use of that term standing for a functor on certain categories of monoids. The \emph{direct product} of $A$-generated inverse monoids $M_1,\dots, M_k$ is the $A$-generated inverse monoid with the diagonal mapping 
$a\mapsto (i_{M_1}(a),\dots,i_{M_k}(a))$ of $A$ into the full Cartesian
product as the natural assignment function; we will denote this product by $M_1\times \cdots\times M_k$, which in general is a proper submonoid of the full Cartesian product. 
The special case of $A$-generated groups will play a significant r\^ole in this paper. This concerns, in particular, the free $A$-generated group which will be denoted $F$, throughout the paper. The elements of $F$ will usually be identified with the reduced words over $\til A^*$.

As already mentioned, the assignment function is not necessarily injective, and, what is more, some generators may even be sent to the identity element of $M$.
This is not a deficiency, but rather is adequate
in our context,
since we want the quotient of an $A$-generated structure to be again
$A$-generated.

The assignment function $i_M$ is usually not explicitly mentioned; it uniquely extends 
to a homomorphism $[\ \ ]_M\colon \tilde{A}^*\to M$ (of involutory monoids). For every word $p\in \tilde{A}^*$, $[p]_M$  is the \emph{value of  $p$ in $M$} or simply the \emph{$M$-value of $p$}. For two words $p,q\in \til{A}^*$, the $A$-generated inverse monoid $M$ \emph{satisfies the relation} $p=q$ if $[p]_M=[q]_M$, in which case the words $p$ and $q$ are \emph{$M$-related} or \emph{$M$-equivalent}, while $M$ \emph{avoids the relation} $p=q$ if $[p]_M\ne[q]_M$.

\subsection{Graphs}\label{subsec:graphs}
In this paper, we consider {the Serre definition~\cite{Serre}} of {graph}
structures, admitting multiple directed edges between pairs of
vertices and  directed loop edges at individual vertices. In the literature, such
structures are often called \emph{{multidigraphs,} directed multigraphs} or
\emph{quivers}.
The following formalisation is convenient for our purposes. 
A \emph{graph} $\mK$ is a structure of the form
$(V\cup K;\alpha,\omega,{}^{-1})$ with $V$ its set of \emph{vertices}, $K$ its set of \emph{edges} (disjoint from $V$),  with \emph{incidence functions} $\alpha\colon K\to V$ and $\omega\colon K\to V$, 
specifying, for each edge $e$ (including loop edges), the \emph{initial} vertex $\alpha e$ and the \emph{terminal} vertex $\omega e$,
and \emph{involution} ${}\inv\colon K\to K$ satisfying $\alpha e=\omega e\inv$, $\omega e=\alpha e\inv$ and $e\ne e^{-1}$ for every edge $e\in K$. Instead of \emph{initial/terminal vertex} the terms \emph{source/target} are also used in the literature. One should think of an edge $e$ with initial vertex $\alpha e=u$ and terminal vertex $\omega e=v$ in ``geometric'' terms as $e\colon \underset{u}{\bullet}\!\overset{}{\rbl}\underset{v}{\bullet}$ and of its inverse $e\inv\colon \underset{u}{\bullet}\!\overset{}{\lbl}\underset{v}{\bullet}$ as ``the same edge but traversed in the opposite direction''. 
The pair $\{e,e\inv\}$ consisting of an edge $e$ and its inverse $e\inv$ is sometimes called a \emph{geometric edge} of the graph $\mK$ and is denoted ${e^\bullet}$, 
the set of all geometric edges is denoted  ${K^\bullet}$. The \emph{degree} of a vertex $v$ in a graph $\mK$ is the number of geometric edges in $\mK$ incident with $v$; this is the same as the number of edges $e$ having initial vertex $\alpha e=v$ and also the number of edges $e$ having terminal vertex $\omega e=v$. For a given graph $\mK$ we sometimes denote the set of its vertices by $\mathrm{V}(\mK)$, the set of its edges by $\mathrm{K}(\mK)$.
A graph $(V\cup K;\alpha,\omega,{}^{-1})$ is \emph{oriented} if the edge set $K$ is partitioned as $K=E\cup E^{-1}=\til{E}$ such that every ${}^{-1}$-orbit contains exactly one element of $E$ and one of $E^{-1}$; the edges in $E$ are the \emph{positive} or \emph{positively oriented} edges, those in $E^{-1}$ the \emph{negative} or \emph{negatively oriented} ones.

A \emph{subgraph} of the graph $\mK$ is a
substructure that is induced over a 
subset of $V\cup K$ which is
closed under the operations $\alpha$ and ${}\inv$ (and therefore also
under $\omega$). In particular, every non-empty subset $S\subseteq V\cup K$
generates a unique
minimal subgraph $\langle S\rangle$ of $\mK$ containing $S$,
which is the \emph{subgraph of $\mK$ spanned by $S$}. Likewise, for graphs $\mK=(V\cup K;\alpha,\omega,\inv)$ and $\mL=(W\cup L;\alpha,\omega,\inv)$ a \emph{morphism} $\varphi\colon \mK\to \mL$ is a mapping $\varphi\colon V\cup K\to W\cup L$ mapping vertices to vertices, edges to edges that is compatible with the operations $\alpha$ and $\inv$ (and therefore also $\omega$). A \emph{congruence}
$\Theta$ on the  graph $\mK=(V\cup K;\alpha,\omega,{}^{-1})$ is an equivalence relation $\Theta=\Theta_V\cup \Theta_K$ with $\Theta_V$ and $\Theta_K$ being equivalence relations on $V$ and $K$, respectively,  compatible with the operations $\alpha$ and ${}^{-1}$ (therefore also $\omega$):
\[e\mathrel{\Theta_K}f\Longrightarrow \alpha e\mathrel{\Theta_V}\alpha f,\  \omega e\mathrel{\Theta_V}\omega f,\ e^{-1}\mathrel{\Theta_K}f^{-1}\mbox{ for all }e,f\in K\]
and the concept of the quotient graph $\mK/\Theta$ is the obvious one, with the usual homomorphism theorem.
A \emph{non-empty path} $\pi$ in $\mK$ is a sequence $\pi=e_1e_2\cdots e_n$ ($n\ge 1$) of \emph{consecutive} edges (that is $\omega e_i=\alpha e_{i+1}$ for all $1\le i< n$); we set $\alpha\pi:=\alpha e_1$ and $\omega \pi=\omega e_n$ (denoting the initial and terminal vertices of the path $\pi$); the \emph{inverse path} $\pi\inv$ is the path $\pi\inv:=e_n\inv\cdots e_1\inv$; it has initial vertex $\alpha\pi\inv=\omega\pi$ and terminal vertex $\omega\pi\inv=\alpha\pi$. A path $\pi$ is \emph{closed} or a \emph{cycle} [at the vertex $v$] if $\alpha\pi=\omega\pi\ [= v]$. A path is \emph{simple} if no vertex occurs twice, except that possibly the initial and terminal vertices may coincide, in which case we have a simple closed path.  We also consider, for each vertex $v$, the \emph{empty path at $v$}, denoted $\varepsilon_v$ for which we set $\alpha\varepsilon_v=v=\omega \varepsilon_v$ and $\varepsilon_v\inv=\varepsilon_v$ (it is convenient to identify $\varepsilon_v$ with the vertex $v$ itself). We say that $\pi$ is a path \emph{from $u=\alpha \pi$ to $v=\omega\pi$}, and we  will also say that $u$ and $v$ are \emph{connected by $\pi$} (and likewise by $\pi\inv$). A path $\pi$ is \emph{reduced} if it does not contain any segment of the form $ee\inv$; successive deletion of such segments until there is no more such segment leads to the \emph{reduced form} $\mathrm{red}(\pi)$ of $\pi$, which is again a path and has the same initial and terminal vertices as the original path $\pi$.

A graph is \emph{connected} if any two vertices can be connected by some path.    A connected graph $\mT$ is a \emph{tree} if it does not admit non-empty reduced closed paths; a graph  is a \emph{forest} if each connected component is a tree. The subgraph $\langle \pi\rangle$ spanned by the non-empty path $\pi$ is the graph spanned by the edges of $\pi$; it coincides with $\langle \pi\inv\rangle$; the graph spanned by an empty path $\varepsilon_v$ simply is $\{v\}$ (one vertex, no edge). An \emph{arc} 
is a graph spanned by a simple non-closed path; a \emph{cycle} or a \emph{cycle graph} is a graph spanned by a simple closed path. The \emph{length} of a cycle or an arc is the number of its geometric edges. The \emph{cyclomatic number} of a finite  graph $\mK=(K\cup V;\alpha,\omega,\inv)$ is $|{K^\bullet}|-|V|+c$ where $c$ is the number of the connected components of $\mK$. This number is non-negative and is equal to $0$ if and only if $\mK$ is a forest.

Let $A$ be a finite set; a \emph{word labelling over the alphabet $A$} of the graph
$\mathcal{K}=(V\cup K;\alpha,\omega,{}^{-1})$ is a mapping $\ell\colon
K\to\til{A}^*$ respecting the involution: $\ell(e^{-1})=\ell(e)^{-1}$
for all $e\in K$.
If every label $\ell(e)$ is an element of $\til A$ we call it an \emph{$A$-labelling}. Every $A$-labelling 
$\ell\colon K\to \til{A}$ gives rise
to an \emph{orientation} of $\mK$: setting $E:=\{e\in K \colon \ell(e)\in A\}$
({positive} edges) and  $E^{-1}:=\{e\in K\colon \ell(e)\in A^{-1}\}$
({negative} edges),
it follows that $E\cap E^{-1}=\varnothing$ and we get $K=\til{E}$. For a path $e_1\cdots e_k$ in a word-labelled graph over $\til A$ as well as in an $A$-labelled graph $\mathcal{K}$, its label  is $\ell(e_1\cdots e_k):=\ell(e_1)\cdots\ell(e_k)$, which is a word over $\til{A}$. 

We consider $A$-labelled graphs as structures $\mK=(V\cup
K;\alpha,\omega,{}^{-1},\ell,A)$ in their own right. By a
\emph{subgraph of an $A$-labelled graph} we mean just a subgraph
with the induced labelling.
Morphisms $\varphi\colon \mK\to \mL$ of $A$-labelled graphs and congruences $\Theta$ on an $A$-labelled graph $\mK$ are defined accordingly. They must be morphisms, respectively congruences, on the underlying graphs and must in addition respect the labelling. For a morphism $\varphi\colon \mK\to\mL$ this means that $\ell(\varphi(e))=\ell(e)$ for every edge $e$ of $\mK$. For a congruence $\Theta$ on $\mK$ this means that $e\mathrel{\Theta}f\Longrightarrow \ell(e)=\ell(f)$ for all edges $e,f\in \mK$. This guarantees that the labelling $\ell$ is well-defined on the quotient graph $\mK/\Theta$. Again we have the usual homomorphism theorem.

An $A$-labelled graph is an
\emph{$A$-graph} if every vertex $u$
has, for every label $a\in \til{A}$, at most one edge
with initial vertex $u$ and label $a$. In the literature, such graphs occur under a variety of different names, such as \emph{folded graph}~\cite{KapMas} or \emph{inverse automaton}~\cite{AS,auinbors}, to mention just two. Every $A$-labelled graph  $\mK$ admits a unique smallest congruence $\Theta$ such that the quotient $\mK/\Theta$ is an $A$-graph.
In case $\mK$ is finite, this quotient can be computed nicely by the following procedure:  if there is a vertex from which there emerge two edges $e_1,e_2$ having the same label $a\in \til A$, then identify these edges and their terminal vertices: 
\begin{tikzpicture}[xscale=0.4,yscale=0.3,  baseline = {(0,-0.1)}]
	\draw(3.5,0)node{$\leadsto$};
	\draw[-latex](0.1,0.1)--(1.9,0.9);
	\draw[-latex](0.1,-0.1)--(1.9,-0.9);
	\filldraw(0,0)circle(2pt);
	\filldraw(2,1)circle(2pt);
	\filldraw(2,-1)circle(2pt);
	\draw (1,0.95)node{$a$};
	\draw(1,-0.95)node{$a$};
	()	\filldraw(5,0)circle(2pt);
	\filldraw(7,0)circle(2pt);
	\draw[-latex](5.15,0)--(6.85,0);
	\draw(6,0.45)node{$a$};
\end{tikzpicture} ,
and continue this process until no such pair of edges is left. This algorithm is called \emph{Stallings foldings} and goes back to Stallings~\cite{stallings}. It is useful to have in mind that two vertices $p_1$ and $p_2$ are identified by this process (that is, by the congruence $\Theta$) if and only if there is a path $p_1\longrightarrow p_2$ labelled by a Dyck word. It follows that in an $A$-graph, if the label $\ell(\pi)$ of a path $\pi$ is a Dyck word over $A$ then the path itself is a Dyck word over $E$ (with $E$ the set of positive edges of $\mK$, the orientation of $\mK$ being induced by the $A$-labelling). 
In an $A$-graph $\mathcal{K}$,
for every word $p\in \til{A}^*$ and every vertex $u$ there is at most
one path $\pi=\pi_u^\mathcal{K}(p)$ with initial vertex $\alpha\pi=u$ and
label $\ell(\pi)=p$. We set $u\cdot p:=\omega(\pi_u^\mathcal{K}(p))$, the terminal vertex of the path starting at $u$ and being labelled $p$ if such a path $p$ exists; then,  for every $p\in \til{A}^*$, the partial mapping 
$[p] \colon V\to V$, $u\mapsto u\cdot p$ is a partial bijection of the
vertex set $V$ of $\mathcal{K}$. 
Thus the involutory monoid $\til{A}^*$ acts on $V$ by partial bijections on the right.
The inverse monoid 
\begin{equation}\label{eq:transition inverse monoid}
	\mathrsfs{T}(\mathcal{K}):=\{[p]\colon p\in \til{A}^*\}
\end{equation}
of partial bijections obtained this way, is called the \emph{transition inverse monoid} $\mathrsfs{T}(\mathcal{K})$ of the $A$-graph $\mathcal{K}$. We note that the empty word $1\in \til A^*$ induces the identity mapping on all of $V$.
The transition inverse monoid $\mathrsfs{T}(\mathcal{K})$ is an $A$-generated
inverse monoid  in a natural way, the letter $a\in \til{A}$ induces the
partial bijection $[a]$ on the vertex set that maps every vertex $u$ to the terminal vertex
$\omega\pi_u(a)$ of the edge $\pi_u(a)$, which is the unique edge with
initial vertex $u$ and label $a$ if such an edge exists.
 
Similarly as for $A$-generated inverse monoids, morphisms of $A$-graphs are essentially unique: given two connected $A$-graphs $\mK$ and $\mL$, there is, for any pair of vertices $u\in \mathrm{V}(\mK)$ and $v\in \mathrm{V}(\mL)$,
at most one morphism $\mK\to \mL$ of $A$-graphs mapping $u$ to $v$. In case such a morphism exists and is surjective we address it as \emph{canonical morphism} $\mK\twoheadrightarrow \mL$.
 
 An $A$-graph $\mathcal{K}$ is called \emph{complete} or  a
\emph{group action graph} (also called  \emph{permutation
	automaton}) if every vertex $u$
has, for every label $a\in \til{A}$ exactly one edge $f$ with initial vertex $\alpha f=u$ and label $\ell(f)=a$. In this case, for every word $p\in \til{A}^*$ and every vertex $u$ there exists exactly one path $\pi=\pi_u^\mathcal{K}(p)$ starting at $u$ and having label $p$. Hence the transition inverse monoid $\cT(\mK)$ is a group of permutations of $V$.

\subsection{Cayley graphs of $A$-generated groups, Sch\"utzenberger graphs of $A$-generated inverse monoids}\label{subsec:caley}
\label{CayleyMMsubsec}
Given an $A$-generated group $Q$, the \emph{Cayley graph}
$\mQ$ of $Q$ has vertex set $Q$ and, for $g,h\in Q$ and $a\in \til A$ there is an edge $g\longrightarrow h$ labelled $a$ if any only if $g[a]_Q=h$. More formally we define the Cayley graph $\mQ$ of the $A$-generated group $Q$ as follows:
\begin{itemize}
	\item[--] the set of vertices of $\mQ$ is $Q$,
	\item[--] the set of edges of $\mQ$ is $Q\times \til{A}$,
	and, for $g\in Q,\ a\in \til{A}$, the incidence functions, involution
	and labelling are
	defined according to 
	\[
	\begin{array}{l@{\;\;:=\;\;}l}
		\alpha(g,a) & g,
		\\
		\hnt
		\omega(g,a) &g[a]_Q,
		\\
		\hnt
		(g,a)\inv \!\!&(g[a]_Q,a\inv),
		\\
		\hnt
		\ell(g,a) & a.
	\end{array}
	\]
\end{itemize}               
The edge $(g,a)$ should be thought of as
$\underset{g}{\bullet} \overset{a}{\rbl}
\!\underset{ga}{\bullet}$, its
inverse as
$\underset{g}{\bullet}\overset{a\inv}{\lbl}\!\underset{ga}{\bullet}$,
where $ga$ stands for $g[a]_Q$.  We note that $Q$ acts on $\mQ$ by
left multiplication as a group of automorphisms via
\[g\longmapsto {}^hg:=hg \quad \mbox{ and } \;\;
(g,a) \longmapsto {}^h(g,a) := (hg,a)
\]
for all $g,h\in Q$ and $(g,a)\in Q\times\til{A}$, where $h$, as a group element, acts on the vertex $g$ of $\mQ$  and on the edge $(g,a)$  of $\mQ$. Clearly, the Cayley graph $\mQ$ is not an invariant of the plain group $Q$ but depends on the assignment function $i_Q$. 
The group $Q$ itself is isomorphic to the automorphism group of its 
Cayley graph viewed as an $A$-graph.

One could equally well define the Cayley graph of an $A$-generated inverse monoid $M$. It has turned out however, that a collection of certain subgraphs of that graph --- the Sch\"utzenberger graphs --- is more important than the entire Cayley graph itself. Recall that Green's equivalence relation $\mR$ is defined on any monoid $M$ by setting, for $m,n\in M$, $m\mathrel{\mR}n$ if and only if $mM=nM$ (that is, $m$ and $n$ generate  the same right ideal of $M$). The $\mR$-class of an element $m$ is usually denoted $R_m$. For any inverse monoid $M$, we have that $m\mathrel{\mR}n$ if and only if $mm\inv =nn\inv$. Equivalently, $m\mathrel{\mR}n$ if and only if  $m=nx$ and $n=mx\inv$ for some $x$. In particular, every $\mR$-class $R_m$ of $M$ contains exactly one idempotent, namely $mm\inv$. {Sch\"utzenberger graphs} have been introduced by Stephen~\cite{stephen} and are defined as follows. For an $A$-generated inverse monoid $M$ and an element $m\in M$ the \emph{Sch\"utzenberger graph $\mM_m$ of $M$ with respect to $m$} has vertex set $R_m$ and for $p,r\in R_m$ and $a\in \til A$ there is an edge $p\longrightarrow r$ labelled $a$ if and only if $p[a]_M=r$. 
Important for us is the following result.
\begin{Lemma}[Theorem~3.2 in~\cite{stephen}]\label{lem:stephen} For every $A$-generated inverse monoid $M$ and any $m\in M$, the Sch\"utzenberger graph $\mM_m$ is an $A$-graph, and for any word $w\in \til A^*$ we have $[w]_M\ge m$ if and only if $w$ labels a path $mm\inv \longrightarrow m$ in $\mM_m$
	\end{Lemma}
\subsection{Extending partial bijections and automorphisms of relational structures}\label{subsec:extending graphs} We are concerned with a (usually) finite set $X$ and a set of partial bijections on $X$, usually denoted $A$, and we are interested in a superset $Y\supseteq X$ and, for every $a\in {A}$ a bijection $\wh{a}$ of $Y$ that extends $a$. Later we shall assume that $X$ and $Y$ carry  relational structures and the involved (partial) mappings are (partial) automorphisms of the structure in question. The mappings $a\in A$ and their extensions $\wh a$ (often also simply denoted $a$) are assumed to act on the right, written as $x\mapsto x\cdot a$ (partial mappings) or $x\mapsto xa$ (total mappings).

In the paper it will be necessary to consider certain algebraic structures generated by a set of symbols whose elements are in bijective correspondence with the elements of $A$; likewise we will need to consider graphs whose edges are labelled by such symbols. In order to avoid overly cumbersome notation we will  not distinguish notationally between such symbols and partial bijections on $X$ (just as in~\cite{HL}). Hence $A$ and $a\in A$ will denote both, depending on the context: 
a set of partial bijections on a set $X$ and at the same time a set of symbols. 
In both cases we denote by $\til A$ the set $A\cup A\inv$,  and we note that in case of abstract symbols we always have $A\cap A\inv=\varnothing$, while in case of partial bijections, an element $a$ may be self inverse: $a\in A\cap A\inv$.

\subsubsection{Embedding $A$-graphs into complete $A$-graphs and unary algebras} Suppose we are given a finite set $X$ and a set $A$ of partial bijections on $X$.
The data $(X,\til A)$ can be encoded in an $A$-graph $\mP$ in an obvious way: the vertices are the elements of $X$ and, 
for every pair $x_1,x_2\in X$, there is an edge $e\colon x_1\longrightarrow x_2$ with label $\ell(e)=a\in \til A$ if and only if $x_1\cdot a=x_2$. Note that $x_1\cdot a=x_2$ if and only if $x_2\cdot a\inv=x_1$. Connected with this is the transition inverse monoid $\cT(\mP)$ whose action will be denoted $x\mapsto {x}\cdot [w]_{\cT(\mP)}$ for the partial bijection $[w]_{\cT(\mP)}$ and $w\in \til A^*$. 
To have a superset $Y\supseteq X$ and bijections $\wh{a}$ on $Y$ that  extend the partial bijections $a\in \til A$ is the same as to have an embedding $\mP\hookrightarrow \ol{\mP}$ of the $A$-graph $\mP$ into a complete $A$-graph $\ol{\mP}$ on the vertex set $Y$. In this situation, the complete graph $\ol{\mP}$ is 
\emph{some completion} of $\mP$. 
While there are many ways to complete an incomplete $A$-graph $\mP$, there is an obviously unique maximal one, which is infinite.
In order to construct it, take the Cayley graph $\mF$ of the free $A$-generated group $F$, choose some letter $a\in \til A$ and remove from the graph all edges $e$ (and their inverses) with initial vertex $1_F$ whose label is not $a$. The graph $\mF$ splits into $2|A|$ subtrees of $\mF$, exactly one of which contains the vertex $1_F$. Call this tree the \emph{free $a$-cone} $[a\mF)$ with base vertex $1_F$, now re-denoted $1_a$. Next take the incomplete $A$-graph $\mP$ and attach a free $a$-cone $[a\mF)$ to every vertex $v\in \mathrm{V}(\mP)$ that has no edge $e$ with this initial vertex $\alpha e=v$ and label $\ell(e)=a$. Formally, this is achieved by forming the disjoint union $\mP\sqcup [a\mF)$ and factoring by the congruence identifying the vertices $v$ and $1_a$ and leaving all other vertices and edges unidentified. By doing so for all vertices $v\in \mathrm{V}(\mP)$ and all letters $a\in \til A$ we get a complete $A$-graph. We denote this  graph as $\mP\mathrm F$ and call it the \emph{free completion of} $\mP$. 

A more intuitive approach to get the free completion $\mP\mathrm{F}$ of an $A$-graph $\mP$ is as follows. For every vertex $v$ and every letter $a\in \til A$ for which there is no edge in $\mP$ with initial vertex $v$ and label $a$,  add a new vertex $va$, say, and a new edge $e\colon v\overset{a}{\longrightarrow} va$, and the inverse of this new edge. 
Every newly added vertex requires the addition of $2|A|-1$ further vertices;
indeed, only one edge emerges from $va$, namely $e\inv\colon va\overset{a\inv}{\longrightarrow} v$ carrying the label $a\inv$. Hence, for every $b\in \til A\setminus\{a\inv\}$ we have to further add a new vertex $vab$, say, and the edge $va\overset{b}{\longrightarrow} vab$ (and its inverse), and we have to continue to do so ad infinitum. Having done so, we have have attached a free $a$-cone $[a\mF)$ to the vertex $v$.  The result $\mP\mathrm{F}$ may thus be obtained as the union $\bigcup \mP_n$ of an ascending chain $\mP\subseteq \mP_1\subseteq\mP_2\subseteq  \cdots\subseteq \mP_n\subseteq \cdots$ of  extensions of $\mP$ where $\mP_n\supseteq \mP$ is comprised of the vertices of $\mP\mathrm{F}$ which have distance at most~$n$ 
from $\mP$ and the edges (and their inverses) required to access them from $\mP$.

The process of successively adding vertices $va_1, va_1a_2,\dots$ leads to the universal algebraic  view of the situation. The letters $a\in \til A$ may be viewed as \textsl{algebraic operations} successively applied to the elements $v$ of the given \textsl{generating set} $X$ (the vertices of $\mP$). To make this more precise, let $Y_\mP:=\mathrm{V}(\mP\mathrm{F})$ be the set of vertices of $\mP\mathrm{F}$; every letter $a\in \til A$ induces a bijective mapping $y\mapsto ya$ on $Y_\mP$. Hence $(Y_\mP;\til A)$ may be interpreted as a \emph{unary algebra} that is, an algebraic structure all of whose fundamental operations are unary functions.   Let us call  any unary algebra $(Z;\til A)$ whose fundamental operations are denoted by the elements of $\til A$ an \emph{$\til A$-set}. In the context of $\til A$-sets $(Z;\til A)$ we shall use the following notation: for a word $w=a_1\cdots a_n\in \til A ^*$ and $z\in Z$ set
\[zw:=((\dots((za_1)a_2)\dots )a_{n-1})a_n.\]
Words $w\in \til A$ of length at least $2$ do not belong to the algebraic signature but can be seen as \emph{term operations} on $\til A$-sets.
We mention two important features of the $\til A$-set $(Y_\mP;\til A)$ obtained from $\mP\mathrm{F}$ in this way: 
\begin{itemize}
\item[(i)] 
it is generated, as an $\til A$-set, by the subset $X=\mathrm{V}(\mP)$ of $Y_\mP$ (for every $y\in Y_\mP$ there exists $x\in X$ and $w\in \til A^*$ such that $y=xw$, that is, $Y_\mP$ is the closure of $X$ under the operations $\til A$ of $(Y_\mP;\til A)$), and 
\item[(ii)] 
it satisfies every identity (or law) of the form $yw=y$ for every Dyck word $w$. 
\end{itemize}
Let us call $\til A$-sets which satisfy (ii) \emph{\tilA{F}-sets}. 
For later use we note that every \tilA{F}-set $(Y;\til A)$ gives rise to a complete $A$-graph whose set of vertices is $Y$, and for $y_1,y_2\in Y$ and $a\in \til A$ there is an edge $y_1\longrightarrow y_2$ having label $a$ if an only if $y_1a=y_2$. We call this graph the \emph{$A$-graph associated with the \tilA{F}-set $(Y;\til A)$} or the \emph{functional graph} of $(Y;\til A)$.

The class of all \tilA{F}-sets forms a variety in the sense of universal algebra, hence this class admits free objects. 
A realisation (or model) of the free \tilA{F}-set generated by $X$ is given by the base set $X\times F$ endowed with  the operations $(x,g)a:=(x,g[a]_F)$ for every $x\in X$, $g\in F$ and $a\in \til A$, and $X$ is embedded in $X\times F$ via $x\mapsto (x,1_F)$. It will be convenient to denote the base set 
$X\times F$ of this free \tilA{F}-set simply by $XF$; for its elements we write $xg$ for $(x,g)$ (if $g\ne 1_F$) and just $x$ for $(x,1_F)$.
We see that the $X$-generated \tilA{F}-set $(Y_\mP;\til A)$, which arises from the free completion $\mP\mathrm{F}$ of the graph $\mP$, 
is not a free \tilA{F}-set (unless every $a\in \til A$ represents the empty transformation, or, in other words, the graph $\mP$ which encodes the data $(X,\til A)$ has no edges). Conversely, this means that the edges of the graph $\mP$, or, equivalently, the partial mappings $a\in \til A$, if non-empty, can be interpreted as \textsl{non-trivial relations} which provide a \textsl{presentation} of the $X$-generated \tilA{F}-set $(Y_\mP;\til A)$: 
\[(Y_\mP;\til A)\cong (XF;\til A)/\Theta(\mP)\]
where $\Theta(\mP)$ is the smallest $\til A$-set congruence on $(XF;\til A)$ containing all pairs $(x_1[a]_F,x_2)$ for $x_1,x_2\in X$ and $a\in \til A$ for which $x_1\overset{a}{\longrightarrow}x_2$ is an edge of $\mP$. Thus the \tilA{F}-set $(XF;\til A)/\Theta(\mP)$ is the \tilA{F}-set \emph{defined by the generating set $X$ subject to the defining relations $\mP$} and will be denoted $(X_\mP F;\til A)$. The $A$-graph associated with this freely presented \tilA{F}-set $(X_\mP F;\til A)$ is then just the free extension $\mP\mathrm{F}$ we started with.

In the following,  we will no longer distinguish notationally between the base set of an algebraic structure and the algebraic structure itself. This is usual when dealing with algebraic structures and we did do so for groups, inverse monoids, etc. 

We may generalise the concept of an \tilA{F}-set. Instead of $F$ take any $A$-generated group $G$; an \tilA{F}-set $(Y;\til A)$ (or just $Y$) is called a \tilA{G}-set if it satisfies all identities of the form $yw=y$ for all $w\in \til A^*$ for which $[w]_G=1_G$. Again the class of all \tilA{G}-sets forms a variety  of $\til A$-sets and hence admits free objects. Again the free \tilA{G}-set generated by $X$, denoted $XG$, can be realised as the base set $X\times G$ endowed with the operations $(x,g)a:=(x,g[a]_G)$ for all $x\in X$,  $g\in  G$ and $a\in \til A$. The elements of $XG$ will again be denoted as $xg$ etc.~(instead of $(x,g)$ etc.). And again the relations given by the partial mappings $\til A$ on $X$ or by its 
graph realisation $\mP$, respectively, provide non-trivial relations which give rise to the \tilA{G}-set $X_\mP G$ \emph{generated by $X$ 
subject to the defining relations} $\mP$ which is defined as\footnote{{This construction occurs in many papers, see e.g.~ \cite{Herwig, Hrushovski, bitterlichdiss}, but usually has not been revealed as a free presentation of a unary algebra.}}
\begin{equation}\label{eq:GAsubjectomP}
X_\mP G:=XG/\Theta(\mP)
\end{equation}
where now $\Theta(\mP)$ is the smallest $\til A$-set congruence on $XG$ containing all pairs $(x_1[a]_G,x_2)$ for all $x_1,x_2\in X$ and $a\in \til A$ for which $x_1\overset{a}{\longrightarrow}x_2$ is an edge in $\mP$.

In this construction it may happen that distinct elements $x_1,x_2\in
X$ are related by $\Theta(\mP)$. Obviously, two distinct elements
$x_1,x_2\in X$ will be related by $\Theta(\mP)$ if there is a path
$\pi\colon x_1\longrightarrow x_2$ in $\mP$ whose label $\ell(\pi)$
has $G$-value $1_G$. That it happens only in this case follows from
the next lemma which may be seen as a solution to the word problem of
the presentation of the algebraic structure $X_\mP G$.

\begin{Lemma}\label{lem:extending G}
	Let $\mP$ be an $A$-graph with vertex set $X$ and let $G$ be an $A$-generated group. Then  $xg\mathrel{\Theta(\mP)} yh$ holds for $xg,yh\in XG$ if and only if there exists a word $w\in \til A^*$ such that ${x}\cdot[w]_{\cT(\mP)}=y$  and $gh\inv=[w]_G$.
\end{Lemma}
\begin{proof}
	Consider the binary relation $\Upsilon$ on $XG$ defined by
	\[xg\mathrel{\Upsilon}yh:\Longleftrightarrow \exists w\in \til A^*\colon {x}\cdot[w]_{\cT(\mP)}=y\mbox{ and }gh\inv=[w]_G.\]
	It is routine to verify that (i) $\Upsilon\subseteq \Theta(\mP)$, (ii) $\Upsilon$ is an equivalence relation and (iii) $\Upsilon$ is compatible with the $\til A$-set structure, that is, it is a congruence. Since $\Theta(\mP)$ is the smallest congruence which contains all pairs $(x[w]_G,y)$ for which the equality $x\cdot [w]_{\cT(\mP)}=y$ holds, it follows that $\Upsilon=\Theta(\mP)$.
\end{proof}

Of particular importance is when no two distinct elements $x_1, x_2\in X$ are related under $\Theta(\mP)$; that is,  $\Theta(\mP)\restriction X = \mathrm{id}_X$. For later use we state this  concept as a definition.

\begin{Def}\label{def: G extends mP} An $A$-generated group $G$ \emph{extends} an $A$-graph $\mP$ if \[\Theta(\mP)\upharpoonright \mathrm{V}(\mP)=\mathrm{id}_{\mathrm{V}(\mP)}.\] 
\end{Def} 
The following, related concept will also be used throughout the paper.
\begin{Def}\label{def: G extends M}	
An $A$-generated group $G$ \emph{extends} an $A$-generated inverse monoid $M$ 
if, for all words $w\in \til A$, $[w]_G=1_G$ implies that $[w]_M$ is idempotent.
\end{Def}

In semigroup-theoretic terms, an $A$-generated group $G$ extends an $A$-generated inverse monoid $M$ in the sense of Definition~\ref{def: G extends M} if and only if the $A$-generated product $M\times G$ is an $E$-unitary cover of $M$. It is clear that $G$ extends the $A$-graph $\mP$ in the sense of Definition~\ref{def: G extends mP} if and only if $G$ extends the transition inverse monoid monoid $\cT(\mP)$ in the sense of Definition~\ref{def: G extends M}. 

The value $[w]_M$ of any Dyck word $w$  in 
an $A$-generated inverse monoid $M$ is an idempotent --- this is easily 
seen by induction on the length of $w$ and taking into account that 
$x\inv ex$ is an idempotent for every element $x$ and every 
idempotent $e$. 
Therefore, the free group $F$ extends every $A$-generated inverse monoid $M$. The \textsl{existence} of a \textsl{finite} $A$-generated group $G$ extending a given finite $A$-generated inverse monoid $M$ is less obvious, but still well-known and easy to see: take any representation of $M$ as a monoid of partial bijections on some finite set $X$, take any extensions to total bijections on $X$ (or on some  finite superset $Y$ of $X$) of the partial bijections represented by the elements of $A$  and consider the group $G$ of permutations generated by these bijections --- then any such group $G$ of permutations of $X$ (or of $Y$) extends $M$ in the sense of Definition~\ref{def: G extends M}. In particular, for a given finite $A$-generated inverse monoid $M$ there is no canonical choice for a finite $A$-generated group $G$ to extend $M$ --- but there exist 
\textsl{definable} choices for given $M$.
By this we mean that we could single out one of several concrete specifications
that yield an unambiguous, explicitly definable output $G(M,A;\iota)$ for any given $M$, depending just on $M$ 
and the assignment of generators $\iota \colon a \mapsto [a]_M$ 
in an isomorphism-invariant manner even up to permutations of the set $A$
(in the sense that $(M,A;\iota)  \cong (M',A;\iota')$ 
implies $G(M,A;\iota) \cong G(M',A;\iota')$). 
One obvious such specification uses the above recipe based on the unique natural representation of $M$ as a monoid of partial bijections on the set $M$ itself, and the direct product of all the extensions to groups that arise from different combinations of extensions of the partial to global bijections of $X := M$.
(This is but one manner of replacing all arbitrary choices in the general recipe by
explicitly definable choices and/or parallelisation of all possible choices.) 
The potentially interesting advantage
of such a concrete and explicit specification with definable outcome precisely lies 
in its isomorphism-invariant nature, which in particular guarantees compatibility with internal symmetries. 
If, for instance, $[a]_M$ and $[a']_M$ are 
related by an automorphism of $(M,A;\iota)$, 
then that translates, just by isomorphism-invariance, 
into the existence of an automorphism of such $G$  
that correspondingly relates $[a]_G$ and $[a']_G$. This fundamental feature,
which is treated in some detail in~\cite{otto4}, relies on explicit definability,
while `free choices' as in the general recipe above can break
internal symmetries (so that, e.g., $[a]_G$ and $[a']_G$ could have
different order in $G$).

Also note that $G$ extends $\mP$ in the sense of Definition~\ref{def: G extends mP} if and only if the natural mapping $X\to X_\mP G$, $x\mapsto x\Theta(\mP)$ provides an embedding of the $A$-graph $\mP$ into the complete $A$-graph associated with the $\til A$-set $X_\mP G$ mentioned earlier.
An immediate consequence of Lemma~\ref{lem:extending G} is the following  well-known fact.
\begin{Lemma}\label{lem:extending P}
	Let $\mP$ be an $\til A$-graph on vertex set $X$ and let $G$ be an $A$-generated group; then $G$ extends $\mP$ if and only if for each word $w\in \til A^*$ for which $[w]_G=1_G$ the partial transformation $[w]_{\cT(\mP)}$ on $X$ is a partial identity, that is, $x\cdot [w]_{\cT(\mP)}=x$ for all $x\in X$ for which it is defined. This is equivalent to the fact that $[w]_{\cT(\mP)}$ is an idempotent in the transition inverse monoid $\cT(\mP)$.
\end{Lemma}

\subsubsection{Relational structures}\label{subsec:relational structures}  We now are ready for the next step. 
We think of $X$ as the base set of a relational structure, i.e.~as being equipped with interpretations of the relation 
symbols in a fixed finite relational signature $\sigma$,
and of the given partial bijections $\til A$ of $X$ as
\textsl{partial automorphisms} of that structure. 
So, let us fix a finite relational signature $\sigma$ consisting of symbols $R$, say, such that each such $R$ has an \textsl{interpretation} as a relation $R^\bX$ on $X$, of some arity $r$. That is, every symbol $R$ comes with a well-defined subset $R^\bX\subseteq X^r$ ($r$ the arity of $R$); the structure $(X;(R^\bX)_{R\in \sigma})$ will be denoted $\bX$. 
For $Z\subseteq X$, a $\sigma$-structure $\mathfrak{Z}=(Z;(R^\mathfrak{Z})_{R\in \sigma})$ is a \emph{weak substructure} of $\mathfrak{X}$ if $R^\mathfrak{Z}\subseteq R^\mathfrak{X}$ for all $R\in \sigma$,  while $\mathfrak{Z}$ is an \emph{induced substructure} of $\mathfrak{X}$  if $R^\mathfrak{Z}=R^\mathfrak{X}\cap Z^r$ for all $R\in \sigma$ where $r$ denotes the arity of $R$. A $\sigma$-structure $\mathfrak{Y}$ is an \emph{extension} of $\mathfrak{X}$ if $\mathfrak{X}$ is an induced substructure of $\mathfrak{Y}$.

We are interested in  extensions $\bY$ of $\bX$ such that the partial automorphisms of $\bX$ extend to total automorphisms of $\bY$.
In this context it is convenient to introduce some further terminology. Suppose we are given a $\sigma$-structure $\bX=(X;(R^\bX)_{R\in \sigma})$ and a set $\til A$ of partial permutations of $X$ that, in fact, are partial automorphisms of $\bX$. The pair $(\bX,\til A)$  is called an \emph{extension problem}. A pair $(\bY,\{\wh{a}\colon a\in \til A\})$ 
is a \emph{solution} for the extension problem $(\bX,\til A)$ if $\bY$ is an extension of $\bX$, $\wh a$ is an automorphism of $\bY$, $a\subseteq \wh a$ and $\wh{a}\inv= \wh{a\inv}$ 
for every $a\in \til A$. If the automorphisms $\wh a$ are clear from context, 
we occasionally refer to the structure $\bY$ itself as a solution for 
the extension problem $(\mathfrak{X}, \til A)$.

We shall apply the concepts of the preceding subsection. Let $G$ be a not necessarily finite $A$-generated group extending the graph $\mP$ given by the data $(X,\til A)$ according to Definition~\ref{def: G extends mP}. 
In order to equip all of $Y:=X_\mP G$ with a relational structure such that every $g\in G$ is not only a bijection on $Y$ but an automorphism of $\bY$ it is natural to add to the interpretations of the various symbols $R$ on $\bX$ whatever is required to make $g$ an automorphism on all of $\bY$. For every $R\in \sigma$ set
\[R^\bY:=\bigcup\{R^\bX g\colon g\in G\}\]
where $R^\bX g=\{(x_1g,\dots,x_rg)\colon (x_1,\dots,x_r)\in R^\bX\}$.
In this way, $\bY=:\bX_\mP G$ becomes a $\sigma$-structure and all $g\in G$ act as automorphisms of $\bY$. The structure $\bX$ is clearly a weak substructure of $\bY$. We give a sufficient condition on $G$ for $\bX$ to be an induced substructure of $\bY$ , that is, for $\bY$ being a solution for the 
extension problem $(\bX,\til A)$.

\begin{Lemma}\label{lem:X is induced} Let $\bX$ be a finite $\sigma$-structure on base set $X$ and $\til A$  a set of partial automorphisms of $\bX$, encoded in the $A$-graph $\mP$ on the vertex set $X$. 
Let $G$ be an $A$-generated group which extends $\mP$ and let  $M=\cT(\mP)$.
Suppose that $G$ satisfies the following condition:
for each $s\le |X|$ and each collection of words $w_1,\dots,w_s\in \til A^*$ for which  $[w_1]_G=\cdots=[w_s]_G$ there exists a word $w\in \til A^*$  such that $[w_1]_M,\dots,[w_s]_M\le [w]_M$. Then $\bX$ is an induced substructure of $\bX_\mP G$.  
\end{Lemma}
\begin{proof}
Let $x_1,\dots, x_r\in X$ and suppose that $(x_1,\dots, x_r)\in R^{\bX_\mP G}$ for some $R\in \sigma$. Then there are $y_1,\dots, y_r\in X$ and $g\in G$ such that $(y_1,\dots,y_r)\in R^\bX$ and $x_i=y_ig$ for $i=1,\dots, r$. By Lemma~\ref{lem:extending G} there exist words $w_1,\dots, w_r\in \til A^*$ such that $y_i\cdot[w_i]_M=x_i$ and $[w_i]_G=g$. 
If $r>|X|$ then the $r$-tuple $(y_1,\dots, y_r)$ nevertheless contains at most $|X|$ distinct elements. In any case, we may assume that the collection of words $w_1,\dots, w_r$ contains at most $|X|$ distinct words. By the assumption of the Lemma there exists a word $w\in \til A^*$ for which $[w_i]_M\le [w]_M$ for all $i$. This implies $y_i\cdot[w]_M=x_i$ for all $i$. Since $[w]_M$ is a partial automorphism of $\bX$ it follows that $(x_1,\dots,x_s)\in R^\bX$. 
\end{proof} 
\begin{Rmk}\label{rmk:width vs |X|}
	If in Lemma~\ref{lem:X is induced} the largest arity $r$ of $\sigma$ is smaller than the cardinality $|X|$ then the bound $|X|$ may be replaced by~$r$.
\end{Rmk}
Readers familiar with semigroup theory will observe that the condition on $G$, in relation to the inverse monoid $M$, formulated in Lemma~\ref{lem:X is induced} just says that $G$  witnesses the group-pointlike sets of $M$ of order up to $|X|$.

\subsubsection{Classes defined by forbidden homomorphisms and free extensions of relational structures} 
Let $\sigma$ be a finite relational signature and $\Sigma$  a finite set of finite $\sigma$-structures; let  $\mathbf{Excl}(\Sigma)$ be the class of all $\sigma$-structures $\bS$ which do not admit a homomorphism $\bT\to \bS$ for 
any $\bT\in \Sigma$. A structure $\bS$ which does not admit such a homomorphism is said to be \emph{$\Sigma$-free} or to \emph{avoid} (all members of) $\Sigma$. The members of $\Sigma$ are called the structures \emph{forbidden for the class} $\mathbf{Excl}(\Sigma)$. A class $\mathbf{C}$ of $\sigma$-structures is said to have the \emph{extension property for partial automorphisms} (EPPA for short) if 
every extension problem $(\mathfrak{S}, \til A)$ for finite $\mathfrak{S} \in \mathbf{C}$ that admits a solution in $\mathbf{C}$ also admits a finite solution in $\mathbf{C}$. 
Given a finite $\sigma$-structure $\bX=(X;(R^\bX)_{R\in \sigma})$ and a set of partial automorphisms $\til A$ of $\bX$ encoded in the finite $A$-graph $\mP$, we can extend $\bX$ by means of the free $A$-generated group $F$, as explained in  Section~\ref{subsec:relational structures} ($F$ satisfies the condition formulated in Lemma~\ref{lem:X is induced}). The $\sigma$-structure on $X_\mP F$ so obtained is called 
the \emph{free extension} $\bX_\mP F$ (\emph{subject to }$\mP$) of
$\bX$. Suppose that $(\bX,\til A)$ has some solution $(\bY,\{\hat a\colon a\in \til A\})$; then $G:=\langle \wh{a}\colon a\in A\rangle$, is an $A$-generated subgroup of the automorphism group $\mathrm{Aut}(\bY)$ of $\bY$. If $\bY$ is $\bT$-free for some finite $\sigma$-structure $\bT$ then so is $[\bX]G$, the substructure of $\bY$ induced on the orbit of $\bX$ under the action of $G$. It follows from principles of universal algebra that the identity mapping $X\to X$, $x\mapsto x$ extends to surjective homomorphisms of $\til A$-sets $X_\mP F\twoheadrightarrow X_\mP G\twoheadrightarrow [X]G$ (the latter being the base set of $[\bX] G$). From the definition of the $\sigma$-structures on $\bX_\mP F$ and $\bX_\mP G$, respectively, it follows that we even have surjective homomorphisms of the relational structures $\bX_\mP F\twoheadrightarrow \bX_\mP G\twoheadrightarrow [\bX] G$. 
Summing up, this means that, if $\bX$ has \textsl{some} $\bT$-free extension then the free extension $\bX_\mP F$ is also $\bT$-free.

\subsection{Stallings graphs and Schreier graphs}
\label{subsec:stallings} We introduce the Stallings and Schreier graphs assigned to a finitely generated subgroup $H$ of the free $A$-generated group $F$; we shall also consider the appropriate concepts for cosets $Hg$. 

Let $H$ be a finitely generated subgroup of $F$ and let $w_1,\dots, w_n$ be generating elements, realised as reduced words over $\til A$. 
Consider the disjoint union of cycle graphs $\mC_1,\ldots,\mC_n$, 
where $\mC_i$ is the graph spanned by 
a closed simple path of length~$|w_i|$  
labelled~$w_i$ at a base vertex $\iota_i$. 
Next factor this disjoint union of labelled cycle graphs $\mC_i$ by the smallest $A$-graph congruence that identifies all base vertices $\iota_i$. 
The equivalence class containing all vertices $\iota_i$ then is a distinguished vertex $\iota$ of the quotient graph. 
This pointed graph with distinguished vertex $\iota$ is called the \emph{Stallings graph $\mS_H$ of the group $H$}.
It depends just on $H$, not on the choice of the generating elements $w_1,\dots, w_n$ (see Proposition~\ref{prop:Stallingsgraph}). 
It can be obtained in a nice manner from 
the cycle graphs $\mC_i$ by a procedure as follows: first identify all base vertices $\iota_1,\dots,\iota_n$ with each other and denote the resulting vertex as $\iota$ and then apply Stallings foldings until the result is an $A$-graph; for more details the reader may consult~\cite{KapMas}.

The essential feature of the Stallings graph $\mS_H$ is this: the group $H$ is comprised exactly of the (reduced) words 
labelling reduced paths in $\mS_H$ that are closed at the vertex $\iota$. 
The Stallings graph $\mS_H$ is an invariant of the group $H$ in the strong sense of Proposition~\ref{prop:Stallingsgraph} below. In order to formulate it we need another concept: for a finite connected $A$-graph $\mP$ with distinguished vertex $\iota$ let $L(\iota,\mP)$ denote the set of all words over $\til A$ that are labels of reduced paths in $\mP$ that are closed at $\iota$. All words in $L(\iota,\mP)$ are reduced, 
and it is easy to see that $L(\iota,\mP)$ is a finitely generated subgroup of $F$. By construction, $L(\iota,\mS_H)=H$. The uniqueness of the Stallings graph $\mS_H$ is characterised by the following proposition.
\begin{Prop}[\cite{KapMas}, Theorem 5.2]\label{prop:Stallingsgraph} Let $H$ be a finitely generated subgroup of $F$ and $\mP$ a finite connected $A$-graph with distinguished vertex $\iota$ such that every vertex of $\mP$ except possibly $\iota$ has degree at least $2$. If $L(\iota,\mP)=H$ then there is a unique isomorphism of pointed graphs $\mS_H\twoheadrightarrow \mP$.
\end{Prop}

We can generalise this concept to cosets $Hg$, already discussed in Section~4.17 in the monograph~\cite{q-theory}. Let $H$ be given, $\mS_H$ be its Stallings graph and let $g\in F$ be represented as a reduced word $g=a_1\cdots a_n\in \til A^*$. Take the graph $\mL$ spanned by a path $p_0\underset{e_1}{\longrightarrow} p_1\underset{e_2}{\longrightarrow}\cdots \underset{e_n}{\longrightarrow} p_n$ where $n=|g|$ is the length of the word $g$. 
We label the edge $e_i$ by the letter $a_i$ so that the label of the path $e_1\cdots e_n$ is just $g$ (with inverse edges labelled accordingly). 
We form the disjoint union of $\mS_H$ and $\mL$ and factor by the smallest $A$-graph congruence that identifies the distinguished vertex $\iota$ of $\mS_H$ and the initial vertex $p_0$ of $\mL$. In other words, we identify the  vertex $p_0$ with the  vertex $\iota$ of $\mS_H$ and perform Stallings foldings until the result is an $A$-graph. 
Denoting the vertex corresponding to the terminal vertex $p_n$ of the path 
in this graph by $\tau$, we obtain a resulting, $2$-pointed graph $\mS_{Hg}$, with $\iota$ and $\tau$ as distinguished vertices. 
We call $\iota$ the \emph{initial} vertex  and $\tau$ the \emph{terminal} vertex of $\mS_{Hg}$, and refer to the $2$-pointed graph $\mS_{Hg}$ as the \emph{Stallings graph of the coset} $Hg$. This term is justified since $Hg$ is comprised of all words that are labels of reduced paths running in $\mS_{Hg}$ with initial vertex $\iota$ and terminal vertex $\tau$. It may happen that $\iota=\tau$; this is exactly the case when $H=Hg$, that is, when $g\in H$. 
For Stallings graphs of cosets a characterisation analogous to Proposition~\ref{prop:Stallingsgraph} holds. For a finite connected $A$-graph $\mP$ with two distinguished vertices $\iota$ and $\tau$ let $L(\iota,\mP,\tau)$ be the set of all words that are labels of reduced paths $\pi\colon \iota\longrightarrow \tau$ running in $\mP$. Every word of $L(\iota,\mP,\tau)$ is reduced and $L(\iota,\mP,\tau)$ is a coset $Hg$ with $H$ finitely generated and $g$ a word labelling a path $\iota\longrightarrow\tau$ in $\mP$. The proof of the next proposition is completely analogous to the one of Proposition~\ref{prop:Stallingsgraph}.
\begin{Prop}\label{prop:StallingsCoset} Let $H$ be a finitely generated subgroup of $F$ and $g\in F$; let $\mP$ be a finite connected $A$-graph with distinguished vertices $\iota$ and $\tau$ such that all vertices of $\mP$ except possibly $\iota$ and $\tau$ have degree at least $2$. If $L(\iota,\mP,\tau)=Hg$ then there is a unique isomorphism of $2$-pointed graphs $\mS_{Hg}\twoheadrightarrow \mP$.
\end{Prop}

Moreover, in a Stallings graph $\mS_{Hg}$ with $\iota\ne \tau$ we may swap the r\^oles of $\iota$ and $\tau$, that is, declare $\iota$ to be the terminal vertex and $\tau$ to be the initial vertex. By doing so we get the Stallings graph of the coset $({}^{g\inv}\!\!H)g\inv$  where ${}^{g\inv}\!\!H$ is the conjugate group $g\inv Hg$ (note that the right coset $({}^{g\inv}\!\!H)g\inv$ coincides with the left coset $g\inv H$). Indeed, the labels of reduced paths in $\mS_{Hg}$ starting at $\tau$ and ending at $\iota$ are exactly the inverses of the labels of the reduced paths starting at $\iota$ and ending at $\tau$.

Stallings graphs $\mS_{Hg}$ of cosets $Hg$ will play a crucial r\^ole in the proofs of the implications  Herwig--Lascar $\Longrightarrow$ Ribes--Zalesskii and Herwig--Lascar $\Longrightarrow$ Ash. 
The free completion $\mS_{Hg}\mathrm{F}$ of $\mS_{Hg}$, 
which we call the \emph{Schreier graph} of the coset $Hg$, 
will also be important in this context.
We shall denote it as $\bm{\mS}_{Hg}$ and consider it also as a $2$-pointed graph with initial vertex $\iota$ and terminal vertex $\tau$. The behaviour of $\bm\mS_{Hg }$, compared to that of $\mS_{Hg}$, can be summarised as follows: for a \textsl{reduced} word $w\in \til A^*$, the equality $\iota\cdot w=\tau$ holds in $\mS_{Hg}$ if and only if $w\in Hg$, while for \textsl{any} word $w\in \til A^*$ the equality $\iota w=\tau$ holds in $\bm\mS_{Hg}$ if and only if $\mathrm{red}(w)\in Hg$. In contrast to the Schreier graph $\bm\mS_{Hg}$, the Stallings graph $\mS_{Hg}$ is not complete (unless $H$ is of finite index~\cite[Proposition 8.3]{KapMas}). It may happen that for a word $w$ for which $\mathrm{red}(w)$ belongs to $Hg$, $\iota\cdot w$ is not defined in $\mS_{Hg}$.

Usually, the Schreier graph of a subgroup $H$ of $F$ is defined to be the $A$-graph with vertex set $\{Hg\colon g\in F\}$, the set of all right cosets with respect to $H$, and edges $e\colon Hg\overset{a}{\longrightarrow} Hg[a]_F$, labelled $a$,  for all $g\in F$ and $a\in \til A$. It is therefore sometimes called the \emph{coset graph} or \emph{Schreier coset graph} with respect to $H$. It is known from~\cite{KapMas} that $\mS_H\mathrm{F}$ and the latter graph are isomorphic as pointed $A$-graphs, via a unique isomorphism mapping $\iota\mapsto H$, see~\cite{KapMas}. The Stallings graph $\mS_H$ is the subgraph of the Schreier graph $\bm\mS_H$ spanned by all reduced paths closed at the distinguished vertex $\iota$. The latter could be used as an alternative definition of the Stallings graph based on the alternative (usual) definition of the Schreier graph as coset graph with respect to $H$. For our purposes the view of the Schreier graph as $\bm{\mS}_H=\mS_H\mathrm{F}$ seems to be more appropriate.

We close this section by an easy and well-known, but important 
observation. Let $M$ be a finite $A$-generated inverse monoid; for $m\in M$ we let
\begin{equation}\label{eq:canonrelmor}
\varphi(m):=\bigl\{[w]_F\colon  w\in \til A^*, 
[w]_M=m \bigr\}. 
\end{equation}
Obviously, $\varphi$ is a mapping from $M$ to the powerset $2^F$ of $F$, 
which, in semigroup-theoretic terms, is 
the \emph{canonical relational morphism} $M\to F$.
\begin{Prop}\label{prop:canonicalrelmorMto F}
	Let $M$ be a finite $A$-generated inverse monoid and $m\in M$; then $\varphi(m)=Hg$ for some finitely generated subgroup $H$ of $F$ and some $g\in F$.
\end{Prop}
\begin{proof}
	Let $m\in M$ and recall from Section~\ref{subsec:caley} the Sch\"utzenberger graph $\mM_m$  with respect to $m$. Its set of vertices is $R_m$, the $\mR$-class of $m$, and for $p,r\in R_m$ and $a\in \til A$ there is an edge $p\longrightarrow r$ labelled $a$ if and only if $r=p[a]_M$. Let $H$ be the set of all reduced words labelling (necessarily reduced) paths in $\mM_m$ that are closed at the vertex $mm\inv$. 
    Since $\mM_m$ is finite, $H$ is a finitely generated subgroup of $F$ and, by Proposition~\ref{prop:Stallingsgraph}, $\mS_H$ is a subgraph of $\mM_m$ with distinguished vertex $\iota=mm\inv$. Moreover, $\mS_{H}$ contains every edge of $\mM_m$ that is contained in a non-trivial cycle subgraph of $\mM_m$. In other words, $\mM_m$ 
    possibly differs from $\mS_{H}$ just by finite trees attached to $\mS_H$. Let $g$ be  a reduced word labelling a (necessarily reduced) path $mm\inv \longrightarrow m$. Then by Proposition~\ref{prop:StallingsCoset}, the Stallings graph $\mS_{Hg}$ of the coset $Hg$ is also a subgraph of $\mM_m$ with distinguished vertices $\iota=mm\inv$ and $\tau=m$, respectively. Again, $\mM_m$ and $\mS_{Hg}$ differ at most by finite trees 
    attached to $\mS_{Hg}$. Let $u\in \til A^*$ be such that $[u]_F\in \varphi(m)$, that is, $[u]_M=m$; then $u$ labels a path $mm\inv \longrightarrow m$ in $\mM_m$ because $m=mm\inv m=mm\inv[u]_M$. The  reduced form $\mathrm{red}({u})$ still labels a path $mm\inv \longrightarrow m$, but this path now runs in the subgraph $\mS_{Hg}$ (by the aforementioned difference between $\mM_m$ and $\mS_{Hg}$), hence $[u]_F=\mathrm{red}({u})\in Hg$. Let conversely $f\in Hg$ be represented as a reduced word over $\til A$. Then $f$ labels a path $mm\inv \longrightarrow m$ in $\mS_{Hg}$ and hence also in $\mM_m$. According to Lemma~\ref{lem:stephen}, $m\le [f]_M$. Let again $u$ be a word with $[u]_M=m$. Then $[uu\inv f]_M=m$ and $[uu\inv f]_F=f$, hence $f\in \varphi(m)$ as required.
\end{proof}

\section{Theorems Revisited} \label{sec: revisited}
We formulate the three theorems in question.
\begin{Thm}[Ribes--Zalesskii~\cite{RibesZalesskii}]\label{thm: Ribes--Zalesskii}  For every  positive integer $n$, the  product $H_1\cdots H_n$ of every $n$-tuple of finitely generated subgroups $H_1, \dots, H_n$ of the $A$-generated free group $F$  is closed in the profinite topology of $F$.
\end{Thm}
Recall that a basis for the open sets of the profinite topology of $F$ is given by the set of all (left or right) cosets $Hg$ with $g\in F$ and $H$ a finite index subgroup of $F$. In fact, it suffices to consider cosets $Ng$ with $N$ a normal subgroup of finite index. From this it follows that a subset $W\subseteq F$ is closed if and only if for every element $g\notin W$ there exists a finite index normal subgroup $N$ of $F$ such that $Ng\cap W=\varnothing$. 
\begin{Thm}[Herwig--Lascar~\cite{HL}]\label{thm: Herwig--Lascar} Let $\sigma$ be a finite relational signature;  every class $\mathbf{C}$ of  $\sigma$-structures that is defined in terms of forbidden homomorphisms of finitely many finite $\sigma$-structures enjoys the extension property for partial automorphisms.
\end{Thm}
In order to formulate Ash's Theorem we need some further terminology. Given a graph $\mE$, an $A$-generated group $G$ and an $A$-generated inverse monoid $M$. We say that a word labelling $v\colon \mathrm{K}(\mE)\to \til A^*$  \emph{commutes} over the group 
$G$ (is $G$-commuting, for short) if, for every closed path $\pi$ in $\mE$ the $G$-value of the label of $\pi$ is the identity element:  $[\ell(\pi)]_G=1_G$; we say two labellings $u,v\colon \mathrm{K}(\mE)\to \til A^*$ are $M$-related if $[u(e)]_M=[v(e)]_M$ for every edge $e\in \mathrm{K}(\mE)$. The following is essentially the same as Theorem~5.1 in~\cite{Ash}, but slightly reformulated. However, this is what was literally proved in Sections~5 -- 7 in~\cite{Ash}. The original, different formulation of Ash's Theorem 5.1 came from the originally intended purpose of this theorem.
\begin{Thm}[Ash~\cite{Ash}, version for inverse monoids]\label{thm: ash} For every finite $A$-generated inverse monoid $M$ and every finite 
graph $\mE$ there exists a finite $A$-generated group $G$ such that every $G$-commuting word labelling of $\mE$ admits an $M$-related relabelling that  commutes over the free group $F$.
\end{Thm}

If a word labelling of $\mE$ does not admit an $M$-related and $F$-commuting relabelling, then, according to the theorem, this labelling is not $G$-commuting. Hence $G$ may be interpreted
as a finite device for eliminating
labellings that are ``bad'' with respect to 
having $M$-related relabellings that commute over $F$.

\subsection{Herwig--Lascar implies Ribes--Zalesskii}\label{subsec:HL>RZ}  Let $n\ge 2$; we consider relational structures over the signature $\sigma_n=\{R_1,\dots, R_n\}$ consisting of $n$ binary relational symbols. For elements $y_1,y_2$ in such a structure $\mathfrak{Y}=(Y;R_1^\mathfrak{Y},\dots,R^\mathfrak{Y}_n)$ in question we shall often write $y_1\mathrel{R_i}y_2$ instead of $(y_1,y_2)\in R_i^\mathfrak{Y}$ (in particular, we often shall omit the superscript ${}\mathfrak{Y}$ in $R_i^\mathfrak{Y}$ whenever there is no danger of confusion). Such structures can be seen as \textsl{edge-coloured directed graphs}, though they are not graphs in the sense of Section~\ref{subsec:graphs} and we do not use the term graph in this context in order to avoid confusion.

We consider the class $\mathbf{CCF}_n := \mathbf{Excl}(\bC_n)$
of all $\sigma_n$-structures $\mathfrak{C}$
that admit no homomorphism $\bC_n\to \bC$ from the $\sigma_n$-structure 
$\bC_n$ consisting of 
 $n$ distinct elements $c_1,\dots,c_n$ with
\[R_i^{\mathfrak{C}_n}=\{(c_i,c_{i+1})\}\mbox{ for }i=1,\dots, n-1\mbox{ and }R_n^{\mathfrak{C}_n}=\{(c_n,c_1)\}.\]
(The structure $\mathfrak{C}_n$ may be seen as the directed cycle graph of 
length $n$ whose edges have colours $1,\dots,n$.)
According to Theorem~\ref{thm: Herwig--Lascar}, $\mathbf{CCF}_n$ has the \emph{extension property for partial automorphisms}: 
 every extension problem $(\mathfrak{S}, \til A)$ 
 for finite $\mathfrak{S}$ in $\mathbf{CCF}_n$ 
 that admits a solution in $\mathbf{CCF}_n$
 also has a finite solution in $\mathbf{CCF}_n$.

The following proof of the implication Herwig--Lascar $\Longrightarrow$ Ribes--Zalesskii is perhaps not the shortest one, but it tries to prepare the ground for a direct proof of the implication Herwig--Lascar $\Longrightarrow$ Ash, which will be the content of Section~\ref{subsec:HL===>Ash}. The proof is inspired by the proof of the Ribes--Zalesskii Theorem given in~\cite{HL} but, in a sense, embeds those ideas of~\cite{HL} into the environment of Stallings and Schreier graphs (which seem to be \textsl{the} natural framework for this topic). Crucial will be the following result.
\begin{Prop}\label{prop:HL>RZ} Let $n\ge 2$, let $H_1,\dots, H_n$ be finitely generated subgroups of $F$ and $g_1,\dots,g_n\in F$.  Then there exists a finite index normal subgroup $N$ of $F$  such that: if $1_F\notin H_1g_1\cdots H_ng_n$ then $1_F\notin NH_1g_1\cdots NH_ng_n$.
	\end{Prop}
\begin{proof}
	We derive this result from the extension property for partial 
    automorphisms for the class $\mathbf{CCF}_n$.
    So, let $H_1,\dots, H_n$ be finitely generated subgroups of $F$ and $g_1,\dots, g_n\in F$ and, for each $i$, let $\mS_i\subseteq \bm\mS_i=\mS_i\mathrm{F}$ be the Stallings and Schreier graphs of the coset $H_ig_i$ (see Section~\ref{subsec:stallings}) with distinguished vertices $\iota_i$ and $\tau_i$, respectively.
	
	Consider the disjoint unions
	\[\mS_1\sqcup\cdots \sqcup\mS_n\subseteq \bm{\mS}_1\sqcup\cdots \sqcup \bm{\mS}_n. \]
The sets of vertices of these graphs will be the base sets for the structures to be defined below. The labelled edges shall indicate (partial) automorphisms of the structures in question. We first note that the free group $F$ acts by bijections on the vertex set of $\bm{\mS}:=\bm{\mS}_1\sqcup\cdots\sqcup \bm{\mS}_n$  on the right via $\rho\mapsto \rho w$ for every vertex $\rho$ and every $w\in F$, and for each individual~$i$
the vertex set of $\bm\mS_i$ is invariant under this action. 

We define binary relations $R_1,\dots,R_n$ on the vertex set $\mathrm{S}=\mathrm{V}(\bm{\mS})$ of $\bm{\mS}$ in order to produce a structure in $\mathbf{CCF}_n$ provided that $1_F\notin H_1g_1\cdots H_ng_n$. We set
\begin{equation}\label{def:defR_1...R_n}
	\tau_1\mathrel{R_1}\iota_2, \tau_2\mathrel{R_2}\iota_3,\dots, 
	\tau_{n-1}\mathrel{R_{n-1}}\iota_n,
	\tau_n\mathrel{R_n}\iota_1
\end{equation}
and add to this whatever is required so that the aforementioned action $\rho\mapsto \rho w$ by $F$ is an action by automorphisms of the structure. More precisely:  \[R_i=\{(\tau_i w,\iota_{i+1} w)\colon w\in F\}\mbox{ for }i=1,\dots,n-1\]
and
\[R_n=\{(\tau_n w,\iota_1 w)\colon w\in F\}.\]
We first show that, if $1_F\notin H_1g_1\cdots H_ng_n$ then $\bS:=(\mathrm{S};R_1,\dots, R_n)$ is a member of $ \mathbf{CCF}_n$. 
Suppose that $\mathfrak{S}=(\mathrm{S};R_1,\dots,R_n)\notin\mathbf{CCF}_n$. This means that there is a ``coloured cycle'': elements $\nu_1,\dots,\nu_n$ such that
\begin{equation}\label{eq:colored cycle}
\nu_1\mathrel{R_1} \nu_2\mathrel{R_2}\nu_3\cdots \nu_{n-1}\mathrel{R_{n-1}}\nu_n\mathrel{R_n}\nu_1.
\end{equation}
We note that $\nu_i\in \mathrm{V}(\bm\mS_i$) for every $i$, because every set $\mathrm{V}(\bm\mS_i)$ is invariant under the action of $F$. By definition of the relations $R_i$, there are $w_1,w_2,\dots,w_n\in F$ such that
\[
\begin{aligned}
	\ &\nu_1=\tau_1w_1\\
	\iota_2w_1=\ &\nu_2=\tau_2w_2\\
	&\ \vdots &\\
	\iota_nw_{n-1}=\ &\nu_n=\tau_n w_n\\
	\iota_1 w_n =\ &\nu_1
\end{aligned}
\]
This situation is depicted in Figure~\ref{fig:structureS}.
\begin{figure}[h]
\begin{tikzpicture}[scale=1.23]
	\filldraw(-1,1)circle(1.5pt);
	\draw(-1,1)[above]node{$\nu_1$};
	\filldraw(-1,2)circle(1.5pt);
	\draw(-1,2)[above]node{$\iota_1$};
	\filldraw(-3,1)circle(1.5pt);
	\draw(-2.95,1)[left]node{$\tau_1$};
	\filldraw(1,2)circle(1.5pt);
	\draw(1,2)[above]node{$\tau_n$};
	\filldraw(1,1)circle(1.5pt);
	\draw(1,1)[above]node{$\nu_n$};
	\filldraw(3,1)circle(1.5pt);
	\draw(3,1)[right]node{$\iota_n$};
	\filldraw(-4,-1)circle(1.5pt);
	\draw(-4,-1)[left]node{$\iota_2$};
	\filldraw(-2,-1)circle(1.5pt);
	\draw(-2,-1)[left]node{$\nu_2$};
	\filldraw(2,-1)circle(1.5pt);
	\draw(2,-1)[right]node{$\nu_{n-1}$};
	\filldraw(4,-1)circle(1.5pt);
	\draw(4,-1)[right]node{$\tau_{n-1}$};
	\filldraw(-3,-2)circle(1.5pt);
	\draw(-3,-2)[below]node{$\tau_2$};
	\filldraw(3,-2)circle(1.5pt);
	\draw(3,-2)[right]node{$\iota_{n-1}$};
	\draw plot [smooth cycle] coordinates  {(-0.5,2.4)(-0.7,0.5)(-3.5,0.5)(-2.5,2)};
	\draw plot [smooth cycle] coordinates  {(0.5,2.4)(0.7,0.5)(3.5,0.5)(2.5,2)};
	\draw plot [smooth cycle] coordinates  {(-2.5,-2.4)(-1.4,-0.5)(-4.5,-0.8)(-4,-2)};
	\draw plot [smooth cycle] coordinates  {(2.5,-2.4)(1.4,-0.5)(4.7,-0.8)(4,-2)};
	\draw[-latex](0.95,1)--(-0.95,1);
	\draw[-latex](0.95,2)--(-0.95,2);
	\draw[->,dashed](-0.1,1.9)--(-0.1,1.1);
	\draw(0.15,1.5)node{$w_n$};
	\draw(0,2)[above]node{$R_n$};
	\draw[-latex](-3.03,0.95)--(-3.97,-0.95);
	\draw[-latex](-1.03,0.95)--(-1.97,-0.95);
	\draw[->,dashed](-3.5,-0.1)--(-1.6,-0.1);
	\draw(-3.5,0)[left]node{$R_1$};
	\draw(-2.5,-0.2)[above]node{$w_1$};
	\draw[-latex](-3,-2)--(-1.5,-3.5);
	\draw[-latex](-2,-1)--(-0.5,-2.5);
	\draw(-2.4,-2.9)node{$R_2$};
	\draw[->,dashed](-1.95,-2.95)--(-1.05,-2);
	\draw(-1.7,-2.4)node{$w_2$};
	\draw[-latex](4,-1)--(3,0.93);
	\draw[-latex](2,-1)--(1,0.93);
	\draw(0,-3)node{$\dots$};
	\draw[->,dashed](3.5,-0.1)--(1.6,-0.1);
	\draw(2.5,-0.2)[above]node{$w_ {n-1 }$};
	\draw(4.4,0)[left]node{$R_ {n-1}$};
	\draw[latex-](2.95,-2.05)--(1.5,-3.5);
	\draw[latex-](1.95,-1.05)--(0.5,-2.5);
	\draw(2.6,-2.9)node{$R_{n-2}$};
	\draw[->,dashed](1.95,-2.95)--(1.05,-2);
	\draw(1.8,-2.4)node{$w_{n-2}$};
\end{tikzpicture}
\caption{The structure $\mathfrak{S}$}\label{fig:structureS}
\end{figure}
Recall the r\^oles of the vertices $\iota_i$ and $\tau_i$ in the Schreier graph of the coset $H_ig_i$:  from the equality $\iota_iw_{i-1}w_i\inv=\tau_i=\iota_ig_i$ (which holds in $\bm\mS_i$ for $i=2,\dots, n$) and $\iota_1w_nw_1\inv=\tau_1=\iota_1g_1$ (which holds in $\bm\mS_1$)
it follows that
\[w_nw_1\inv\in H_1g_1,\ w_1w_2\inv \in H_2g_2,\dots, 
\ w_{n-1}w_n\inv\in H_ng_n\]
and so
\[1_F=w_nw_1\inv w_1w_2\inv\cdots  w_{n-1}w_n\inv \in H_1g_1 H_2g_2\cdots H_ng_n.\]
Altogether, if $1_F\notin H_1g_1\cdots H_ng_n$ then 
${\bS}=(\mathrm{S};R_1,\dots,R_n)\in\mathbf{CCF}_n$.

Now consider the substructure $\bV$ of $\bS$ induced on the finite subset $\mathrm{V}=\mathrm{V}(\mS)$  where 
\begin{equation}\label{eq:substructure}
	\mS:=\mS_1\sqcup\cdots \sqcup \mS_n.
\end{equation}
First note that \eqref{def:defR_1...R_n} holds in $\bV$.
The $A$-graph structure on $\mS$ indicates partial automorphisms of the $\sigma_n$-structure $\bV$. Since $\bS$ is an extension in $\mathbf{CCF}_n$ in which all these partial automorphisms extend to total automorphisms, there exists a finite such extension, say $\ol{\bV}=(\ol{V}; R_1,\dots, R_n)\in \mathbf{CCF}_n$. Recall that we intend to show that there exists a finite index normal subgroup $N$ of $F$ such that $1_F\notin NH_1g_1\cdots NH_ng_n$.

We obtain the normal subgroup $N$ from the finite $A$-generated group $G$ of automorphisms of $\ol{\mathfrak{V}}$ that extend the partial automorphisms given by $\til A$ on $\mathfrak V$,   
by putting $F/N := G$. 
The free group $F$ acts on $\ol{\bV}$ on the right by automorphisms via $\rho w:=\rho[w]_G$. For $\rho\in \bV$ and $w\in F$, if  $\rho\cdot w$  is defined in $\bV$ (where $\rho\mapsto \rho\cdot w$ denotes the partial bijection on $\mathrm{V}=\mathrm{V}(\mS)$ induced by $w$) then it has the same value in $\ol{\bV}$, that is $\rho\cdot w=\rho[w]_G=\rho w$. 

We claim that $1_F\notin NH_1g_1\cdots NH_ng_n$. Suppose towards a contradiction this were not the case. Then, for every $i=1,\dots, n$ there is $k_i\in NH_ig_i$ such that $1_F=k_1k_2\cdots k_n$. We note that $\iota_i k_i=\tau_i$ for every $i$: indeed, $k_i=n_ih_ig_i$ for some $n_i\in N$ and $h_i$ in $H_i$. Now $\iota_in_i=\iota_i$ since $[n_i]_G=1_G$. Moreover, $\iota_i\cdot h_i=\iota_i$ and $\iota_i\cdot g_i=\tau_i$ hold in $\mathfrak{V}$.  By use of \eqref{def:defR_1...R_n}, which holds in $\bV$, the following hold in $\ol{\bV}$:

\[
\begin{aligned}
\tau_1 k_2\cdots k_n & \mathrel{R_1} \iota_2 k_2\cdots k_n\\
	=\tau_2 k_3\cdots k_n &\mathrel{R_2}\iota_3k_3\cdots k_n\\
		&\vdots\\
	=\tau_{n-1}k_n & \mathrel{R_{n-1}} \iota_nk_n\\ =\tau_n&\mathrel{R_n}\iota_1.
		\end{aligned}
\]
Since $\iota_1=\iota_1k_1\cdots k_n=\tau_1k_2\cdots k_n$ this  leads to a ``coloured cycle'' as in \eqref{eq:colored cycle}, a contradiction to the fact that $\ol{\mathfrak{V}}\in \mathbf{CCF}_n$.
\end{proof}
The Ribes--Zalesskii Theorem is an immediate consequence: let $H_1,\dots, H_n$ be finitely generated subgroups of $F$ and $g\in F$. Suppose that $g\notin H_1\cdots H_n$; this is equivalent to $1_F\notin H_1\cdots H_ng\inv$. By Proposition~\ref{prop:HL>RZ} this implies $1_F\notin NH_1\cdots NH_ng\inv=NH_1\cdots H_ng\inv$ for some finite index normal subgroup $N$ of $F$. Hence the set $NH_1\cdots H_ng\inv$  is a union of finitely many cosets $Nf_1,\dots, Nf_s$ (since $N$ is finite index) for which $Nf_i\ne N$ for all $i$. It follows that   $N\cap H_1\cdots H_ng\inv\subseteq  N\cap  NH_1\cdots H_ng\inv=\varnothing$. Multiplication by $g$ on the right gives $Ng\cap H_1\cdots H_n=\varnothing$.

\subsection{Ribes--Zalesskii implies Ash}	
The following is a consequence of the Ribes--Zalesskii Theorem:
\begin{Thm}\label{thm:ash for cycles} For every finite $A$-generated inverse monoid $M$ and every positive integer $n$ there exists a finite $A$-generated group $G$ such that for every $n$-tuple $v_1,\dots,v_n$ of words in $\til A^*$ for which $[v_1\cdots v_n]_G=1_G$ there exists an $n$-tuple of  words $u_1,\dots,u_n$ in $\til A^*$ such that
	\begin{enumerate}
		\item $[u_i]_M=[v_i]_M$ for all $i$, and
		\item $[u_1\cdots u_n]_F= 1_F$.
	\end{enumerate}
\end{Thm}
\begin{proof}
	We derive this claim from the Theorem of Ribes--Zalesskii.   We note that for every finitely generated subgroup $H$ of $F$, every conjugate subgroup ${}^w\!H:=wHw\inv$ is also finitely generated. Hence every subset of $F$ of the form \[H_1w_1\cdot H_2w_2\cdots H_nw_n\] for the $H_i$ finitely generated subgroups of $F$ and  $w_i\in F$ is also closed. Indeed, one has
	\[H_1w_1\cdot H_2w_2\cdots H_nw_n=H_1\cdot{}^{w_1}\!H_2\cdot{}^{w_1w_2}\!H_3\cdots{}^{w_1\cdots w_{n-1}}\!H_n w_1w_2\cdots w_n,\] 
	so the right hand side is a product of $n$ finitely generated subgroups of $F$ multiplied by the element $w_1\cdots w_n$. Since right multiplication $x\mapsto xw$ for some fixed element $w\in F$ is a homeomorphism of $F$ with respect to the profinite topology, the latter set is also closed.

Let $\varphi\colon M\to 2^F$ be the mapping defined by
	\[\varphi(m):=\big\{[w]_F\colon  w\in \til A^*,
    [w]_M=m\big\}.\]  Proposition~\ref{prop:canonicalrelmorMto F} tells us that for
 every $m\in M$, $\varphi(m)=Hw$ for some finitely generated subgroup $H$ of $F$ and some $w\in F$.

For every $n$-tuple $(m_1,m_2,\dots,m_n)$ of elements $m_i\in M$ consider the corresponding set
	\[\varphi(m_1)\cdot\varphi(m_2)\cdots \varphi(m_n)=H_1w_1\cdot H_2w_2\cdots H_nw_n,\]
which is a closed subset of $F$. Such a set may or may not contain the identity element $1_F$. For every such set $H_1w_1\cdot H_2w_2\cdots H_nw_n$ that does not contain $1_F$ there exists a finite index normal subgroup $N(m_1,\dots,m_n)$ of $F$ such that $N(m_1,\dots,m_n)\cap H_1w_1\cdot H_2w_2\cdots H_nw_n=\varnothing$. 
The intersection of all such (finitely many) normal subgroups is again a finite index normal subgroup $N$ of $F$. By definition, $N$ satisfies the following. For every $n$-tuple
$(m_1,\dots,m_n)\in M^n$, if $1_F\notin
	\varphi(m_1)\cdots\varphi(m_n)$, then $N\cap
	\varphi(m_1)\cdots\varphi(m_n)=\varnothing$. The finite group $G:=F/N$
fulfils the requirements of the theorem. Indeed, let $v_1,\dots,
	v_n\in \til A^*$ be such that $[v_1\cdots v_n]_G=1_G$,
which means that $[v_1\cdots v_n]_F\in N$. We let $m_i:=[v_i]_M$ for
$i=1,\dots, n$. Since $[v_i]_F\in \varphi(m_i)$, we also
have \[[v_1\cdots v_n]_F=[v_1]_F\cdots [v_n]_F\in \varphi(m_1)\cdots
	\varphi(m_n).\]
So in particular $N\cap \varphi(m_1)\cdots \varphi(m_n)\ne \varnothing$,
from which it follows by construction of $N$ 
that $1_F\in \varphi(m_1)\cdots \varphi(m_n)$. Consequently,  $1_F=g_1\cdots g_n$ for some $g_i\in \varphi(m_i)$.
By definition of $\varphi$ there exist words $u_1,\dots u_n$ such that $[u_i]_F=g_i$ and $[u_i]_M=m_i=[v_i]_M$ for every $i$,
and 
$1_F=g_1\cdots g_n=[u_1\cdots u_n]_F$.\end{proof}
    
It should be mentioned that, conversely, Theorem~\ref{thm:ash for cycles} implies the Ribes--Zalesskii Theorem. Indeed let $H_1,\dots, H_n$ be finitely generated subgroups of $F$ and let $g\in F$ (considered as a reduced word over $\til A$). Take the graph $\mS$ in \eqref{eq:substructure} with $\mS_i$ the stallings graph of $H_i$ for $i=1,\dots, n-1$ and $\mS_n$ the Stallings graph of the coset $H_ng\inv$, and let $M=\cT(\mS)$.\footnote{The graph $\mS$ plays a crucial r\^ole in~\cite{ASconstructiveRZ}. }  We have that $g\in H_1\cdots H_n$ if and only if $1_F\in H_1\cdots H_ng\inv$ and it is easy to see that the latter holds if and only if there are words $u_1,\dots, u_n\in \til A^*$ such that $[u_1\cdots u_n]_F=1_F$ and $\iota_i\cdot u_i=\iota_i$ for $i=1,\dots,n-1$ and $\iota_n\cdot u_n=\iota_n\cdot g\inv=\tau_n$ holds in $\mS$. Let $G$ be a finite $A$-generated group satisfying Theorem~\ref{thm:ash for cycles} for the inverse monoid $M$. Let $\varphi\colon F\twoheadrightarrow G$ be the canonical morphism and suppose that $\varphi(g)\in \varphi(H_1\cdots H_n)$. Then there exist (reduced) words $h_1,\dots,h_n$ with $h_i\in H_i$ for all $i$ such that $\varphi(g)=\varphi(h_1\cdots h_n)$, that is, $[g]_G=[h_1\cdots h_n]_G$ so that $[h_1\cdot \ldots\cdot h_ng\inv]_G=1_G$. By Theorem~\ref{thm:ash for cycles} there exist words $u_1,\dots, u_n$ such that $[u_1\cdots u_n]_F=1_F$, $[u_i]_M=[h_i]_M$ (for $i=1,\dots, n-1$) and $[u_n]_M=[h_ng\inv]_M$. The latter condition implies $\iota_i\cdot u_i=\iota_i$ (for $i=1,\dots, n-1$) and $\iota_n\cdot u_n=\iota_n\cdot g\inv=\tau_n$, which in turn gives $g\in H_1\cdots H_n$. Consequently, if $g\notin H_1\cdots H_n$, then $g$ and $H_1\cdots H_n$ are separated by the morphism $\varphi$.

Observe that Theorem~\ref{thm:ash for cycles} is Theorem~\ref{thm: ash} for the special case of cycle graphs. {We have thus proved the equivalence of Ash's Theorem for cycle graphs with the Theorem of Ribes--Zalesskii (a more long-winded proof for this equivalence may be found in~\cite{steinberginverseautomata}.)} We now derive Theorem~\ref{thm: ash} from Theorem~\ref{thm:ash for cycles}.
The following concept will be convenient for the proof. 
It is essentially the same as that of an ``associated tree'' in~\cite[Section 7]{Ash}.

\begin{Def}[critical group]\label{def: critical group}
Let $M$ be a finite $A$-generated inverse monoid, $\mE$ a finite connected graph  without loop edges and $G$ a finite $A$-generated group extending $M$. We say that $G$ is \emph{critical for the pair} $(M,\mE)$ if the following condition is satisfied. For every $G$-commuting labelling $v\colon\mathrm{K}(\mE)\to \til A^*$ there exists a finite tree $\mT$ of diameter at most $|\mathrm{V}(\mE)|-1$ such that $\mathrm{V}(\mE)\subseteq \mathrm{V}(\mT)$, and a labelling $w\colon \mathrm{K}(\mT)\to \til A^*$ such that for every edge $e\in \mathrm{K}(\mE)$, if $\pi_e$ denotes the (unique) reduced path ${\alpha e}\longrightarrow {\omega e}$ in $\mT$ then $[v(e)]_M\le[w(\pi_e)]_M$.
\end{Def}

Note that every word labelling of a finite tree commutes over $F$. A group $G$ critical for $(M,\mE)$ guarantees that every labelling of $\mE$ 
that is acknowledged by $G$ to commute admits an $M$-related $F$-commuting relabelling, as we will see immediately after the statement of the following key lemma of this subsection.  
\begin{Lemma}[key lemma]\label{lem: key lemma RZ->Ash}
For every finite $A$-generated inverse monoid $M$  and every finite connected graph  $\mE$  without loop edges there exists a finite $A$-generated group $G$ which is critical for the pair $(M,\mE)$.
\end{Lemma}
We first show how this lemma implies Ash's Theorem.
\begin{proof}[Proof of Theorem~\ref{thm: ash} by use of Lemma~\ref{lem: key lemma RZ->Ash}]
	Let $M$ be a finite $A$-generated inverse monoid and let $\mE$ be a finite connected graph. Let $\mD$ be the graph obtained from $\mE$ by removing all loop edges (if any). Suppose that $G$ is critical for the pair $(M,\mD)$ and let $v\colon \mathrm{K}(\mE)\to \til A^*$ be a $G$-commuting labelling. Since every expansion of  $G$ is also critical, we may assume that $G$ extends $M$ in the sense of Definition~\ref{def: G extends M}, i.e.~$[x]_M=[xx\inv]_M$ for every word $x$ for which $[x]_G=1_G$. For every edge $e$ of $\mE$ we define $u(e):=v(e)v(e)\inv w(\pi_e)$ ($w$ and $\pi_e$ as in Definition~\ref{def: critical group}). If $e$ is a loop edge then $\pi_e$ is an empty path, and $[u(e)]_F=1_F$ and $[u(e)]_M=[v(e)]_M$ since $[v(e)]_G=1_G$ implies that $[v(e)]_M$ is idempotent. For non-loop edges (that is, edges of $\mD$), $[v(e)]_M\le[w(\pi_e)]_M$ implies $[v(e)]_M=[v(e)v(e)\inv w(\pi_e)]_M=[u(e)]_M$. 
    For a closed  path $e_1\cdots e_t$ in $\mE$ then
	\[[u(e_1\cdots e_t)]_F=[w(\pi_{e_1}\cdots \pi_{e_t})]_F=1_F,
	\]  
	since the path $\pi_{e_1}\cdots \pi_{e_t}$ is a closed path in the tree $\mT$ and therefore the corresponding reduced path is the empty path at $\alpha\pi_{e_1}={\alpha e_1}$.

    So far we have obtained Theorem~\ref{thm: ash} for finite connected graphs. Now let $\mE$ be any finite graph with connected components $\mE_1,\dots,\mE_n$. For every $i$ there is a finite $A$-generated group $G_i$ serving as in Theorem~\ref{thm: ash} for $M$ and $\mE_i$. Then the $A$-generated direct product $G:=G_1\times\cdots\times G_n$ serves as required for the inverse monoid $M$ and the graph $\mE=\mE_1\sqcup\cdots\sqcup \mE_n$.
\end{proof}
\begin{Rmk}
	The above proof of Theorem~\ref{thm: ash} does not use 
		the bound on the diameter of the tree $\mT$ 
	in  Definition~\ref{def: critical group}  (critical group).
	That condition will be used, however, for the inductive procedure in the proof of Lemma~\ref{lem: key lemma RZ->Ash}. 
\end{Rmk}

\begin{proof}[Proof of Lemma~\ref{lem: key lemma RZ->Ash} by use of Theorem~\ref{thm:ash for cycles}]
   This proof is from \cite[Section 7]{Ash}, but polished, with many details included and adjusted to the present context.
   In the first part we establish the claim for cycle graphs 
   of length up to $n$ in the r\^ole of $\mE$, 
   uniformly for all $M$. In the second part, this is lifted to the general
   case of connected graphs $\mE$ on $n$ vertices, essentially by induction  
   on the cyclomatic number: surgery that eliminates one simple cycle of 
    $\mE$ reduces the claim to the claims for lower cyclomatic number, 
    cycle graphs and trees.
Altogether we derive Lemma~\ref{lem: key lemma RZ->Ash} from Theorem~\ref{thm:ash for cycles}, hence from Theorem~\ref{thm: Ribes--Zalesskii}. 

\medskip
\emph{To the first part.} For any positive integer $n\ge 2$ 
 let $\mC_n$ be the cycle graph of length $n$. From Theorem~\ref{thm:ash for cycles} we get the existence of a critical group for every pair $(M,\mC_n)$. Indeed, let $G$ be as in the statement of Theorem~\ref{thm:ash for cycles} and let $v\colon \mathrm{K}(\mC_n)\to \til A^*$ be a $G$-commuting labelling of $\mC_n$. Denote the vertices of $\mC_n$ by $o_1,\dots, o_n$ and consider  edges $e_i\colon o_{i-1}\longrightarrow o_{i}$ (indices to be taken $\bmod\ n$). 
Then $[v_1\cdots v_n]_G=1_G$ for $v_i=v(e_i)$. 
According to Theorem~\ref{thm:ash for cycles} there are words $u_1,\dots,u_n$ such that $[u_i]_M=[v_i]_M$ and [$u_1\cdots u_n]_F=1_F$. The claim that $G$ is critical for the pair $(M,\mC_n)$ follows from the following lemma. Note that considering labellings $v$ of $\mC_n$ satisfying $v(e_{t+1})=\cdots =v(e_n)=1$ (the empty word) 
we  find that $G$ is critical also for every pair $(M,\mC_t)$ with $t<n$. 

\begin{Lemma}\label{lem:1stpart} Let $u_1,\dots, u_n$ be words such that $[u_1\cdots u_n]_F=1_F$ and let $o_1,\dots, o_n$ be the vertices of the cycle graph $\mC_n$ of length $n$.
Then there exists a finite tree $\mT$ of diameter at most $n-1$ with $\mathrm{V}(\mT)\supseteq\{o_1,\dots,o_n\}$ and a labelling $w\colon\mathrm{K}(\mT)\to \til A^*$ such that $[u_i]_M\le [w(\pi_i)]_M$ where $\pi_i$ is the reduced path $o_{i-1}\longrightarrow o_{i}$ in $\mT$ (indices to be taken $\bmod\ n$). 
\end{Lemma}

\begin{proof}
Consider first the special case that the elements $[u_1\cdots u_i]_F$ ($i=1,\dots, n$) are pairwise distinct; this means that for every proper factor $u_i\cdots u_j$ of $u_1\cdots u_n$ we have $[u_i\cdots u_j]_F\ne 1_F$.  
Let $\mC_n$ be the cycle graph of length $n$ with edges $e_i\colon o_{i-1}\longrightarrow o_i$ (indices to be taken $\bmod\ n$) and labelling $u\colon e_i\mapsto u_i$. For every $i$ let $v_i:=\mathrm{red}(u_i)$. Note that no word $v_i$ is empty. Now form a new  graph $\mR$  by replacing in $\mC_n$ every edge $e_i\colon o_{i-1}\longrightarrow o_i$ by an arc of the form 
\begin{equation}\label{eq:new arc}
\vartheta_i\colon o_{i-1}\longrightarrow p_{i1}\longrightarrow\cdots\longrightarrow p_{i,k_{i}-1}\longrightarrow o_i
\end{equation} 
of length $k_i$ the length of the word $v_i$, with new vertices $p_{ij}$  and the edges being labelled by the letters of $v_i$ so that $\ell(\vartheta_i)= v_i$. Assume that distinct arcs $\vartheta_i$ and $\vartheta_j$ have no edges and no intermediate vertices $p_{ir}$ and $p_{js}$ in common. The result is an $A$-labelled cycle graph of length $|v_1\cdots v_n|$. Next apply Stallings foldings to get the maximum $A$-graph quotient $\mW:=\mR/\Theta$ of $\mR$. No two distinct vertices $o_i$ and $o_j$ of $\mR$ are identified by the folding process since in $\mR$ there is no path $o_i\longrightarrow o_j$ labelled by a Dyck word. Hence $\{o_1,\dots,o_n\}$ can indeed be realised as a subset of $\mathrm{V}(\mW)$. In the folding process no two vertices $p_{is}, p_{it}$ of any one of the arcs $\vartheta_i$~\eqref{eq:new arc} are identified with each other since no two such vertices are connected by a path labelled by a Dyck word. It follows that in $\mW$ we have a path $o_{i-1}\longrightarrow o_i$ labelled $v_i$ for every $i$. 
Since $v_1\cdots v_n$ is a Dyck word, the folding process transforms the cycle graph $\mR$ into a tree: this follows easily by induction on the length of the word $v_1\cdots v_n$. 
That is, $\mW$ is a finite tree in which
the label of the reduced path  ${o_{i-1}}\longrightarrow{o_{i}}$ is  $v_i=\mathrm{red}(u_{i})$, and clearly 
we have $[u_{i}]_M\le[\mathrm{red}(u_{i})]_M$. We need to modify $\mW$ a bit in order to get the required bound on the diameter of $\mT$. For this purpose, consider all arcs spanned by simple paths
\begin{equation}\label{eq:bad path Ash}
	q_0\underset{f_1}{\longrightarrow}q_1\underset{f_2}{\longrightarrow}\cdots\underset{f_{t-1}}{\longrightarrow} q_{t-1}
	\underset{f_t}{\longrightarrow}q_t 
\end{equation} of consecutive edges $f_1,\dots, f_t$ in $\mW$ such that the intermediate vertices $q_1,\dots, q_{t-1}$ are all of degree $2$ and none is of the form ${o_i}$ while each of the vertices $q_0$ and $q_t$  is of degree at least $3$ or is of the form ${o_i}$ for some $i$. If some reduced path between any two vertices $o_i$ and $o_{j}$ in $\mW$ traverses one of the edges $f_s$ of \eqref{eq:bad path Ash} then it needs to traverse the entire arc \eqref{eq:bad path Ash}. Hence we may contract the path \eqref{eq:bad path Ash} to a single edge, that is, replace every arc spanned by a path of the form \eqref{eq:bad path Ash} by the single edge  $f\colon q_0\longrightarrow q_t$ (and its inverse) and label $f$ accordingly, that is, $w(f)$ is the string of the (letter) labels of the edges $f_1,\dots,f_t$; for all other edges $f$, let $w(f)$ be the label induced by the labelling of $\mW$. 
In the resulting tree $\mT$ then, the label of every reduced path ${o_{i-1}}\longrightarrow {o_{i}}$ is the same as in $\mW$  (namely $v_i=\mathrm{red}(u_i)$). The tree $\mT$ now has at most $n$ vertices of degree up to $2$, hence has diameter at most $n-1$ as required by Definition~\ref{def: critical group} (critical group). 

The general case, namely that some of the elements  $[u_1\cdots u_i]_F$ ($i=1,\dots,n$) coincide, can be done by induction on $n$ where the case $n=1$ is trivial. So let $n>1$ and 
assume the claim of the lemma to be true 
for all products $v_1\cdots v_k$  of length $k$ smaller than $n$, and let $u_1,\dots, u_n$ with $[u_1\cdots u_n]_F=1_F$ be given. 
If $[u_i\cdots u_j]_F\ne 1_F$ for every proper factor $u_i\cdots u_j$ of $u_1\cdots u_n$, then we are done, by the argument above. 
So let us assume that $[u_i\cdots u_j]_F=1_F$ for some proper factor $u_i\cdots u_j$ of $u_1\cdots u_n$. It is sufficient --- possibly after some cyclic permutation of the words $u_1,\dots, u_n$ and some change of indices --- to consider the case $[u_1\cdots u_k]_F=1_F$ for some $k<n$. Then also $[u_{k+1}\cdots u_n]_F=1_F$.  By the inductive assumption, there exist finite trees $\mT_1$ and $\mT_2$ with $\mathrm{V}(\mT_1)\supseteq \{o_1,\dots,o_k\}$ and $\mathrm{V}(\mT_2)\supseteq \{o_{k+1},\dots, o_n\}$, of diameters at most $k-1$ and $n-k-1$, respectively,  and labellings $w_i\colon \mathrm{K}(\mT_i)\to \til A^*$ ($i=1,2$) according to the claim. This means that (i) for $1\le i\le k$,  if $\varrho_i$ is the reduced path $o_{i-1}\longrightarrow o_i$ in $\mT_1$ (indices to be taken $\bmod\ k$), then $[u_i]_M\le [w_1(\varrho_i)]_M$ and (ii) for $k+1\le i\le n$, if $\varsigma_i$ is the reduced path $o_{i-1}\longrightarrow o_i$ in $\mT_2$ (indices to be taken $\bmod\ n-k$), then $[u_i]_M\le [w_2(\varsigma_i)]_M$.  We define the tree $\mT$ by 
\[\mT:=\mT_1\sqcup \mT_2\sqcup \{f,f\inv\}\]
with $f\colon o_k\longrightarrow o_n$ a new edge connecting the vertices $o_k$ of $\mT_1$ and $o_n$ of $\mT_2$. We set $w_f(f^{\pm 1}):=1$ and \[w:=w_1\sqcup w_2\sqcup w_f\] to get a labelling of $\mT$. We note that in case (i) above the reduced path $o_{i-1}\longrightarrow o_i$ in $\mT$ coincides with the one in $\mT_1$ provided that $i>1$ while in case (ii) the reduced  path $o_{i-1}\longrightarrow o_i$ in $\mT$ coincides with the one in $\mT_2$ provided that $i>k+1$. This implies that $[u_i]_M\le[w(\pi_i)]_M$ for $i=2,\dots, k,k+2,\dots, n$ where $\pi_i$ is the reduced path $o_{i-1}\longrightarrow o_i$ in $\mT$. The reduced path $\pi_1\colon o_n\longrightarrow o_1$ in $\mT$ is given by $\pi_1=f\inv\varrho_1$ while the reduced path $\pi_{k+1}\colon o_k\longrightarrow o_{k+1}$ in $\mT$ is $\pi_{k+1}=f\varsigma_{k+1}$. The inequality $[u_i]_M\le[w(\pi_i)]_M$ 
therefore holds also for $i=1$ and $i=k+1$, as required. 
Since, by the inductive assumption,  $\mT_1$ has diameter at most $k-1$ and $\mT_2$ has diameter at most $n-k-1$, $\mT$ has diameter at most $n-1$.
\end{proof}

\emph{To the second part.}
We consider all connected graphs $\mE$ without loop edges on $n$ vertices. We show the existence of a critical group for every pair $(M,\mE)$ by induction on the number $l$ of geometric edges of $\mE$, simultaneously for all $M$. (In effect, the induction is on the cyclomatic number.) 
The base case 
of the induction is $l=n-1$, in which case the graph itself is a tree; we may take $\mT=\mE$, $w=v$ and $G$ any group, in particular any group $G$ that extends $M$ in the sense of Definition~\ref{def: G extends M}. So let $l>n-1$ and suppose that the claim is true for all connected graphs $\mD$ on $n$ vertices with fewer than $l$ geometric edges and all finite $A$-generated inverse monoids $N$. Let $\mE$ be given with $l$ geometric edges, let $k$ be an edge of $\mE$ that is contained in a non-trivial cycle subgraph of $\mE$, and let $\mD=\mE\setminus\{k,k\inv\}$. Note that $\mathrm{V}(\mE)=\mathrm{V}(\mD)$. 
With the first part of the proof we obtain a group $G$  critical for $(M,\mC_n)$ and hence, as remarked at the end of the last paragraph before Lemma~\ref{lem:1stpart}, also for $(M,\mC_t)$ for all $t\le n$. 
By the inductive assumption, the claim of the lemma is true for any pair $(N,\mD)$. 
So let $H$ be an expansion of $G$ that is critical for the pair $(M\times G,\mD)$. 
We claim that $H$ is critical for $(M,\mE)$; for later reference we state:\footnote{{In his proof in~\cite{Ash}, Ash used the Margolis--Meakin expansion of $G$ (to be defined in the discussion after the proof of Lemma~\ref{lem:main lemma}) instead of the direct product $M\times G$. Apart from making the proof more transparent, the replacement of the Margolis--Meakin expansion of $G$ by the direct product $M\times G$ is essential for the understanding of the arguments in Section~\ref{sec:groups for Ash}.}}
\begin{equation}\label{eq:Hcritical for (MtimesG)}
\begin{split}
 G \mbox{ critical for } (M,\mC_n) 
    \;\mbox{ and }\; H\mbox{ critical for } (M\times G,\mD) \; \\ 
    \Longrightarrow \;\; H\mbox{ critical for }(M,\mE).   
\end{split}
\end{equation}

Let $v\colon \mathrm{K}(\mE)\to \til A^*$ be an $H$-commuting labelling of $\mE$. We choose edges $k_1,\dots,k_t$ in $\mD=\mE\setminus\{k,k\inv\}$ so that $k_1\cdots k_t\colon \omega k\longrightarrow \alpha k$ is a simple path (that is, there is no repetition of vertices)  in $\mD$ and thus $kk_1\cdots k_t$ is a simple closed path in $\mE$. Since $v$ is $H$-commuting, 
\begin{equation}\label{eq: compatible H}
	[v(k)v(k_1\cdots k_t)]_H=1_H.
\end{equation}
Since the labelling $v\restriction \mD$ is also $H$-commuting, there exists,
by the inductive assumption, a finite tree $\mX$ of diameter at most $n-1$ such that
\begin{equation}\label{eq: mapping VDtoVX}
\mathrm{V}(\mD)\subseteq \mathrm{V}(\mX)
\end{equation}
 and a labelling $w_\mX\colon \mathrm{K}(\mX)\to\til A^*$ such that for every edge $e$ of $\mD$, if $\pi_e$ is the reduced path ${\alpha e}\longrightarrow{\omega e}$ in $\mX$, then $[v(e)]_{M\times  G}\le[w_\mX(\pi_e)]_{M\times G}$. This holds in particular for the edges $k_i$: $[v(k_i)]_{M\times G}\le[w_\mX(\pi_{k_i})]_{M\times G}$ for all $i$. As a consequence we get
\begin{equation}\label{eq: equality G}
[v(k_i)]_{G}=[w_\mX(\pi_{k_i})]_{G}\mbox{ for }	i=1,\dots,t,
\end{equation}
since the relation $[v]_{M\times G}\le [w]_{M\times G}$ implies $[v]_G\le [w]_G$ ($G$ is a quotient of $M\times G$) and $\le$ is the identity relation in the group $G$.

Next, we look more carefully at the reduced path ${\omega k}\longrightarrow {\alpha k}$ in $\mX$. We write it explicitly as
\begin{equation*}
	{\alpha k_1}=p_0\underset{f_1}{\longrightarrow}p_1\underset{f_2}{\longrightarrow}p_2\cdots p_{s-1}
	\underset{f_s}{\longrightarrow}p_s={\omega k_t}
\end{equation*} 
and record that $s\le n-1$ by the assumption on the diameter of $\mX$. Let $\mC$ be the cycle graph spanned by the path $f_1\cdots f_s$ augmented by an extra edge $f_0\colon p_s={\alpha k}\longrightarrow {\omega k}=p_0$ (and its inverse), so that $\mC$ is a cycle graph of length $s+1\le n$. We endow $\mC$ with the labelling $w_\mC$ defined by

\begin{equation*}
	w_\mC(f_i)=  
	\begin{cases} w_\mX(f_i) &\mbox{ for }i=1,\dots, s\\
		v(k) &\mbox{ for }i=0.
	\end{cases} 
\end{equation*}
The path $f_1\cdots f_s$ is the reduced form of the path $\pi_{k_1}\cdots \pi_{k_t}$, and hence 
\[[w_\mX(f_1\cdots f_s)]_F=[w_\mX(\pi_{k_1}\cdots \pi_{k_t})]_F\] so that
\begin{equation}\label{eq: equal G}
[w_\mX(f_1\cdots f_s)]_G=[w_\mX(\pi_{k_1}\cdots \pi_{k_t})]_G.
\end{equation}
Furthermore, \eqref{eq: compatible H} implies $[v(k)v(k_1\cdots k_t)]_G=1_G$ since $H\twoheadrightarrow G$. 
In combination with~\eqref{eq: equality G} and~\eqref{eq: equal G} this 
implies that the labelling $w_\mC$ is $G$-commuting. It follows that there exists a finite tree $\mY$ such that $\mathrm{V}(\mC)\subseteq \mathrm{V}(\mY)$,  and a labelling $w_\mY\colon \mathrm{K}(\mY)\to \til A^*$ such that, if $\rho_i$ is the reduced path ${p_{i-1}}\longrightarrow {p_i}$ in $\mY$, then
\begin{equation}\label{eq: ineq M}
	[w_\mC(f_i)]_M\le[w_\mY(\rho_i)]_M \mbox{ for }i=0,\dots,s
\end{equation}
(indices to be taken $\bmod\ s+1$). We now remove the edges $f_i^{\pm1}$ from the tree $\mX$ and form the disjoint union with $\mY$, the result of which is the forest
\[\mU:=\mY\sqcup\left.\mX\!\setminus\! \{f_1^{\pm 1},\dots,f_s^{\pm1}\}\right.\]
consisting of the connected components $\mY$ and $\mX_0,\dots,\mX_s$,  where $\mX_i$ is the component of $p_i$ in the forest $\mX\!\setminus\!\{f_1^{\pm 1},\dots,f_s^{\pm1}\}$. Let $\mZ:=\mU/\Theta$ where $\Theta$ is the congruence on $\mU$ that identifies every vertex $p_i$ of $\mX_i$  with the vertex ${p_i}$ in $\mY$. Restricted to the tree $\mY$ and each of the trees $\mX_i$, $\Theta$ is the identity relation. Then $\mZ$ is a tree and we may consider $\mY$ and every $\mX_i$ as subgraphs of $\mZ$ and such that $\mY\cap \mX_i=\{p_i\}$, and we identify each vertex with its $\Theta$-class. In particular, we have the inclusions
\[\mathrm{V}(\mE)=\mathrm{V}(\mD)\subseteq \mathrm{V}(\mX)\subseteq\mathrm{V}(\mZ).\] Figure~\ref{fig:trees} gives an overview of the various graphs involved. 
\begin{figure}[ht]
	\begin{tikzpicture}
		\draw[->](0,0)--(0.95,0);
		\draw[->](1.05,0)--(1.95,0);
		\draw[->](2.05,0)--(2.95,0);
		\draw[->](3.05,0)--(3.95,0);
		\draw[dotted](4,0)--(5,0);
		\draw[->](5.05,0)--(5.95,0);
		\draw(0,1)--(0,0);
		\draw(1,1)--(1,0);
		\draw(2,1)--(2,0);
		\draw(3,1)--(3,0);
		\draw(4,1)--(4,0);
		\draw(5,1)--(5,0);
		\draw(6,1)--(6,0);
		\draw([above](0,1.2)node{$\mathcal{X}_{0}$};
		\draw([above](1,1.2)node{$\mathcal{X}_{1}$};
		\draw([above](2,1.2)node{$\mathcal{X}_{2}$};
		\draw([above](3,1.2)node{$\mathcal{X}_{3}$};
		\draw([above](4,1.2)node{$\mathcal{X}_{4}$};
		\draw([above](5,1.2)node{$\mathcal{X}_{s-1}$};
		\draw([above](6,1.2)node{$\mathcal{X}_{s}$};
		\filldraw(0,0)circle(1.5pt);
		\filldraw(1,0)circle(1.5pt);
		\filldraw(2,0)circle(1.5pt);
		\filldraw(3,0)circle(1.5pt);
		\filldraw(4,0)circle(1.5pt);
		\filldraw(5,0)circle(1.5pt);
		\filldraw(6,0)circle(1.5pt);
		\draw[below](0.5,0)node{$f_1$};
		\draw[below](1.5,0)node{$f_2$};
		\draw[below](2.5,0)node{$f_3$};
		\draw[below](3.5,0)node{$f_4$};
		\draw[below](5.5,0)node{$f_s$};
		\draw[left](0,0)node{${\omega k}$};
		\draw[right](6,0)node{${\alpha k}$};
		\draw(2.5,0.5)node{$\mathcal{X}$};
	\end{tikzpicture}
	\begin{tikzpicture}
		\draw[thick, dotted](6,-0.1) .. controls (5.5,-2.7) and (0.5,-2.7) .. (0,-0.05);
		\draw[thick, dotted](0,0)--(0.95,0);
		\draw[thick, dotted](1.05,0)--(1.95,0);
		\draw[thick, dotted](2.05,0)--(2.95,0);
		\draw[thick, dotted ](3.05,0)--(3.95,0);
		\draw[dotted](4,0)--(5,0);
		\draw[thick, dotted](5.05,0)--(5.95,0);
		\filldraw(0,0)circle(1.5pt);
		\filldraw(1,0)circle(1.5pt);
		\filldraw(2,0)circle(1.5pt);
		\filldraw(3,0)circle(1.5pt);
		\filldraw(4,0)circle(1.5pt);
		\filldraw(5,0)circle(1.5pt);
		\filldraw(6,0)circle(1.5pt);
		\draw(0,0)--(0.3,-0.5)--(1.5,-1.3)--(2.5,-1.8)--(4.5,-1.5)--(5.5,-0.5)--(6,0);
		\draw(0.3,-0.5)--(1,0);
		\draw(1.5,-1.3)--(2,0);
		\draw(3,0)--(3,-0.5)--(2.5,-1.8);
		\draw(3,-0.5)--(4,0);
		\draw(5.5,-0.5)--(5,0);
		\draw(2.5,-0.5)node{$\mathcal{Y}$};
		\draw(0,-2)node{$\phantom{}$};
		\draw(5,-1.8)node{$\mathcal{C}$};
	\end{tikzpicture}
    \begin{tikzpicture}
		\draw(0,1)--(0,0);
		\draw(1,1)--(1,0);
		\draw(2,1)--(2,0);
		\draw(3,1)--(3,0);
		\draw(4,1)--(4,0);
		\draw(5,1)--(5,0);
		\draw(6,1)--(6,0);
		\draw([above](0,1.2)node{$\mathcal{X}_{0}$};
		\draw([above](1,1.2)node{$\mathcal{X}_{1}$};
		\draw([above](2,1.2)node{$\mathcal{X}_{2}$};
		\draw([above](3,1.2)node{$\mathcal{X}_{3}$};
		\draw([above](4,1.2)node{$\mathcal{X}_{4}$};
		\draw([above](5,1.2)node{$\mathcal{X}_{s-1}$};
		\draw([above](6,1.2)node{$\mathcal{X}_{s}$};
		\filldraw(0,0)circle(1.5pt);
		\filldraw(1,0)circle(1.5pt);
		\filldraw(2,0)circle(1.5pt);
		\filldraw(3,0)circle(1.5pt);
		\filldraw(4,0)circle(1.5pt);
		\filldraw(5,0)circle(1.5pt);
		\filldraw(6,0)circle(1.5pt);
		\filldraw(0,0)circle(1.5pt);
		\filldraw(1,0)circle(1.5pt);
		\filldraw(2,0)circle(1.5pt);
		\filldraw(3,0)circle(1.5pt);
		\filldraw(4,0)circle(1.5pt);
		\filldraw(5,0)circle(1.5pt);
		()\filldraw(6,0)circle(1.5pt);
		\draw(4.5,0)node{$\dots$};
		\draw(0,0)--(0.3,-0.5)--(1.5,-1.3)--(2.5,-1.8)--(4.5,-1.5)--(5.5,-0.5)--(6,0);
		\draw(0.3,-0.5)--(1,0);
		\draw(1.5,-1.3)--(2,0);
		\draw(3,0)--(3,-0.5)--(2.5,-1.8);
		\draw(3,-0.5)--(4,0);
		\draw(5.5,-0.5)--(5,0);
		\draw(2.6,-0.9)node{$\mathcal{Y}$};
		\draw(3,-2.2)node{$\mathcal{Z}$};
	\end{tikzpicture}
	\caption{The trees $\mX, \mY$ and $\mZ$ and the cycle graph $\mathcal{C}$ (dotted)}\label{fig:trees}
\end{figure}%

On the set of edges, $\Theta$ is the identity relation, whence we may write $\mathrm{K}(\mZ)=\mathrm{K}(\mY)\cup\bigcup_{i=0}^s \mathrm{K}(\mX_i)$ (identifying each edge of $\mU$ with its $\Theta$-class).
 Next we endow $\mZ$ with its natural labelling: for any $f\in \mathrm{K}(\mZ)$ let
\begin{equation*}
	w_\mZ(f)=\begin{cases}
		w_\mX(f)&\mbox{ if }f\in \mX_i\mbox{ for some }i\\
		w_\mY(f)&\mbox{ if }f\in \mY.
		\end{cases}
\end{equation*} 
We want to show that the tree $\mZ$ with the 
labelling $w_\mZ$ satisfies the requirements 
of Definition~\ref{def: critical group} for $\mE$ (with the exception of the condition on the diameter of the required tree).
For this we show that, for every edge $e\in \mathrm{K}(\mE)$
and the associated reduced path $\varsigma_e \colon
{\alpha
	e}\longrightarrow{\omega e}$
 in $\mZ$, we have 
$[v(e)]_M\le[w_\mZ(\varsigma_e)]_M$.
We consider two cases: $e=k$ and $e\ne k^{\pm1}$, that is, $e\in \mD$.

Case~1: $e=k$. In this case, ${\alpha k}=p_s$ and ${\omega k}=p_0$ and the reduced path $\varsigma_k\colon{\alpha k}\longrightarrow{\omega k}$ in $\mZ$ coincides with the one in $\mY$, which was earlier denoted $\rho_0$. Since $v(k)=w_\mC(f_0)$ we have
\[[v(k)]_M=[w_\mC(f_0)]_M\le[w_\mY(\rho_0)]_M=[w_\mZ(\varsigma_k)]_M\]
by \eqref{eq: ineq M}.

Case~2: $e\ne k^{\pm 1}$. First recall that
\begin{equation}\label{eq: ineq for X}
	[v(e)]_{M\times G}\le[w_\mX(\pi_e)]_{M\times G}
\end{equation}
by the construction of $\mX$ where $\pi_e$ is the reduced path ${\alpha e}\longrightarrow{\omega e}$ in $\mX$. 
If $\pi_e$ does not use any of the edges $f_1^{\pm 1},\dots, f_s^{\pm 1}$, then the reduced path $\varsigma_e\colon \alpha e\longrightarrow \omega e$ in $\mZ$ coincides with $\pi_e$ and we are done. 
But it may happen that the path $\pi_e$ traverses a segment of the arc spanned by $f_1\cdots f_s$, that is,
\begin{equation}\label{eq:options}
	\pi_e=\pi_1f_r\cdots f_t\pi_2 \;\;\mbox{ or }\;\; \pi_e=\pi_2f_t\inv\cdots f_r\inv\pi_1
\end{equation}
for some $r\le t$ and some reduced paths $\pi_1$ and $\pi_2$ in $\mX$. It follows that $\pi_1$ runs in $\mX_{r-1}$ and $\pi_2$ runs in $\mX_t$ (recall that $f_i\colon p_{i-1}\longrightarrow p_i$); it is sufficient to deal with the first 
case in~\eqref{eq:options}. The reduced path $\varsigma_e\colon {\alpha e}\longrightarrow{\omega e}$ in $\mZ$ then is given by
$\varsigma_e=\pi_1\rho\pi_2$ where $\rho$ is the reduced path ${p_{r-1}}\longrightarrow{p_t}$ in $\mY$. Recall that the reduced path ${p_{i-1}}\longrightarrow{p_i}$ in $\mY$ is denoted $\rho_i$. Then the path $\rho$ is the reduced form of the product $\rho_r\cdots \rho_t$, 
so that
\begin{equation}\label{eq: inequal I}
	[w_\mY(\rho_r\cdots\rho_t)]_I\le[w_\mY(\rho)]_I
\end{equation}
 in every $A$-generated inverse monoid $I$. For any $r\le i\le t$ we have $w_\mX(f_i)=w_\mC(f_i)$ and by construction of $\mY$, $[w_\mC(f_i)]_M\le[w_\mY(\rho_i)]_M$, hence $[w_\mX(f_i)]_M\le[w_\mY(\rho_i)]_M$. Using this and \eqref{eq: inequal I} we get
 \[[w_\mX(f_r\cdots f_t)]_M\le[w_\mY(\rho_r\cdots \rho_t)]_M\le[w_\mY(\rho)]_M=[w_\mZ(\rho)]_M.\]
Since $\pi_1$ and $\pi_2$ are reduced paths in $\mX$, $w_\mX(\pi_i)=w_\mZ(\pi_i)$ for $i=1,2$. Putting it all together we get
\begin{equation*}
\begin{aligned}
	{[w_\mZ(\varsigma_e)]_M} &= [w_\mZ(\pi_1)\cdot w_\mZ(\rho)\cdot w_\mZ(\pi_2)]_M
	=[w_\mX(\pi_1)\cdot w_\mY(\rho)\cdot w_\mX(\pi_2)]_M\\ &=[w_\mX(\pi_1)]_M[w_\mY(\rho)]_M[w_\mX(\pi_2)]_M\\
	&\ge [w_\mX(\pi_1)]_M[w_\mX(f_r\cdots f_t)]_M[w_\mX(\pi_2)]_M\\
	&=[w_\mX(\pi_1f_r\cdots f_t\pi_2)]_M = [w_\mX(\pi_e)]_M\ge [v(e)]_M,
\end{aligned}
\end{equation*}
where the last inequality follows from \eqref{eq: ineq for X} since $M\times G\twoheadrightarrow M$.

For $H$ to be critical for the pair $(M,\mE)$ one little detail is still missing: we need to modify $\mZ$ so that the diameter of the resulting tree is at most $n-1$. 
To this end
we proceed just as we did at the beginning of the proof: consider in $\mZ$ all maximal arcs spanned by simple paths of the form \eqref{eq:bad path Ash} with intermediate vertices of degree $2$ but none of it being 
some $p\in \mathrm{V}(\mE)$; replace every such arc by a single edge $e\colon q_0\longrightarrow q_t$ and define the label of each such edge by $w_\mT(e):=w_\mZ(e_1)\cdots w_\mZ(e_t)$; for every edge $e$ of $\mZ$ not involved in such an arc we just set $w_\mT(e):=w_\mZ(e)$. The resulting tree $\mT$ with the labelling $w_\mT$ then fulfils all conditions required for Definition~\ref{def: critical group} (critical group).
\end{proof}

\subsection{Ash implies Herwig--Lascar} \label{subsec: Ash===>HL}

Let $\sigma$ be a finite relational signature, and let us call an element $y$ of a  $\sigma$-structure $\bY$ \emph{relation-free} if it is not contained in any relational tuple $(y_1,\dots, y_r)$ (for any $r\ge 1$) of $\bY$ and let $\mathrm{rf}(\bY)$ denote the set of all relation-free elements of $\bY$. We define the \emph{weight} $\mathrm{w}(\bY)$ of $\bY$ as the number of relation-free elements of $\bY$ plus the number of the relational tuples of $\bY$; that is, $\mathrm{w}(\bY):=|\mathrm{rf}(\bY)|+\sum_{R\in \sigma}|R^\bY|$. Let $\mathbf{C} := \mathbf{Excl}(\Sigma)$
be the class of $\sigma$-structures defined by forbidden homomorphisms from a finite collection $\Sigma$ of finite $\sigma$-structures $\bT \in \Sigma$. 
In the proof given below that $\mathbf{C}$ is an EPPA-class, the weight of these forbidden structures $\bT$ will give a bound on the number of vertices of the graphs used in the application of Theorem~\ref{thm: ash}. Before we are ready for the proof we still have to discuss how to find a group suitable for the situation of Lemma~\ref{lem:X is induced}.
\begin{Lemma}\label{lem:X is induced holds} Let $M$ be a finite $A$-generated inverse monoid, $r\ge 1$ and let $\mE$ be the  graph consisting of two vertices $o\ne p$ 
linked by 
$r$ geometric edges.
Let $G$ be a finite $A$-generated group according to 
Theorem~\ref{thm: ash} with respect to $M$ and $\mE$. Then for every $r$-tuple $w_1,\dots, w_r$ of words for which $[w_1]_G=\cdots=[w_r]_G$ there exists a word $w$ such that $[w_1]_M,\dots, [w_r]_M\le [w]_M$. 
\end{Lemma}
\begin{proof} Let $k_1,\dots,k_r\colon o\longrightarrow p$ be the edges from $o$ to $p$. Label these edges according to $\ell(k_i^{\pm 1})=w_i^{\pm 1}$ for all $i$. 
The assumption that $[w_1]_G=\cdots=[w_r]_G$ implies that this 
labelling commutes over $G$.  
By Theorem~\ref{thm: ash},  there exist  words $v_i$ ($i=1,\dots, r$) such that (i) $[v_iv_j\inv]_F=1_F$ and (ii) $[w_i]_M=[v_i]_M$ for all $i,j$. From (i) it follows that $[v_i]_F=[v_j]_F$ for all $i,j$ so that all words $v_i$ have the same reduced form: $\mathrm{red}(v_i)=\mathrm{red}(v_j)=:w$. Hence $[w]_M\ge [v_i]_M=[w_i]_M$ for all $i$.
\end{proof}

We are now ready to derive the Theorem of Herwig--Lascar from Ash's Theorem.

\begin{proof}[Proof of Theorem~\ref{thm: Herwig--Lascar} by use of Theorem~\ref{thm: ash}] Let $\mathbf{C}$ be a class of \hyphenation{structu-res} $\sigma$-structures defined by finitely many forbidden homomorphisms of finite $\sigma$-structures, let $\bX=(X;(R^\bX)_{R\in \sigma})$ be a finite member of $\mathbf{C}$  equipped with some set $\til A$ of partial automorphisms, and let $\mP$ be the $A$-graph with vertex set $X$ encoding these partial mappings. For the set $\Sigma$ of forbidden structures $\bT$ let $n$ be the maximum of the weights of all members of $\Sigma$, that is,  $n= \max_{\bT\in \Sigma}\mathrm{w}(\bT)$. Now choose a finite $A$-generated group $G$ according to Theorem~\ref{thm: ash} for the inverse monoid $M=\cT(\mP)$ and for the (finite) graph $\mP$ consisting of the disjoint union of all (isomorphism types of) connected graphs having at most $n$ vertices and, between any two distinct vertices, at most $|X|$ many geometric edges. By Lemma~\ref{lem:X is induced holds} in combination with Lemma~\ref{lem:X is induced},  $\bX$ is an induced substructure of $\bX_\mP G$. We shall prove that, if $\bX_\mP G$ is not $\bT$-free for some $\bT\in \Sigma$, then neither is the free extension $\bX_\mP F$. 
	
Let $\bT\in \Sigma$ and let $\phi\colon \bT\to \bX_\mP G$ be a homomorphism. The weight of $\phi(\bT)$ as a weak substructure of $\bX_\mP G$ is at most $n$. So there exists a subset $J\subseteq G$ of size $n$ such that every relational tuple $(\phi(t_1),\dots,\phi(t_r))$ of $\phi(\bT)$ is completely contained in some single translate $Xg$ for some $g\in J$ (meaning that every entry $\phi(t_i)$ of this tuple is contained in that translate) and that $\phi(\bT)\subseteq \bigcup_{g\in J}Xg$. In other words, by the definition of the weight, the size $n$ of $J$ guarantees that (i) the members of every relational tuple of $\phi(\bT)$ are contained in at least one single translate $Xg$ of $X$ (one translate for each tuple) and (ii) the entire homomorphic image $\phi(\bT)$ is contained in the union $\bigcup_{g\in J}Xg$.\footnote{At this point we would not require to 
consider relation-free elements, since we could redefine the homomorphism $\phi$ to map the relation-free elements wherever we wish. The present pedantic approach is in a sense rather a preparation for the proof of Theorem~\ref{thm: HL strengthened}.} For convenience, we consider $X$ as a subset of $X_\mP G$ as well as of $X_\mP F$, in fact we assume that $X= X_\mP G\cap X_\mP F$. We want to establish a homomorphism $\phi(\bT)\to \bX_\mP F$. The idea is to find a function $\psi\colon \bigcup_{g\in J}Xg\to X_\mP F$ which is a homomorphism of $\sigma$-structures when restricted to $\phi(\bT)$.  More precisely, for each $g\in J$ we would like to choose some $f_g\in F$ and consider the mapping: \begin{equation}\label{eq:desired mapping}
\begin{array}{rcccl}   Xg & \longrightarrow &  X &\longrightarrow &Xf_g \subseteq X_\mP F\\
	z  & \longmapsto    & zg\inv& \longmapsto & (zg\inv)f_g \end{array}
\end{equation}
where $g\inv$ acts on $X_\mP G$ and maps $Xg$ to $X$, and $f_g$ acts on $X_\mP F$ and maps $X$ to $Xf_g$. Suppose we had a mapping $\psi\colon \bigcup_{g\in J}Xg\to X_\mP F$  whose restriction $\psi\upharpoonright Xg$ to any translate $Xg$ can be described by \eqref{eq:desired mapping}. Such a mapping would be compatible with the relational structure of $\phi(\bT)$, hence a homomorphism when restricted to $\phi(\bT)$. Indeed, any relational tuple $(\phi(t_1),\dots,\phi(t_r))\in R^{\bX_\mP G}$ (for some $R\in \sigma$) of $\phi(\bT)$ is contained in some translate $Xg$. The mapping described in \eqref{eq:desired mapping} certainly maps the elements $\phi(t_1),\dots,\phi(t_r)$ to elements of $Xf_g$, and 
the image tuple under $\psi$ belongs to $R^{\bX_\mP F}$ --- this follows from the definition of the $\sigma$-structures $\bX_\mP G$ and $\bX_\mP F$.

 The difficulty lies in the well-definedness of the map $\psi$ according to~\eqref{eq:desired mapping}.
Distinct translates $Xg$ and $Xh$ for $g,h\in J$ may overlap:  the choices $f_g, f_h\in F$ should satisfy
\[(zg\inv)f_g=(zh\inv)f_h\]
for all $g,h\in J$ for which $Xg\cap Xh\ne\varnothing$ and all $z\in Xg\cap Xh$. This would make such a mapping $\psi$ well-defined. The situation is depicted in Figure~\ref{fig:ashtohl}.
\begin{figure}[ht]
\begin{tikzpicture}[xscale=1.6,yscale=0.8]
	\usetikzlibrary{arrows.meta} 
	\filldraw(0,2)circle(1pt);
	\filldraw(0,1.5)circle(1pt);
	\filldraw(0,1)circle(1pt);
	\filldraw(0,-2)circle(1pt);
	\filldraw(0,-1.5)circle(1pt);
	\filldraw(0,-1)circle(1pt);
	\filldraw(-2,0.5)circle(1pt);
	\filldraw(-2,0)circle(1pt);
	\filldraw(-2,-0.5)circle(1pt);
	\filldraw(2,0.5)circle(1pt);
	\filldraw(2,0)circle(1pt);
	\filldraw(2,-0.5)circle(1pt);
	\draw(0,-2.2)[below]node{$y_i$};
	\draw(0,2.2)[above]node{$x_i$};
	\draw(-2.1,0)[left]node{$z_i$};
	\draw plot [smooth cycle] coordinates  {(0.5,2.3)(-0.5,2.3)(-0.5,-2.3)(0.5,-2.3)};
	\draw plot [smooth cycle] coordinates  {(-1.5,4.3)(-2.5,4.3)(-2.5,-0.3)(-1.5,-0.3)};
	\draw plot [smooth cycle] coordinates  {(2.5,4.3)(1.5,4.3)(1.5,-0.3)(2.5,-0.3)};
	\draw plot [smooth cycle] coordinates  {(-1.5,0.3)(-2.5,0.3)(-2.5,-4.3)(-1.5,-4.3)};
	\draw plot [smooth cycle] coordinates  {(2.5,0.3)(1.5,0.3)(1.5,-4.3)(2.5,-4.3)};
	\draw[-Latex,blue](0,2)--(-2,0.5);
	\draw[-Latex,blue](0,2)--(2,0.5);
	\draw[-Latex,blue](0,1.5)--(-2,0);
	\draw[-Latex,blue](0,1.5)--(2,0);
	\draw[-Latex,blue](0,1)--(-2,-0.5);
	\draw[-Latex,blue](0,1)--(2,-0.5);
	\draw[-Latex,red](0,-2)--(-2,-0.5);
	\draw[-Latex,red](0,-2)--(2,-0.5);
	\draw[-Latex,red](0,-1.5)--(-2,0);
	\draw[-Latex,red](0,-1.5)--(2,0);
	\draw[-Latex,red](0,-1)--(-2,0.5);
	\draw[-Latex,red](0,-1)--(2,0.5);
	\draw[blue](-2,-3)node{$Xg$};
	\draw[red](-2,3)node{$Xh$};
	\draw[blue](2,-3)node{$Xf_g$};
	\draw[red](2,3)node{$Xf_h$};
	\draw[red](-1,2.7)node{$\cdot h$};
	\draw[red](1,2.65)node{$\cdot f_h$};
	\draw[blue](-1,-2.7)node{$\cdot g$};
	\draw[blue](1,-2.65)node{$\cdot f_g$};
	\draw[>-{Latex[]Latex},red](-0.2,2.2)--(-2,4);
	\draw[>-{Latex[]Latex},red]	(0.2,2.2)--(2,4);
	\draw[>-{Latex[]Latex},blue](-0.2,-2.2)--(-2,-4);
	\draw[>-{Latex[]Latex},blue](0.2,-2.2)--(2,-4);
	\draw(0,0)node{$X$};
	\end{tikzpicture}	
\caption{Well-definedness of the mapping~\eqref{eq:desired mapping}}\label{fig:ashtohl}
\end{figure}

This is the point where we apply Ash's Theorem. Define a graph $\mJ$ with set of vertices $\mathrm{V}(\mJ)=J$ in which two vertices  $g$ and $h$ are connected by one or more edges if and only if $Xg\cap Xh\ne \varnothing$. In fact, there will be $|Xg\cap Xh|$ many edges from $g$ to $h$, and their inverses. 
Suppose that $Xg\cap Xh=\{z_1,\dots, z_r\}$; then there are $x_1,\dots, x_r, y_1\dots,y_r\in X$ such that $x_ig=z_i=y_ih$ for all $i$, hence $x_igh\inv =y_i$ for all $i$. According to Lemma~\ref{lem:extending G} there exist words $v_1,\dots,v_r\in \til A^*$ such that $x_igh\inv =y_i={x_i}\cdot[v_i]_M$ for all $i$ and $[v_i]_G=gh\inv$. In this situation we impose edges $e_1,\dots, e_r\colon g\longrightarrow h$ with labels $v_1,\dots, v_r$, respectively. The corresponding inverse edges $e_i\inv\colon h\longrightarrow g$ then have labels $v_i\inv$ and satisfy $y_ihg\inv = x_i = {y_i}\cdot[v_i\inv]_M$. 
Doing this for all pairs $g,h\in J$ for which $Xg\cap Xh\ne\varnothing$ 
yields a word labelling of $\mJ$ that commutes over $G$: 
 since for every edge $e\colon g\longrightarrow h$, the $G$-value of its label is $gh\inv$. Hence the label of any closed path has $G$-value $1_G$. According to Theorem~\ref{thm: ash} there exists an $M$-related relabelling that commutes over the free group $F$. For the edges $e_1,\dots, e_r\colon g\longrightarrow h$ and their labels $v_1,\dots, v_r$ (mentioned above) this means that there are  words $w_1,\dots, w_r\in \til A^*$ such that $[v_i]_M=[w_i]_M$ and the overall labelling by these words commutes over $F$. In particular, for fixed $g$ and $h$, the words $w_i$ (of the edges $e_i\colon g\longrightarrow h$) all have the same reduced form, say $w$, and $[v_i]_M=[w_i]_M\le [w]_M$ so that also ${x_i}\cdot[w]_M=y_i$ for all $i$. We relabel the graph once more by replacing every $w_i$ by this reduced word $w=:w^{(g,h)}$. Again this is done for all pairs $g, h$ for which $Xg\cap Xh\ne \varnothing$. 

Suppose first that $\mJ$ is connected.  Choose a base vertex $o$ and define, for any vertex $g$, an element $f_g$ by setting $f_g:=[p]_F$ where $p$ is the label of some path $g\longrightarrow o$; this labelling is meant with respect to the last relabelling where every edge $e\colon g\longrightarrow h$ is labelled $w^{(g,h)}$. Since the labelling commutes over $F$, the element $f_g$ is well-defined and does not depend on the choice of the path $g\longrightarrow o$. Now let $z\in Xg\cap Xh$ for some $g,h\in J$. Then there are $x,y\in X$ such that $xg=z=yh$, hence $x=zg\inv$, $y=zh\inv$ and $xgh\inv = y = {x}\cdot [w^{(g,h)}]_M$. From the definition of $f_{\underline{\ }}$ it follows that $f_g=[w^{(g,h)}]_Ff_h$. Altogether we arrive at
\[(zg\inv)f_g=xf_g=x[w^{(g,h)}]_Ff_h=({x}\cdot[w^{(g,h)}]_M)f_h=yf_h=(zh\inv) f_h,\]
which shows that the mapping $\psi\colon\bigcup_{g\in J}Xg\to X_\mP F$ 
according to~\eqref{eq:desired mapping} is well-defined.
If the graph $\mJ$ is not connected, then the same argument can be applied to each connected component of $\mJ$ separately. It follows that the free extension $\bX_\mP F$ is not $\bT$-free, and the theorem is proved.
\end{proof}

\subsection{Herwig--Lascar implies Ash}\label{subsec:HL===>Ash} 
In this subsection we give a direct proof of the implication Herwig--Lascar $\Longrightarrow$ Ash, which seems to be interesting in its own right. Essentially it is a generalisation of the arguments in Section~\ref{subsec:HL>RZ} from cycle graphs to arbitrary graphs; the proof of Proposition~\ref{prop:HL>RZ} can be seen as a warm up to the proofs of Lemmas~\ref{lem:notChasEPPA} and~\ref{lem:EPPAhasnotC}. It also provides an essential ingredient for the understanding of the discussion of 
the group-theoretic version(s) of the Herwig--Lascar Theorem
in Section~\ref{sec: group-theoretic version}. Otherwise this section is independent of the remainder of the paper.

Let $\mE=(V\cup K;\alpha,\omega,\inv)$ be a graph; a \emph{coset labelling of $\mE$} is a mapping $\ell\colon K\to 2^F$, $k\mapsto \ell(k):=H_kg_k$, with $H_k$ a finitely generated subgroup of $F$ and $g_k\in F$ such that $\ell(k\inv)=\ell(k)\inv$, that is, $H_{k\inv} g_{k\inv}=(H_kg_k)\inv =({}^{g_k\inv}\!\!H_k)g_k\inv$, which is the set of all inverses of the elements of $H_kg_k$. 
A coset labelling $\ell$ of $\mE$ \emph{commutes over the $A$-generated group $G$} (is \emph{$G$-commuting} for short) if there exists 
a $G$-commuting word labelling $l$ 
such that $l(k)\in \ell(k)$ for every $k\in K$
(in this subsection we consider $F$ to be the set of all reduced words in $\til A^*$ and assume that every word labelling takes values in $F$). In a sense, a coset labelling is $G$-commuting if it \emph{covers} (or \emph{contains}) a $G$-commuting word labelling.
The essential achievement in this subsection will be to prove that  Proposition~\ref{prop:HL>Ash} below is  a consequence of the Herwig--Lascar Theorem. 

\begin{Prop}\label{prop:HL>Ash} Let $\mE$ be a finite graph and let
  $\ell$ be a coset labelling of $\mE$. If $\ell$ does not commute over $F$, then $\ell$ does not commute over $G$ for some finite $A$-generated group $G$.
\end{Prop}

Intuitively, the statement means that the non-commutativity with respect to $F$ of a coset labelling of a finite graph  can be witnessed (or recognised) by some finite group $G$; or, to say it the other way round: the free group is well approximable by finite groups with respect to commutativity of coset labellings of finite graphs. 

We first note that Proposition~\ref{prop:HL>Ash} is a direct generalisation of Proposition~\ref{prop:HL>RZ} from cycle graphs to arbitrary graphs. Indeed, let $H_1,\dots,H_n$ be finitely generated subgroups of $F$ and $g_1,\dots, g_n\in F$. The sequence of cosets $H_1g_1,\dots, H_ng_n$ can be seen as a coset labelling of the cycle graph of length $n$. To say that $1\notin H_1g_1\cdots H_ng_n$ is equivalent to saying that this labelling is not $F$-commuting; to say that $1\notin NH_1g_1\cdots NH_ng_n$ is equivalent to saying that the original coset labelling is not $G$-commuting for the group $G=F/N$.

Next we note that it suffices to prove Proposition~\ref{prop:HL>Ash} (i.e.~derive it from Theorem~\ref{thm: Herwig--Lascar}) for finite connected graphs without loop edges. 
Let $\mE$ be a finite connected graph and let $\mD$ be the graph obtained from $\mE$ by removing all loop edges (if any).  Consider a coset labelling $\ell_\mE$ that is not $F$-commuting. If the restriction of $\ell_\mE$ to $\mD$ is also not $F$-commuting we are done. Suppose we have a labelling $\ell_\mE$ of $\mE$ that is not $F$-commuting while its restriction to $\mD$ is $F$-commuting. Then there exists  a loop edge $e$ whose label $\ell_\mE(e)$ is a proper coset $Hg\ne H$ (that is, $1_F\notin Hg$). There is a finite $A$-generated group $G$ for which $[g]_G\ne 1_G$: $G$ can be chosen to be the transition group of any finite completion $\ol{\mS_{Hg}}$ of the Stallings graph $\mS_{Hg}$ (this is Stallings' proof~\cite{stallings} of Hall's Theorem~\cite{hall}).   Hence, given a labelling $\ell_\mE$ of $\mE$, if there exists a finite group $G$ detecting that its restriction to $\mD$ is not $F$-commuting, then there also exists a finite group $H$ detecting that $\ell_\mE$ itself is not $F$-commuting. Finally, a graph $\mE$ satisfies the condition formulated in Proposition~\ref{prop:HL>Ash} if and only every connected component of $\mE$ does.

The intuitive idea of the proof roughly is as follows. Suppose we are given a graph $\mE$ and a coset labelling $\ell$ of $\mE$. Let $e=\{k,k\inv\}$ be a geometric edge with incident vertices $u=\alpha k$ and $v=\omega k$. The idea is to replace $e$ with the plain Schreier graph $\bm\mS_e:=\bm\mS^\circ_{H_kg_k}$ where $H_kg_k=\ell(k)$, and to encode the entire labelling $\ell$ on the disjoint union $\bm\mS:=\bigsqcup \bm\mS_e$ where the union runs over all geometric edges $e$ of $\mE$. 
In order that the latter union encodes the labelling completely we need to mark the two distinguished vertices of $\bm\mS_{H_kg_k}$ within $\bm\mS_e$ 
so as to associate the correct coset $H_kg_k$ with the edge $e$ 
when traversed in the direction from $u$ to $v$.
This can be achieved by \textsl{colouring} the initial vertex $\iota$ of $\bm\mS_{H_kg_k}$ by colour $u$ and the terminal  vertex $\tau$ of $\bm\mS_{H_kg_k}$  by colour $v$. When $e$ is traversed in the opposite direction this colouring already provides the correct coset, namely $({}^{g_k\inv}\! H){g_k\inv}$. Note that
it may happen that a distinguished vertex carries both colours $u$ and $v$ (this happens exactly if the corresponding coset $H_kg_k$ coincides with the group $H_k$.) So far, 
on each connected component $\bm\mS_e$ of $\bm\mS$, two vertices $\iota$ and $\tau$ of $\bm\mS_e$ (possibly identical) 
are --- thanks to the colouring --- associated 
with the two vertices $u$ and $v$ of $\mE$ (namely those incident with the geometric edge $e$). From point of view of the vertices of $\mE$ we have, for every vertex $u$ of $\mE$,  $u$-coloured vertices in $\bm\mS$, one in each component $\bm\mS_e$, for every geometric edge $e$ incident with $u$. Now recall that the free group $F$ acts on $\bm\mS$ leaving invariant every component $\bm\mS_e$. For every vertex $v$ of $\mE$ we consider the collection of all $v$-coloured distinguished vertices  as being ``rigidly entangled'' forming a ``$v$-block'' that is moved ``as a whole'' by the action of $F$ (although they belong to pairwise distinct components $\bm\mS_e$). Later, in the formal proof, these $v$-blocks will be realised as relational tuples whose entries are comprised of the $v$-coloured vertices (the sets of these will be denoted by $M_v$). 

The idea now is that the original labelling $\ell$ commutes over $F$ if and only if the action of $F$ moves the $v$-blocks in such a way that, for every geometric edge $e$, the two distinguished vertices in $\bm\mS_e$ are mapped to a common vertex of $\bm\mS_e$. (A trivial special case of this is when all labels $\ell(k)$ are subgroups: then all pairs of  distinguished vertices in every component $\bm\mS_e$ coincide from the very beginning and the labelling trivially commutes since every label contains the identity element $1_F$.) More precisely, for every vertex $v$ of $\mE$ there is an element  $w_v\in F$ such that the following holds
for the collection $(w_v\colon v\in V)$: 
for every geometric edge $e$ of $\mE$ with incident vertices $u$ and $v$ the image under $w_u$ of the $u$-coloured vertex of $\bm\mS_e$ coincides with the image under $w_v$ of the $v$-coloured vertex. It then follows that $w_uw_v\inv$ labels a path from the $u$-vertex to the $v$-vertex in $\bm\mS_e$. Consequently, the reduced form $\mathrm{red}(w_uw_v\inv)$ belongs to the coset label of the edge $e$ (the label when the edge is traversed from $u$ to $v$). It follows that, from any coset label we may select an element such that for any cycle in the graph $\mE$, the product over these selected elements is the identity element $1_F$. For the special case of a cycle graph $\mE$ this is the essence of what happens in the proof of Proposition~\ref{prop:HL>RZ}, see Figure~\ref{fig:structureS}. This idea will now be formalised and translated into the language of relational structures.

In the following we shall assume that the graph $\mE=(V\cup K;\alpha,\omega,\inv)$ in question is  finite, connected and without loop edges. We introduce a relational signature $\sigma_V=\{R_v\colon v\in V\}$ consisting of relational symbols $R_v$, indexed by the vertices of $\mE$, the arity of $R_v$ being $r_v$, the degree of the vertex $v$. Let $\bC$ and $\bC^\star$ be two $\sigma_V$-structures  defined as follows. The base set of $\bC$ is $K$, the set of edges of $\mE$; for $v\in V$ we let  $K_v:=\{k\in K\colon \omega k=v\}$. 
For every $v\in V$ we choose an enumeration $k_1,\dots, k_{r_v}$ of $K_v$ and set $R_v^\bC:=\{(k_1,\dots,k_{r_v})\}$.  This defines a $\sigma_V$-structure $\bC=(K;(R_v^\bC)_{v\in V})$.
We call  $\bC$ a \emph{$\sigma_V$-structure associated with the graph $\mE$}. 

In the context of $\sigma_V$-structures, for any equivalence relation $\sim$ on the base set $K$, the quotient structure $\bC/{\sim}$ is defined on the quotient set $K/{\sim}$ as follows: a tuple $(\bar{k_1}\dots, \bar{k_{r}})$ of the quotient set $K/{\sim}$ is  in $R_v^{\bC/{\sim}}$  if and only if for all $i$ there are $q_i\in \bar{k_i}$ such that the tuple $(q_1,\dots, q_{r})$ is in $R_v^\bC$.
Finally, given a $\sigma_V$-structure $\bC$ associated with $\mE$ we define the $\sigma_V$-structure $\bC^\star$ by factoring the base set $K$ of $\bC$ by the equivalence relation $\equiv$  which identifies every edge $k$ with its inverse $k\inv$:
\[\bC^\star:=\bC/{\equiv}\]
where for $k_1,k_2\in K$, \[k_1\equiv k_2\Longleftrightarrow k_1=k_2\mbox{ or }k_1=k_2\inv. \]
We call $\bC^\star$ a \emph{consolidated} $\sigma_V$-structure  associated with the graph $\mE$. In addition, any $\sigma_V$-structure $\bC/{\sim}$ (including $\bC$ itself) obtained by factoring the base set $K$ of $\bC$ by an equivalence relation $\sim$ contained in $\equiv$ (that is by identifying \textsl{some}, but not necessarily \textsl{all}, edges $k$ with their inverses $k\inv$) is called an \emph{approximation of $\bC^\star$}. 

Prima facie, every structure $\bC$ associated with a graph $\mE$ depends on the choice of the enumeration of each set $K_v$, but it is easy to see that any two different such structures $\bC_1$ and $\bC_2$ are isomorphic as $\sigma_V$-structures. However, this does not imply that the corresponding consolidated structures $\bC_1^\star$ and $\bC_2^\star$ are also isomorphic. The reason for this is that an isomorphism $\bC_1\to \bC_2$ not necessarily respects the equivalence relation $\equiv$ on $K$. An example $\mE$ for such a behaviour is the cycle graph  on two vertices (that is, $|V|=2$ and $|K|=4$). In this case there exist four possible associated structures $\bC_i$, all of which are isomorphic as $\sigma_V$-structures, and also four possible consolidated structures $\bC^\star_i$ which belong to two different isomorphism classes of $\sigma_V$-structures. This peculiar behaviour is irrelevant for our purpose.

Recall the concepts of the Stallings graph $\mS_{Hg}$ and the Schreier graph $\bm{\mS}_{Hg}$ of a coset $Hg$ 
(see Section~\ref{subsec:stallings}). Both are $2$-pointed $A$-graphs consisting of the \emph{plain} $A$-graph $\mS^\circ_{Hg}$ (respectively $\bm\mS^\circ_{Hg}$) (plain $A$-graph means that there are no distinguished vertices) and the distinguished (initial and terminal) vertices $\iota$ and $\tau$. 
We use the suggestive notation
\[\mS_{Hg}=(\iota,\mS^\circ_{Hg },\tau)\]
and set $\iota(\mS_{Hg}):=\iota$ and $\tau(\mS_{Hg}):=\tau$. 
For the ``inverse'' coset $({}^{g\inv\!\!}H){g\inv}$ we have  $\mS^\circ_{Hg }=\mS^\circ_{({}^{g\inv\!\!}H){g\inv}}$ so that
\[\mS_{({}^{g\inv\!\!}H){g\inv}}=(\tau,\mS^\circ_{Hg},\iota),\]
which means that 
the plain $A$-graph $\mS^\circ_{Hg}$ remains unchanged and just 
the two distinguished vertices $\iota$ and $\tau$ are swapped. This is expressed by the formula
\[\iota(\mS_{Hg})=\tau(\mS_{({}^{g\inv\!\!}H){g\inv}})\mbox{ and }\iota(\mS_{({}^{g\inv\!\!}H){g\inv}})=\tau(\mS_{Hg}),\] 
and the same holds for the corresponding Schreier graphs $\bm\mS_{Hg}$ and $\bm\mS_{({}^{g\inv\!\!}H){g\inv}}$.

Next consider a coset labelling $\ell$ of $\mE$. For every edge $k$ of $\mE$ we have $\ell(k)=H_kg_k$ for some finitely generated subgroup $H_k$ of $F$ and some $g_k\in F$ and such that $\ell(k\inv)=H_{k\inv}g_{k\inv}=({}^{g_k\inv\!\!}H_k)g_k\inv$, that is, $g_{k\inv}=g_k\inv$ and $H_{k\inv}={}^{g_k\inv\!\!}H_k$. 
Recall the notation for geometric edges $k^\bullet:=\{k,k\inv\}$ and $K^\bullet:=\{k^\bullet\colon k\in K\}$. For     $e={k^\bullet}$ let \[\mS_e:=\mS^\circ_{H_kg_k}=\mS^\circ_{({}^{g_k\inv\!\!}H){g_k\inv}}\] be the plain Stallings graph of $H_kg_k$ --- that is, the Stallings graph with the distinguished vertices ignored, which is also the plain Stallings graph of $H_{k\inv}g_{k\inv}=({}^{g_k\inv\!\!}H_k)g_k\inv$ --- and let $\mS:=\bigsqcup_{e\in {K^\bullet}} \mS_e$. We define a mapping $\varpi\colon K\to \mathrm{V}(\mS)$ as follows: for $e={k^\bullet}=\{k,k\inv\}$  let $\varpi(k):=\tau(\mS_{H_kg_k})$ (the terminal vertex of $\mS_{H_kg_k}=(\iota,\mS_e,\tau)$, considered as a vertex of $\mS_e\subseteq \mS$) and $\varpi(k\inv):=\iota(\mS_{H_kg_k})$ (the initial vertex of $\mS_{H_kg_k}=(\iota,\mS_e,\tau)$, considered as a vertex of $\mS_e\subseteq \mS$). The discussion above about plain Stallings graphs implies that the mapping $\varpi$ is well-defined in the sense that $\varpi(k\inv)=\tau(\mS_{({}^{g_k\inv\!\!}H_k)g_k\inv})=\tau(\mS_{H_{k\inv}g_{k\inv}})$, the terminal vertex of the Stallings graph $\mS_{H_{k\inv}g_{k\inv}}$. 

If $k_1$ and $k_2$ are edges that belong to different geometric edges, that is, ${k_1^\bullet}\ne {k_2^\bullet}$, then $\varpi(k_1)\ne \varpi(k_2)$ since $\varpi(k_1)\in \mathrm{V}(\mS_{k_1^\bullet})$ while $\varpi(k_2)\in \mathrm{V}(\mS_{{k_2^\bullet}})$ (which are disjoint). In particular, the restriction $\varpi\upharpoonright K_v$ is injective for every $v\in V$. For later use we set $\varpi_v:=\varpi\upharpoonright K_v$ and let $M_v:=\varpi(K_v)=\varpi_v(K_v)$ and note that $\varpi_v$ is a bijection $K_v\to M_v$ for every vertex $v$. 
Only in case when $k_1$ and $k_2$ belong to the same geometric edge ${k_1^\bullet}={k_2^\bullet}$ the values $\varpi(k_1)$ and $\varpi(k_2)$ may coincide. This happens exactly if  the initial and terminal vertices of $\mS_{\ell({k_1})}$ (and hence also of $\mS_{\ell({k_2})}$) coincide, which happens if and only if the coset $\ell(k_i)$ is actually a group.

Now we consider Schreier graphs; every Stallings graph $\mS_{Hg}$ is naturally embedded in the Schreier graph $\bm\mS_{Hg}$ which, by definition, is the free completion $\mS_{Hg}\mathrm{F}$ of $\mS_{Hg}$. For $e=\{k,k\inv\}$ we set $\bm\mS_e:=\bm\mS^\circ_{H_kg_k}$ and $\bm\mS:=\bigsqcup_{e\in {K^\bullet}}\bm\mS_e$. There are natural inclusions $\mS_e\subseteq \bm\mS_e$ for every $e\in {K^\bullet}$ and hence also $\mS\subseteq \bm\mS$. In particular, the mapping $\varpi$ can be interpreted as a mapping $K\to \mathrm{V}(\bm\mS)$. Recall that the free group $F$ acts by bijections on every set $\mathrm{V}(\bm\mS_e)$ and hence on $\mathrm{V}(\bm\mS)$ on the right via $\rho\mapsto \rho w$, for every $\rho\in\mathrm{V}(\bm\mS)$ and every $w\in F$, and the set $\mathrm{V}(\bm\mS_e)$ is invariant under this action for each $e\in {K^\bullet}$.

As in Section~\ref{subsec:HL>RZ} we define a $\sigma_V$-structure on $\mathrm{V}(\bm\mS)$, that is, for every $v\in V$ we define a relation $R_v$ of arity $r_v$ on $\mathrm{V}(\bm\mS)$ as follows. First choose some $\sigma_V$-structure $\bC$ associated with $\mE$ --- in particular, for every $v\in V$, choose and fix an enumeration $k_1,\dots, k_{r_v}$ for the elements of the set $K_v$; then transfer this enumeration to the corresponding set $M_v$ via the bijection $\varpi_v$:
\begin{equation}\label{eq:ordering of M_v}
\varpi_v k_1,\dots,\varpi_v k_{r_v};    
\end{equation} then set
\begin{equation}\label{eq:initialisation}
(\varpi_v k_1,\dots,\varpi_v k_{r_v})\in R_v
\end{equation}
and add to this whatever is required so that the aforementioned action $\rho\mapsto \rho w$ by $F$  becomes an automorphism of the $\sigma_V$-structure. More precisely
\begin{equation}\label{eq:extend structure}
R_v=\{(\varpi_v k_1,\dots,\varpi_v k_{r_v})w\colon w\in F\},
\end{equation}
where the application of $w$ to an $r_v$-tuple is understood componentwise.
In this way we get a $\sigma_V$-structure $\bS=(\mathrm{S};(R_v)_{v\in V})$ on the base set $\mathrm{S}=\mathrm{V}(\bm\mS)$ with  relations $R_v$ of arity $r_v$ and the free group $F$ acts on this structure by automorphisms on the right. The idea is to initialise a $\sigma_V$-structure  on the subset $\bigcup_{v\in V}M_v$ \eqref{eq:initialisation} and then extend the structure to all of $\mathrm{V}(\bm\mS)$ via the action of $F$ \eqref{eq:extend structure}.
The mapping $\varpi\colon K\to \mathrm{V}(\bm\mS)$ is a homomorphism  of $\sigma_V$-structures $\bC\to \bS$. Its image $\varpi(K)=\bigcup_{v\in V} M_v$, as a (weak) substructure of $\bS$, forms an approximation of $\bC^\star$ since for $k_1,k_2\in K$, $\varpi(k_1)=\varpi(k_2)$ can happen only if $k_1\equiv k_2$, as already mentioned.
We state and prove the essential assertion about $\bS$. 
\begin{Lemma}\label{lem:notChasEPPA}
If there is a homomorphism $\phi\colon\bC^\star\to \bS$, then the coset labelling $\ell$  of $\mE$  commutes over $F$.
\end{Lemma}

\begin{proof}
In the following, the interpretations of the symbols $R_v$ in $\bS$ are written without superscript as in the previous text, while those in $\bC^\star$ are denoted as $R_v^{\bC^\star}$.
For every $v\in V$ we have $R_v=\{(\varpi_v k_1,\dots,\varpi_v k_{r_v})w\colon w\in F\}$.
Let $\phi\colon \bC^\star\to \bS$ be a homomorphism of $\sigma_V$-structures. Consider two vertices $u,v\in V$ connected by an edge in $\mE$, say $k\colon u\longrightarrow v$ and $k\inv\colon v\longrightarrow u$. Let $K_v=\{k_1,\dots,k_{r_v}\}$ and $K_u=\{l_1,\dots,l_{r_u}\}$ with $\{({k_1^\bullet},\dots,{k_{r_v}^\bullet})\}= R_v^{\bC^\star}$ and $\{({l_1^\bullet},\dots,{l_{r_u}^\bullet})\}=R_u^{\bC^\star}$. By the definition of $R_v$ and $R_u$, respectively, there exist $w_v,w_u\in F$ such that
\[\phi({k_1^\bullet},\dots,{k_{r_v}^\bullet})=(\varpi_v k_1,\dots,\varpi_v k_{r_v})w_v\]
and
\[\phi({l_1^\bullet},\dots,{l_{r_u}^\bullet})=(\varpi_u l_1,\dots,\varpi_u l_{r_u})w_u,\]
where the application of $\phi$ and $\cdot w_v$ respectively $\cdot w_u$ are understood componentwise. Note that the ordering of the entries in the respective tuples comes from~\eqref{eq:ordering of M_v}. 
By construction, there are indices $i,j$ such that $k=k_i$ and $k\inv = l_j$. But then ${k_i^\bullet}={l_j^\bullet}$ which implies
\begin{equation}\label{eq:equal in V(S)}
(\varpi k)w_v=(\varpi_vk_i)w_v=\phi({k_i^\bullet})=\phi({l_j^\bullet})=(\varpi_ul_j)w_u =(\varpi k\inv)w_u.
\end{equation}
The equality~\ref{eq:equal in V(S)} takes place in $\mathrm{V}(\bm\mS_{{k^\bullet}})$ where $\bm\mS_{{k^\bullet}}$ is the plain graph of the Schreiergraph $\bm\mS_{H_kg_k}$. By definition of the mapping $\varpi$ we have $\varpi(k)=\tau:=\tau(\mS_{H_kg_k})=\tau(\bm\mS_{H_kg_k})$ and $\varpi(k\inv)=\iota:=\iota(\mS_{H_kg_k})=\iota(\bm\mS_{H_kg_k})$. Equality \eqref{eq:equal in V(S)} then just says that $\iota w_u=\tau w_v$ so that $\tau=\iota w_uw_v\inv$. The latter equality holds in $\bm\mS_{H_kg_k}$, which implies that $\tau=\iota\cdot\mathrm{red}(w_uw_v\inv)$ holds in $\mS_{H_kg_k}$, which in turn implies that $\mathrm{red}(w_uw_v\inv)\in H_kg_k$. Observe that $u=\alpha k$ and $v=\omega k$. 
 So we have shown  that for every edge $k\colon \alpha k\longrightarrow \omega k$, the word $w(k):=\mathrm{red}(w_{\alpha k}w_{\omega k}\inv)$ belongs to $\ell(k)$, and it is clear that the word labelling $k\mapsto w(k)=\mathrm{red}(w_{\alpha k}w_{\omega k}\inv)$ commutes over $F$.
\end{proof}
Next consider again, for $e\in {K^\bullet}$, the plain Stallings and Schreier 
graphs $\mS_e\subseteq \bm\mS_e$ and their disjoint unions 
\[\mS:=\bigsqcup_{e\in {K^\bullet}}\mS_e\subseteq \bigsqcup_{e\in {K^\bullet}}\bm\mS_e=\bm\mS.\]
Then $\mS$ is a finite subgraph of $\bm\mS$ and we let $\mathfrak{V}$ be the substructure of $\bS$ induced on the set $\mathrm{V}=\mathrm{V}(\mS)$. The $A$-labelling of $\mS$ provides partial automorphisms of the structure $\mathfrak{V}$. Recall that the initialisation \eqref{eq:initialisation} of the $\sigma_V$-structure $\bS$ on $\mathrm{V}(\bm\mS)$ actually takes place in $\mathrm{V}(\mS)$. 
This guarantees that an approximation of $\bS^\star$ is sitting in $\bV$.
Let $\ol{\mathfrak{V}}$ be an extension (on base set $\ol{\mathrm{V}}\supseteq \mathrm{V}$) of $\mathfrak{V}$ on which all partial automorphisms of $\mathfrak{V}$ given by $\til A$ extend to total automorphisms. Let $G$ be the corresponding $A$-generated group of automorphisms of $\ol{\mathfrak{V}}$. 
The next result formulates and proves the main assertion about $G$. {Recall that $\ell$ is a coset labelling of the graph $\mE=(V\cup K;\alpha,\omega,\inv)$, that is, $\ell$ is a map $K\to 2^F$ such that $\ell(k)=H_kg_k$ for $H_k$ a finitely generated subgroup of $F$ and $g_k\in F$. Also note that $V$ denotes the set of vertices of $\mE$ while $\mathrm{V}$ 
stands for the base set of $\mathfrak{V}$ which, in fact, is the set of vertices of $\mS$.} 
\begin{Lemma}\label{lem:EPPAhasnotC}
If the coset labelling $\ell$ commutes over $G$, then there exists a homomorphism $\bC^\star\to\ol{\mathfrak{V}}$. 
\end{Lemma} 
\begin{proof}
For every edge $k\in K$ choose a word $w(k)\in \ell(k)$ such that the word labelling $w\colon K\to F$ commutes over $G$. Let
$k\in K$ and $x\in H_kg_k$; then in $\mS_{H_kg_k}$ we have $\iota\cdot x =\tau$ where $\iota$ and $\tau$ are the initial and terminal vertices of $\mS_{H_kg_k}$, respectively. Hence in $\ol{\bV}$ we have $\iota[x]_G=\tau$. This implies that, for any $k\in K$, $\varpi(k\inv)[w(k)]_G=\varpi(k)$ also holds in $\ol{\bV}$. Now fix a vertex $o\in V$, and for every $v\in V$ choose a path $\pi_v\colon v\longrightarrow o$ and let $w(\pi_v)$ be the label of that path. By assumption, the $G$-value $g_v:= [w(\pi_v)]_G$ does not depend on the path $\pi_v$ but only on the vertex $v$.

 We modify the mapping $\varpi\colon K\to \mathrm{V}$  to a mapping $\mu\colon K\to \ol{\mathrm{V}}$ by setting, for $k\in K_{v}$, $\mu(k):= \varpi(k)g_{v}$ (note that $v=\omega k$). Clearly, $\mu$ is a homomorphism $\bC\to \ol{\bV}$ of $\sigma_V$-structures. Again let $k\in K$; the definition of
$v\mapsto g_v$ 
implies $g_{\alpha k}=[w(k)]_Gg_{\omega k}$ and we get
\begin{align*}
\mu(k\inv)&=\varpi(k\inv)g_{\omega k\inv} =\varpi(k\inv)g_{\alpha k}\\ &=\varpi(k\inv)[w(k)]_Gg_{\omega k}=\varpi(k)g_{\omega k}=\mu(k).
\end{align*}
It follows that ${\equiv} \subseteq \mathrm{ker}\mu$ and the homomorphism $\mu\colon \bC\to \ol{\bV}$ factors through $\bC^\star$. 
\end{proof}

The proof of Proposition~\ref{prop:HL>Ash} is now immediate.
\begin{proof}[Proof of Proposition~\ref{prop:HL>Ash} by use of Theorem~\ref{thm: Herwig--Lascar}] Let the finite graph $\mE=(V\cup K;\alpha,\omega,\inv)$ be connected and without loop edges, and let $\ell$ be a coset labelling of $\mE$ which is not $F$-commuting. 
 We need to show that $\ell$ is not $G$-commuting 
  for some suitable finite  $A$-generated group $G$. 

According to the discussion after Proposition~\ref{prop:HL>Ash} it is sufficient to prove the statement for graphs without loop edges. From the graph $\mE$ we obtain the relational signature $\sigma_V$ and 
an associated $\sigma_V$-structure $\bC$ and its consolidated structure $\bC^\star$. Let $\mathbf{C}$ be the class 
$\mathbf{Excl}( \bC^\star )$
of all $\sigma_V$-structures $\bT$ that do not admit a homomorphism $\bC^\star\to \bT$.  Now consider the structure $\bS$; by Lemma~\ref{lem:notChasEPPA} there is no homomorphism $\bC^\star\to \bS$, that is, $\bS\in \mathbf C$. Let $\bV$ be the (finite) substructure of $\bS$ induced on the set $\mathrm{V}=\mathrm{V}(\mS)$. 
 Viewing the partial automorphisms of $\bV$ that are induced by the $A$-graph structure of $\bS$
as an extension problem $(\bV,\til A)$, we have $\bS$ as a solution in $\mathbf{C}$.
By Theorem~\ref{thm: Herwig--Lascar} there exists a finite solution $\ol{\bV}$   in $\mathbf{C}$. Let $G$ be the $A$-generated group generated by the corresponding automorphisms.  Since $\ol{\bV}\in \mathbf{C}$ there is no homomorphism $\bC^\star\to \ol{\bV}$.  By Lemma~\ref{lem:EPPAhasnotC} then, the labelling $\ell$ is not $G$-commuting.
\end{proof}	
The argument that Theorem~\ref{thm: ash} follows from Proposition~\ref{prop:HL>Ash} is now similar to the second part of the proof of Theorem~\ref{thm:ash for cycles}. There we looked at all $n$-tuples of elements of the inverse monoid $M$. An $n$-tuple may be seen as a labelling of a cycle graph of length $n$. So let $\mE=(V\cup K;\alpha,\omega,\inv\!)$ be a finite connected graph. An $M$-labelling of $\mE$ is a mapping $l\colon K\to M$ such that $l(k\inv)=l(k)\inv$ for all $k\in K$. Every such labelling gives rise to a coset labelling $\ell$ by letting $\ell(k)=\varphi(l(k))$ where $\varphi$ is defined in \eqref{eq:canonrelmor}. 
By Proposition~\ref{prop:canonicalrelmorMto F} every $\varphi(l(k))$ is a right coset $H_kg_k$ for some finitely generated subgroup $H_k$ of $F$ and such that $\varphi(l(k\inv))=(\varphi(l(k)))\inv$. Some of these coset labellings may be $F$-commuting. We take those labellings $\ell$ that are not $F$-commuting. For every such there exists a finite group $G_\ell$ such that $\ell$ is not $G_\ell$-commuting. 
Then the $A$-generated direct product $G$ over all such groups $G_\ell$ 
yields 
a finite group $G$ required by Theorem~\ref{thm: ash}.

Finally, we argue similarly as after the proof of Theorem~\ref{thm:ash for cycles} and show that, conversely, (Ash's) Theorem~\ref{thm: ash} implies Proposition~\ref{prop:HL>Ash}. Let $\mE=(V\cup K;\alpha,\omega,\inv)$ be a finite graph and let $\ell\colon K\to 2^F$, $k\mapsto H_kg_k$ be a coset labelling not commuting over $F$. For each geometric edge $e=\{k,k\inv\}$ of $\mE$ let $\mS_e$ be the plain Stallings graph of $H_kg_k$ and let $\mS=\bigsqcup_{e\in K^\bullet}\mS_e$. Let $M=\cT(\mS)$  and let $G$ be a finite $A$-generated group as in Theorem~\ref{thm: ash} with respect to $\mE$ and $M$. We argue that $\ell$ is not $G$-commuting. Suppose, towards a contradiction, this were not the case and $\ell$ 
were $G$-commuting: there exists a $G$-commuting word labelling $l\colon K\to F$ covered by $\ell$. According to Theorem~\ref{thm: ash} there exists an $M$-related and $F$-commuting relabelling $l'$ of $l$. Let $k\in K$ and let $(\iota,\mS_{k^\bullet},\tau)$ be the Stallings graph of $H_kg_k$; since $l(k)\in \ell(k)=H_kg_k$ we have $\iota\cdot[l(k)]_M=\tau$. Since $l$ and $l'$ are $M$-related we also have $\iota\cdot[l'(k)]_M=\tau$, which just means that $l'(k)\in H_kg_k=\ell(k)$. Since $l'$ is $F$-commuting this contradicts the assumption that $\ell$ is not $F$-commuting.

\subsection{Group-theoretic versions of the Herwig--Lascar Theorem}
\label{sec: group-theoretic version} As mentioned in the introduction there is also a ``group-theoretic version'' of the Herwig--Lascar Theorem (Theorem 3.3 in~\cite{HL}). There are  even several variants of this version which we discuss briefly in this subsection; they 
 deal with (non)solvability in the free group $F$ of systems of equations of a certain type. Proposition~\ref{prop:HL>Ash} immediately gives us one of these versions, namely Proposition~\ref{prop:HLgroupformulationABO} below. {We present two further versions: Proposition~\ref{prop:ashalmeida} and Proposition~\ref{prop:HLgroupformulation}, the latter being the original group-theoretic version of the Herwig--Lascar Theorem mentioned in the introduction.}

 However, in order to line up with the literature (Propositions ~\ref{prop:ashalmeida} and~\ref{prop:HLgroupformulation} below) we  have to modify the statement of Proposition~\ref{prop:HL>Ash} a bit: instead of labellings $\ell\colon k\mapsto H_kg_k$ of the edges of a graph $\mE$ by right cosets $Hg$ we now consider labellings by left cosets $gH$. These labellings also have to satisfy $\ell(k\inv)=(\ell(k))\inv$ which now reads as $g_{k\inv}=g_k\inv$ and $H_{k\inv}= H_k^{g_k\inv}$ where $H^{g\inv}:=gHg\inv$, altogether, $\ell(k\inv)=g_{k}\inv H_k^{g_{k}\inv}$. So, from now on ``coset labelling'' means ``left coset labelling''. It is clear that one version is as good as the other since every left coset is a right coset (with respect to some conjugate subgroup) and vice versa.
 The reason why left cosets occur in~\cite{HL} where right cosets occur in our context is that Herwig--Lascar let (partial) automorphisms act on the left while we let them act on the right. We come to the first variant of the group-theoretic version of Herwig--Lascar, which is a consequence of the Proposition~\ref{prop:HL>Ash}. 
 
 \begin{Prop}\label{prop:HLgroupformulationABO} Let $\mE=(V\cup K;\alpha,\omega,\inv)$ be a finite graph, $\ell\colon K\to 2^F$, $k\mapsto g_kH_k$ a coset labelling and $(X_k)_{k\in K}$ variables. Then there exists a finite index normal subgroup $N$ of $F$ such that: if the system of equations
 \begin{equation}\label{eq:c-path-equation}
X_{k_1}\cdots X_{k_t}= 1_F\colon \ k_1\cdots k_t\mbox{ a closed path in }\mE
\end{equation}
 has no solution in $F^K$ subject to the constraints $X_k\in g_kH_k$, then it even has no solution subject to the constraints $X_k\in g_kH_kN$.
 \end{Prop}

 \begin{Lemma}\label{lem:groupversionHL1} Proposition~\ref{prop:HL>Ash} and Proposition~\ref{prop:HLgroupformulationABO} are equivalent.     
 \end{Lemma}
 \begin{proof}
     Let $\mE$ be a finite graph with edge set $K$ and let $k\mapsto g_kH_k$ be a (left) coset labelling. For all $k\in K$, let $x_k\in g_kH_k$ and $n_k\in N$. If for some $k_1,\dots, k_t\in K$, $x_{k_1}n_{k_1}\cdots x_{k_t}n_{k_t}=1_F$ then $x_{k_1}\cdots x_{k_t}\in N$. This guarantees that Proposition~\ref{prop:HLgroupformulationABO} implies Proposition~\ref{prop:HL>Ash}. For the converse assume 
   that $(x_k)_{k \in K}$ is a $G$-commuting word labelling, i.e.\  
     that the choice of elements $x_k\in g_kH_k$ is such that for all closed paths $k_1\cdots k_t$ in $\mE$, the product $x_{k_1}\cdots x_{k_t}$ belongs to $N$. We need a choice of elements $n_k\in N$ for $k\in K$ such that for every such product $x_{k_1}\cdots x_{k_t}$,
     \begin{equation}\label{eq:cycleproductN=1}
      x_{k_1}n_{k_1}\cdots x_{k_t}n_{k_t}=1_F.   
     \end{equation} For every single such product a suitable choice is obvious. But we need a choice that      
     works globally. For every connected component $\mD$ of $\mE$ choose a spanning tree $\mT_\mD$ (that is, a 
     tree subgraph that connects
     all vertices of $\mD$). For each edge $k\in \mT_\mD$ set $n_k=1_F$. For $k\in \mD\setminus \mT_\mD$ let $k_1\cdots k_s$ be the unique reduced path $\omega k\longrightarrow \alpha k$ in $\mT_\mD$. According to our assumption, $x_{k_1}\cdots x_{k_s}x_k=:m_k\in N$. Let $n_k:=m_k\inv$; then, since $n_{k_i}=1_F$ for the edges $k_i$,
     \[x_{k_1}n_{k_1}\cdots k_{k_s}n_{k_s}x_kn_k=x_{k_1}\cdots x_{k_s}x_km_k\inv =1_F\]
     for all closed paths of the mentioned kind. By a standard argument (see e.g.~the proof of Lemma 6.1 in~\cite{KapMas}) it follows that for this choice of $(n_k)_{k\in K}$ the equality~\eqref{eq:cycleproductN=1} holds for every closed path $k_1\cdots k_t$ in $\mE$. 
 \end{proof}
By the standard argument from~\cite{KapMas} mentioned in the preceding proof, the infinite system \eqref{eq:c-path-equation} is equivalent to some finite subsystem consisting of equations where the closed path $k_1\cdots k_t$ runs through a set of reduced ``basic cycles'' 
whose cardinality is 
the cyclomatic number 
of $\mE$,
together with all ``trivial cycles''  $kk\inv$ for $k\in K$. 

Almeida and Delgado discussed the connection between Theorem~3.3 in~\cite{HL} (i.e.~the group-theoretic version of the Herwig--Lascar Theorem) and graph labellings in~\cite{almeidadelgado1,almeidadelgado2} and showed that this theorem  is equivalent to Ash's Theorem for inverse monoids (Theorem~\ref{thm: ash} in the present paper). 
We need some further concepts for our discussion. Let $\mE=(V\cup K;\alpha,\omega,\inv)$ be a finite connected graph. An \emph{extended word labelling} of $\mE$ is a mapping $w\colon V\cup K\to F$ such that $w\upharpoonright K$ is a word labelling. We also define \emph{extended coset labelling} in the obvious way: it is a mapping $\ell\colon V\cup K\to 2^F$ such that $\ell(p)=g_pH_p$ for every $p\in V\cup K$ with all $H_p$ finitely generated subgroups of $F$ and $\ell\upharpoonright K$ is a coset labelling. 

Let $G$ be an $A$-generated group. An extended word labelling $w$ of $\mE$  is \emph{$G$-coherent} if for each edge $k\in K$, the equality \[[w(\alpha k)]_G[w(k)]_G=[w(\omega k)]_G\] holds. This simply means that $ab=c$ for all edges $k\colon \underset{a}{\bullet}\overset{b}{\longrightarrow}\underset{c}{\bullet}$ where $a=[w(\alpha k)]_G$, $b=[w(k)]_G$ and $c=[w(\omega k)]_K$. 
This notion is strongly connected with $G$-commutativity, as Lemmas~\ref{lem:comp-coh-L1} and~\ref{lem:comp-coh-L2} below show. 

In the following, $\mE=(V\cup K;\alpha,\omega,\inv)$ is a finite connected graph. 
\begin{Lemma}\label{lem:comp-coh-L1}
	(1)	Let $w\colon K\to F$ be a $G$-commuting word labelling of $\mE$; then for every $v\in V$ and $g\in G$ there exists a
$G$-coherent extended word labelling $\wh{w}$ of $\mE$
such that $[\wh{w}(v)]_G=g$ and $\wh{w}\upharpoonright K=w$.

(2) Conversely,  for every $G$-coherent extended word labelling $w$ of $\mE$ the restriction $w\upharpoonright K$ commutes over $G$.
	\end{Lemma}
\begin{proof} Let $w$ be $G$-commuting; let $f\in F$ be a word such that $[f]_G=g$ and for every vertex $u\in V$ let $\pi_u\colon v\longrightarrow u$ be a path from $v$ to $u$. For $u\in V$ define $\wh{w}(u):=fw(\pi_u)$ and for $k\in K$ let $\wh{w}(k);=w(k)$. This gives an extended word labelling as required by (1). For (2) let $w$ be a $G$-coherent extended word labelling. 
It follows by induction on the length that for every path $\pi\colon u\longrightarrow v$, $[w(u)]_G[w(\pi)]_G=[w(v)]_G$, hence $[w(\pi)]_G=1_G$ for every closed path $\pi$.
\end{proof}

We use an extension of the given graph $\mE$, 
which in~\cite{almeidahyperdecidability,almeidadelgado1,almeidadelgado2} 
was called the \emph{cone} of $\mE$.
Roughly speaking, we add to the graph $\mE$ a new vertex $v_0$, and for every old vertex $v$ a new edge $\overrightarrow{v}\colon  v_0\longrightarrow v$ and its inverse $\overleftarrow{v}$. More formally, given $\mE=(V\cup K;\alpha,\omega,\inv)$ we let
\[\wh{V}:=V\cup \{v_0\} \mbox{ where }v_0\notin V\] and 
\[\wh{K}:=K\cup \{\overrightarrow{v},\overleftarrow{v}\colon v\in V\} \mbox{ where }\overrightarrow{v}\ne\overleftarrow{v}\mbox{ and }\overrightarrow{v},\overleftarrow{v}\notin K\]
and set $\wh{\mE}=(\wh{V}\cup\wh{K};\alpha,\omega,\inv)$ with incidence functions and inversion on the new edges $\overrightarrow{v}$ and $\overleftarrow{v}$ defined in the obvious way: $\alpha \overrightarrow{v}=v_0$, $\omega \overrightarrow{v}=v$ and $\overrightarrow{v}\inv=\overleftarrow{v}$ (this defines also $\alpha \overleftarrow{v}$ and $\omega\overleftarrow{v}$). The graph $\wh{\mE}$ is designed to carry within a word labelling the information contained in an extended word labelling of $\mE$. Indeed, given an extended word labelling $w\colon V\cup K\to F$ we define the word labelling $\wh{w}$ of $\wh{\mE}$ by $\wh{w}(k)=w(k)$ for every $k\in K$ and $\wh{w}(\overrightarrow{v})=w(v)$ (and hence also $\wh{w}(\overleftarrow{v})=w(v)\inv$). A proof of a statement equivalent to the following one forms a part of the proof of Theorem 7 in~\cite{almeidahyperdecidability}; for the convenience of the reader we include a proof.

\begin{Lemma}\label{lem:comp-coh-L2}
An extended word labelling $w$ of $\mE$ is $G$-coherent if and only if the word labelling $\wh{w}$ of $\wh{\mE}$ commutes over $G$.
\end{Lemma}
\begin{proof}
For the `if' part it suffices to notice that for any edge $k$, the equality \[[w(\alpha k)]_G[w(k)]_G=[w(\omega k)]_G\] holds if and only if $[\wh{w}(\overrightarrow{\alpha k}\cdot k\cdot\overleftarrow{\omega k})]_G=1_G$. For the `only if' part let $e_1e_2\cdots e_n$ be a closed path in $\wh{\mE}$. Immediately from the definition of the labelling $\wh{w}$ it follows that the claim is true for every trivial cycle path $kk\inv$ with $k\in \wh{K}$. 
By induction on the length of the path one can see that, in order to show that the $G$-value of its label is $1_G$ it is sufficient to assume that the path is simple, that is,  with the exception of the initial and terminal vertices, no vertex is met twice. 
Also, instead of the path $e_1\cdots e_n$ we may look at any shifted path $e_i\cdots e_ne_1\cdots e_{i-1}$. Hence we may assume that the vertex $v_0$ is either the initial and the terminal vertex, or else does not occur at all. In the first case, the path is of the form $\overrightarrow{u_1}e_2\cdots e_{n-1}\overleftarrow{u_{n-1}}$ for the vertices  $u_1=\alpha e_2$ and $u_{n-1}=\omega e_{n-1}$. Let $u_2,\dots, u_{n-2}$ be the other vertices of the path $e_2\dots e_{n-1}$, that is $\omega e_i=u_i=\alpha e_{i+1}$. Then $G$-coherence of the labelling $w$ implies that $[\wh{w}(e_i)]_G=[\wh{w}(\overleftarrow{u_{i-1}}\overrightarrow{u_i})]_G$ for every $i$. 
It follows that \[[\wh{w}(\overrightarrow{u_1}e_2\cdots e_{n-1}\overleftarrow{u_{n-1}})]_G=[\wh{w}(\overrightarrow{u_1}\overleftarrow{u_1} \overrightarrow{u_2}\cdots \overleftarrow{u_{n-2}} \overrightarrow{u_{n-1}}\overleftarrow{u_{n-1}})]_G=1_G.\]
If $v_0$ does not occur on the path, let $u_1,\dots, u_n$ be the vertices of the path, that is, $\omega e_n=u_1=\alpha e_1$ and $\omega e_i=u_{i+1}=\alpha e_{i+1}$ for $i=1,\dots, i-1$. Again by use of $[\wh{w}(e_i)]_G=[\wh{w}(\overleftarrow{u_{i}}\overrightarrow{u_{i+1}})]_G$ we get 
\[[\wh{w}(e_1\cdots e_n)]_G=[\wh{w}(\overleftarrow{u_1}\overrightarrow{u_2}\overleftarrow{u_2}\cdots \overrightarrow{u_{n}}\overleftarrow{u_{n}}\overrightarrow{u_1})]_G=[\wh{w}(\overleftarrow{u_1}\overrightarrow{u_1})]_G=1_G.\]
\end{proof}

We now define an extended coset labelling   $\ell\colon V\cup K\to 2^F$, $p\mapsto g_pH_p$ of a graph $\mE$ to be $G$-coherent for $G=F/N$ if there exists an extended word labelling $l\colon V\cup K\to F$ such that $l(p)\in \ell(p)$ for every $p\in V\cup K$ and $l$ is $G$-coherent. This may be formulated as follows: for every $p\in V\cup K$ there exists a choice $x_p\in g_pH_pN$ such that the equality $x_{\alpha k}x_k=x_{\omega k}$ holds
for all $k\in K$. 
This leads us to the Almeida--Delgado version of the (group-theoretic formulation of the) Herwig--Lascar Theorem~\cite{almeidadelgado1,almeidadelgado2}.

\begin{Prop}\label{prop:ashalmeida} Let $\mE=(V\cup K;\alpha,\omega,\inv)$ be a finite graph and let $\ell\colon V\cup K \to 2^F$, $p\mapsto g_pH_p$ be an extended  coset labelling  and let $(X_p)_{p\in V\cup K}$ be variables. Then there exists a finite index normal subgroup $N$ of $F$ such that: if the system
\begin{equation}\label{eq:graph system}
	X_{\alpha k}X_k=X_{\omega k}\colon  k\in K
\end{equation}	
	 has no solution in $F^{K\cup V}$ subject to the constraints $X_p\in g_pH_p$  for all $p\in V\cup K$, then it even has no solution subject to the constraints $X_p\in g_pH_pN$ for all $p\in V\cup K$.  
\end{Prop}

 The equations in~(\ref{eq:graph system}) do not require inversion and have the advantage that they make sense in situations other than groups. This led to the concepts of \emph{hyperdecidabilty} and \emph{tameness} of pseudovarieties of finite semigroups and  monoids~\cite{almeidatamenes,ACZreducibility}. In fact, the formulation of the statement in~\cite{almeidadelgado1,almeidadelgado2} does not require graphs with inversion of edges. The connections between $G$-commutativity and $G$-coherence discussed in Lemmas~\ref{lem:comp-coh-L1} and~\ref{lem:comp-coh-L2} show that Propositions~\ref{prop:HLgroupformulationABO} and~\ref{prop:ashalmeida} are equivalent. More precisely, the implication Proposition~\ref{prop:HLgroupformulationABO} $\Longrightarrow$ Proposition~\ref{prop:ashalmeida} follows from Lemma~\ref{lem:comp-coh-L2} while the reverse implication Proposition~\ref{prop:ashalmeida} $\Longrightarrow$ Proposition~\ref{prop:HLgroupformulationABO} follows from Lemma~\ref{lem:comp-coh-L1}.
 \begin{Lemma}
 	Proposition~\ref{prop:HLgroupformulationABO} and Proposition~\ref{prop:ashalmeida} are equivalent.
 \end{Lemma}

 We now come to the formulation of the group-theoretic version of the Theorem of Herwig--Lascar~\cite[Theorem 3.3]{HL}.
For finite index sets $I$ and $L$ consider tuples $(X_i)_{i\in I}$ of variables for  elements from $F$, $(g_l)_{l\in L}$ of  elements of $F$, and $(H_l)_{l\in L}$ of finitely generated subgroups of $F$. In their paper~\cite{HL} Herwig and Lascar considered \textsl{equations} of the form $X_i\equiv X_jg_l \bmod H_l$ and $X_i\equiv g_l \bmod H_l$, respectively. Elements $x_i,x_j\in F$ \emph{solve} these equations if and only if \[x_iH_l=x_jg_lH_l\mbox{ and }x_iH_l=g_lH_l,\] which is of course equivalent to $ x_j\inv x_i\in  g_lH_l$ and  $x_i\in g_lH_l$, respectively. 
Hence these equations are equivalent to \textsl{constraints} of the form $X_j\inv X_i\in g_lH_l$ and $ X_i\in g_lH_l$ the second of which occurs already in Propositions~\ref{prop:HLgroupformulationABO} and~\ref{prop:ashalmeida}.
Theorem 3.3 in~\cite{HL} (i.e.~the group-theoretic version of the Herwig--Lascar Theorem) then can be formulated as follows.

\begin{Prop}\label{prop:HLgroupformulation} 
	Let $I$ and $L$ be finite index sets, let $(X_i)_{i\in I}$, $(g_l)_{l\in L}$ and $(H_l)_{l\in L}$ be variables, elements of $F$ and finitely generated subgroups of $F$ and let $S\subseteq I\times L$ and $T\subseteq I^2\times L$. Then there exists a finite index normal subgroup $N$ of $F$ such that: if the system of constraints  
\begin{equation} \label{eq:meprop}   X_i\in g_lH_l\colon (i,l)\in S \; \mbox{ and } \; X_j\inv X_i\in g_lH_l\colon (j,i,l)\in T  \end{equation}
    has no solution $(x_i)\in F^I$, then neither has the system
    \[X_i\in g_lH_lN\colon (i,l)\in S \; \mbox{ and } \; X_j\inv X_i\in g_lH_lN\colon (j,i,l)\in T.\] 
\end{Prop}
The Theorem of Ribes--Zalesskii can be expressed this way as follows: let $H_1,\dots, H_n$ be finitely generated subgroups of $F$ and $g\in F$. Consider variables $X_1,\dots, X_{n-1}$ and constraints
\[X_{n-1}\in gH_n,\ X_{n-1}\inv X_{n-2}\in H_{n-1},\dots, X_2\inv X_1\in H_2, X_1\in H_1.\]
One readily checks that this system has a solution if and only if  $1_F\in gH_n\cdots H_1$, which is equivalent to $g\in H_1\cdots H_n$. Hence the theorem implies that, if $g\notin H_1\cdots H_n$, then there exists a finite index normal subgroup $N$ of $F$ such that $g\notin (H_1N)\cdots(H_nN)=H_1\cdots H_nN$, which implies $gN\cap H_1\cdots H_n=\varnothing$. So, Proposition~\ref{prop:HLgroupformulation} can indeed be seen as a generalisation of the Ribes--Zalesskii Theorem.

We shall argue that Propositions~\ref{prop:ashalmeida} and~\ref{prop:HLgroupformulation} are equivalent. The easier direction is that Proposition~\ref{prop:HLgroupformulation} implies Proposition~\ref{prop:ashalmeida}. Suppose we are given the situation described in the latter proposition. Solvability of the equation 
\[X_{\alpha k}X_k=X_{\omega k}\] subject to the constraints
\[X_{\alpha k}\in g_{\alpha k}H_{\alpha k},X_{k}\in g_{ k}H_{k},X_{\omega k}\in g_{\omega k}H_{\omega k}\]
just means that the statement
\begin{equation}\label{eq:eqalmdel}
\exists X_{\alpha k}
\exists X_k
\exists X_{\omega k}
\left[ 
\begin{array}{r@{\;\;}l}
& X_{\alpha k}\in g_{\alpha k}H_{\alpha k} \;\wedge\; X_{k}\in g_{ k}H_{k}
\\
\wedge & X_{\omega k}\in g_{\omega k}H_{\omega k} \;\wedge\; X_{\alpha k}X_k=X_{\omega k}
\end{array}
\right]
\end{equation}
is true. But \eqref{eq:eqalmdel} is equivalent to
\begin{equation}\label{eq:HL}
\exists X_{\alpha k}
\exists X_{\omega k}
\Bigl[ 
X_{\alpha k}\in g_{\alpha k}H_{\alpha k} \;\wedge\; 
X_{\omega k}\in g_{\omega k}H_{\omega k} \;\wedge\; 
X_{\alpha k}^{-1}X_{\omega k} \in g_{k}H_{k}
\Bigr]
\end{equation}
and the equivalence of these statements remains valid if every coset $g_lH_l$ involved is replaced by the coset $g_lH_lN$, for any fixed normal subgroup $N$ of $F$. We see that the conditions occurring in \eqref{eq:HL} are just constraints 
as in Proposition~\ref{prop:HLgroupformulation}. Hence, for every $k\in K$, 
\eqref{eq:eqalmdel} may be replaced 
by the equivalent~\eqref{eq:HL}. We therefore may convert the system of Proposition~\ref{prop:ashalmeida} into an equivalent one of the form in Proposition~\ref{prop:HLgroupformulation}, and the equivalence remains valid if each coset $g_lH_l$ is replaced by the coset $g_lH_lN$, for any fixed (finite index) normal subgroup of $F$. 
Thus Proposition~\ref{prop:HLgroupformulation} implies Proposition~\ref{prop:ashalmeida}.

For the converse we have to be more careful. Consider a finite system as in Proposition~\ref{prop:HLgroupformulation}. We first look at the constraints $X_i\in g_lH_l$ in \eqref{eq:meprop}. If for some variable $X_i$ no such constraint occurs, then we add the trivial one $X_i\in F$. Hence we may assume that \textsl{every} variable $X_i$ occurs in at least one such constraint. So, for each $i\in I$ there exist 
group elements $g_{i1},\dots, g_{in_i}\in F$ 
for some $n_i\ge 1$, 
and finitely generated subgroups $H_{i1},\dots, H_{in_i}$ of $F$ such that the constraint $X_i\in g_lH_l$  occurs in \eqref{eq:meprop} if and only if $g_lH_l=g_{is}H_{is}$ for some $s\in \{1,\dots,n_i\}$. For every $i\in I$ we introduce $n_i$ new variables $X_{i1},\dots,X_{in_i}$ and consider a new system in the variables $V:=\bigcup_{i\in I}\{X_{i1},\dots,X_{in_i}\}$ defined as follows: instead of the constraints
\[X_i\in g_{i1}H_{i1},\dots,X_i\in g_{in_i}H_{in_i}\]
take the constraints
\begin{equation}\label{eq:new equations1}\begin{split}
&X_{i1}\in g_{i1}H_{i1},\dots, X_{in_i}\in g_{in_i}H_{in_i}\mbox{ and }X_{is}\inv X_{it}\in\{1_F\}\\ &\mbox{ for all }s,t=1,\dots, n_i\mbox{ and } i\in I,
\end{split} 
\end{equation}
and, for all $i,j\in I$,  instead of every single constraint $X_j\inv X_i\in g_lH_l$ consider all constraints
\begin{equation}\label{eq:new equations2}
X_{js}\inv X_{it}\in g_lH_l\mbox{ for all }s=1,\dots,n_j\mbox{ and }  t=1\dots, n_i.
\end{equation}
The new system \eqref{eq:new equations1}\&\eqref{eq:new equations2}  
  is again a system of the same type as the one in \eqref{eq:meprop}. What is more, the old system \eqref{eq:meprop} admits a solution if and only if the new system \eqref{eq:new equations1}\&\eqref{eq:new equations2}  admits one. Again this equivalence concerning solvability remains valid if 
each coset $g_lH_l$ involved is replaced by $g_lH_lN$, for any fixed normal subgroup $N$.

An important feature of the new system is that for \textsl{every} variable  $X$ there is exactly one constraint of the form $X\in gH$. We now proceed as in~\cite[p.~415]{almeidadelgado1}\footnote{The argument in~\cite{almeidadelgado1} does not take into account that some variables $X_i$ in \eqref{eq:meprop} may admit more than one constraint $X_i\in g_lH_l$.} and construct a ``graph system'' as in Proposition~\ref{prop:ashalmeida} that is equivalent to \eqref{eq:new equations1}\&\eqref{eq:new equations2}. We define a graph and an extended coset labelling $\ell$ as follows: the set of vertices is the set of variables of the system \eqref{eq:new equations1}\&\eqref{eq:new equations2}, and we set $\ell(X_{it})=g_{it}H_{it}$ for every such vertex (i.e.~variable); further, for two vertices (i.e.~variables) $X_{it}, X_{js}\in V$ and every constraint $X_{js}\inv X_{it}\in g_lH_l$ we impose an edge $k\colon X_{js}\longrightarrow X_{it}$ with label $\ell(k)=g_lH_l$; if necessary (that is, if the corresponding ``inverse'' constraint $X\inv_{it}X_{js}\in {g_l\inv}H_l^{g_l\inv}$ is not already contained in \eqref{eq:new equations2}) we add for every such edge its inverse edge $k\inv\colon X_{it}\longrightarrow X_{js}$ and label it by $\ell(k\inv)=g_l\inv H_l^{g_l\inv}$. 
The graph and labelling so obtained provide a system of equations and a set of constraints as in Proposition~\ref{prop:ashalmeida}. By the same arguments as discussed in connection with the equivalence of \eqref{eq:eqalmdel} and \eqref{eq:HL}, we see that \eqref{eq:new equations1}\&\eqref{eq:new equations2} and the graph system are equivalent, and this is still true if each 
constraint $g_lH_l$ is replaced by $g_lH_lN$ for some fixed normal subgroup $N$ of $F$. Altogether, this means that, if we start with a system of the form \eqref{eq:meprop} that has no solution, and if we assume Proposition~\ref{prop:ashalmeida} to be true, then we end up with the existence of a finite index normal subgroup $N$ of $F$ such that the system obtained from  \eqref{eq:meprop} by replacement of each coset $g_lH_l$ by the coset $g_lH_lN$ still has no solution. This just means that Proposition~\ref{prop:ashalmeida} implies Proposition~\ref{prop:HLgroupformulation}. We have therefore proved
the following.
\begin{Lemma}
	Proposition~\ref{prop:ashalmeida} and Proposition~\ref{prop:HLgroupformulation} are equivalent.
\end{Lemma}

\subsection{Construction of groups {required for} Ash's Theorem}\label{sec:groups for Ash} We briefly discuss the existence of and how to possibly construct, for a given finite graph $\mE$ and a given finite $A$-generated inverse monoid $M$,  a finite $A$-generated
group $G$ as  in Theorem~\ref{thm: ash}.
Let us call such a group $G$ \emph{suitable} for $M$ and $\mE$. 
Starting point is a finite $A$-generated group $G$  that extends $M$ in
the sense of Definition~\ref{def: G extends M}.
Such a group then is suitable for the given inverse
monoid $M$ and the cycle graph $\mC_1$ of length $1$ (that is,
for the graph having one vertex and one geometric edge). Indeed, let $w\colon \{e,e\inv\}\to \til A^*$  be a $G$-commuting labelling; that is, $[w(e)]_G=1_G$. Since $G$ extends $M$, $[w(e)]_M$ is an idempotent and therefore $[w(e)]_M=[w(e)\inv]_M=[w(e)w(e)\inv]_M$; 
then $v(e):=w(e)w(e)\inv$ provides an $M$-related $F$-commuting labelling.

Throughout this subsection denote by $\mC_i$ the cycle graph of length $i$.
We start with the construction of Ash. Based on $G_1:=G$ one can construct an expansion  $G_2$ of $G_1$ that is suitable for $M$ and $\mC_2$, an expansion $G_3$ of $G_2$ which is suitable for $M$ and $\mC_3$, and so on, i.e.\ inductively an expansion $G_k$ of $G_{k-1}$ that is suitable for $M$ and 
$\mC_k$, for all $k$. 
Every group $G_k$ is constructed in order to avoid relations of a certain kind which are satisfied by its predecessor $G_{k-1}$.  The relations to be avoided by  $G_k$ can be described in terms of the Cayley graph $\mG_{k-1}$ of the preceding group $G_{k-1}$. In order to make this more precise it is convenient to introduce the following concepts.

In the following  we shall consider the indices in an $n$-tuple 
$(x_1,\dots,x_n)$ as coming form a cyclic index set, that is, the indices in expressions like $x_{i\pm 1}$  have to be understood $\bmod\ n$; this situation will be expressed by writing the tuple as $(x_i)_{i\in \mathbb{Z}_n}$. Let $G$ be an $A$-generated group with Cayley graph $\mG$. A \emph{$2$-cyclic configuration} $(\mP_i,g_i)_{i\in \mathbb{Z}_2}$ of $\mG$ consists of a pair $\mP_1,\mP_2$ of connected subgraphs of $\mG$ and a pair $g_1,g_2$ of elements of $G$ such that $g_2=1_G$ and $g_1$ and $g_2$ are contained in distinct connected components of $\mP_1\cap \mP_2$. For $n\ge 3$, an \emph{$n$-cyclic configuration} $(\mP_i,g_i)_{i\in\mathbb{Z}_n}$ of $\mG$ consists of an $n$-tuple $\mP_1,\dots,\mP_n$ of connected subgraphs of $\mG$ and an $n$-tuple $g_1,\dots,g_n$ of elements of $G$ such that $g_n=1_G$ and
\begin{itemize}
	\item  $g_i\in \mP_i\cap \mP_{i+1}$ for all $i$ and
	\item $\mP_{i-1}\cap \mP_i\cap\mP_{i+1}=\varnothing$ for all $i$.
\end{itemize}
The definition of $n$-cyclic configurations is reminiscent of the one of $n$-coset cycles (Definition~\ref{def: cosetcycle}) in Section~\ref{sec:strengthened}; however, the context is considerably different.
 Figure~\ref{fig:3cycle} shows examples of $2$-cyclic and $5$-cyclic configurations.
\begin{figure}[ht]
\begin{tikzpicture}[scale=1.3]
	\filldraw (0,-0) circle (1.5pt); 
	\draw[below](0,-0.05) node{$g_2$};
	\draw plot [smooth cycle] coordinates {(-0.5,2.8)(0.5,2.8)(-0.7,1.25)(0.5,-0.3)(-0.5,-0.3)(-1.7,1.25)};
	\draw plot [smooth cycle] coordinates {(0.5,2.8)(-0.5,2.8)(0.5,1.25)(-0.5,-0.3)(0.5,-0.3)(1.5,1.25)};
	\draw(-1.2,1.3)node{$\mathcal{P}_1$};
	\draw(1,1.3)node{$\mathcal{P}_2$};
	\draw[above](0,2.5)node{$g_1$};
	\filldraw(0,2.5)circle(1.5pt);
\end{tikzpicture}\qquad \quad
\begin{tikzpicture}[scale=1.3]
	\draw plot [smooth cycle] coordinates {(-1,2) (0,1.9) (1.4,2.2)( 1.5,3) (-1,3)};
	\draw plot [smooth cycle] coordinates {(-1.7,2.8)(-0.4,3)(-1,0.7)(-1.8,0.7)};
	\draw plot [smooth cycle] coordinates  {(1.7,3)(0.4,3)(1,0.7)(1.9,0.7)};
	\draw plot [smooth cycle] coordinates  {(1,1.2)(2,1.3)(0.7,-0.7)(-0.5,-0.5)};
	\draw plot [smooth cycle] coordinates   {(-1,1.2)(-2,1.3)(-0.7,-0.7)(0.5,-0.5)};
	\filldraw (0,-0.3) circle (1.5pt); 
	\draw[below](0,-0.3) node{$g_5$};
	\filldraw(1.5,1) circle (1.5pt);
	\draw[right](1.4,1.2)node{$g_4$};
	\filldraw(-1.5,1)circle(1.5pt);
	\draw[left](-1.43,1)node{$g_1$};
	\filldraw(-0.8,2.5)circle(1.5pt);
	\draw[above](-0.8,2.5)node{$g_2$};
	\filldraw(0.8,2.5)circle(1.5pt);
	\draw[above](0.8,2.5)node{$g_3$};
	\draw(-0.7,0.3) node{$\mathcal{P}_1$};
	\draw(0.7,0.3) node{$\mathcal{P}_5$};
	\draw(-1.2,1.7)node{$\mathcal{P}_2$};
	\draw(1.2,1.7)node{$\mathcal{P}_4$};
	\draw(0,2.5)node{$\mathcal{P}_3$};
	\end{tikzpicture}	
\caption{$2$-cyclic and $5$-cyclic configurations}\label{fig:2&5cycles}
\end{figure}

An $A$-generated group $H$ \emph{avoids} the $n$-cyclic configuration
 $(\mP_i,g_i)_{i \in {\mathbb{Z}_n}}$ of $\mG$ if the following holds:
 \begin{quotation}
	  \noindent
	 for every
	 $n$-tuple of words $u_1,\dots, u_n$ in $\til A^*$, if
	 $u_i$ labels a path from $g_{i-1}$ to $g_i$ in $\mG$ that runs
	 entirely in $\mP_i$, then
	 $[u_1\cdots u_n]_H\ne 1_H$.
	 \end{quotation}
 In other words, generator sequences that follow the given 
 cyclic configuration of the graphs $\mP_i$ in $\mG$ do not give rise to cycles in the Cayley graph $\mH$ of $H$.

 A finite group $H$ which avoids a given $n$-cyclic configuration can be constructed as follows. We form the disjoint union $\mP_1\sqcup\cdots \sqcup \mP_n$ of its graphs and factor by the smallest $A$-graph congruence that identifies the vertices $g_i^{(i)}$ of $\mP_i$ and $g_i^{(i+1)}$ of $\mP_{i+1}$ for $i=1,\dots, n-1$, but not the vertices $g_n^{(n)}$ of $\mP_n$ and $g_n^{(1)}$ of $\mP_1$. The result is a finite connected $A$-graph $\mP$. The transition group of any finite completion $\ol{\mP}$ of $\mP$ then avoids the given configuration. Since every finite $\mG$ has only finitely many $n$-cyclic configurations we can do this for each one; the $A$-generated direct product of the resulting groups then is a finite $A$-generated group that avoids all $n$-cyclic configurations of $\mG$. The $A$-generated direct product of this group with the original group $G$ then is an expansion $H$ of $G$ that avoids all $n$-cyclic configurations of $\mG$.

We get a  (non-deterministic) procedure that, given as input a finite $A$-generated inverse monoid $M$ and a
positive integer $n$, outputs a finite $A$-generated group $G$
suitable for $M$ and $\mC_n$ (recall that this means that $G$ satisfies the condition of Theorem~\ref{thm: ash} with respect to $M$ and $\mC_n$). This procedure is as follows: starting with a group $G_1$ which is suitable for $M$ and $\mC_1$ (that is, $G_1$ extends $M$ in the sense of Definition~\ref{def: G extends M}) we obtain a series of expansions $G_1\twoheadleftarrow G_2\twoheadleftarrow\cdots\twoheadleftarrow G_{n-1}\twoheadleftarrow G_n\twoheadleftarrow\cdots$ such that, in this series every $G_k$ avoids all $k$-cyclic configurations of the Cayley graph $\mG_{k-1}$ of its predecessor $G_{k-1}$. Then, for every $n$, the $n$th group $G_n$ in this series is suitable for $M$ and $\mC_n$ (and hence for $\mC_m$ for all $m \leq n$). 
Let us call this procedure \emph{$n$-cycle procedure}.

Next we obtain a procedure which on the input of a finite $A$-generated inverse monoid $M$ and two positive integers $k,n$ outputs a group $G$ which is suitable for $M$ and for all connected graphs $\mE$  of cyclomatic number (at most) $k$ on (at most) $n$ vertices.
Indeed, we may iterate 
the $n$-cycle procedure $k$ times. Suppose the inverse monoid $M$ and the numbers $n$ and $k$ are given. We let $G_1$ be a group suitable for $M$ and $\mC_n$, $G_2\twoheadrightarrow G_1$ a group suitable for $M\times G_1$ and $\mC_n$, $G_3\twoheadrightarrow G_2$ a group suitable for $M\times G_1\times G_2$ and $\mC_n$,  and so on, $G_k\twoheadrightarrow G_{k-1}$ a group suitable for $M\times G_1\times \cdots\times  G_{k-1}$ and $\mC_n$. 
Each group $G_i$ is constructed by the application of the $n$-cycle procedure. By part 1 of the \textsl{proof of Lemma~\ref{lem: key lemma RZ->Ash} by use of Theorem~\ref{thm:ash for cycles}} it follows that every group $G_i$ is critical for the pair $(M\times G_1\times\cdots\times G_{i-1},\mC_n)$ in the sense of Definition~\ref{def: critical group}. Recall that $\times$ denotes the $A$-generated direct product. 
Then iterated application of \eqref{eq:Hcritical for (MtimesG)} 
shows that $G_k$ is critical for the pair $(M,\mE)$ and therefore is suitable for $M$ and $\mE$, provided that $\mE$ is connected on $n$ vertices and has cyclomatic number $k$. 
  In order to see this, we consider a chain of connected graphs $\mE_0\subseteq \mE_1\subseteq \cdots\subseteq \mE_{k-1}\subseteq \mE_k=\mE$ where every $\mE_i$ has cyclomatic number $i$ 
  and $\mE_{i+1}$ is obtained from $\mE_i$ by adding an edge (and its inverse). (We note that if $H\twoheadrightarrow G$ then $G\times H\cong H$ for the $A$-generated direct product, so, in fact, $G_1\times G_2\times \cdots \times G_i\cong G_i$ for all $i$, but in the following it is useful to write the products in their redundant form.) 
  Firstly
  we have that $G_k$ is critical for $(M\times G_1\times\cdots\times G_{k-1},\mC_n)$ 
  and $G_k$ is critical for $((M\times G_1\times\cdots\times G_{k-1})\times G_k,\mE_0)$ (the latter holds for trivial reasons since $\mE_0$ is a tree), hence by \eqref{eq:Hcritical for (MtimesG)}, $G_k$ is critical for $(M\times G_1\times \cdots\times G_{k-1},\mE_1)$. Next, $G_{k-1}$ and therefore also $G_k$ is critical for $(M\times G_1\times \cdots\times G_{k-2},\mC_n)$ (since $G_k\twoheadrightarrow G_{k-1}$) and $G_k$ is critical for $((M\times G_1\times \cdots G_{k-2})\times G_{k-1},\mE_1)$, so again by \eqref{eq:Hcritical for (MtimesG)}, $G_k$ is critical for $(M\times G_1\times  \cdots\times G_{k-2},\mE_2)$. We can repeat this argument until we find that $G_k$ is critical  for $(M,\mE_k)$ and hence is suitable $M$ and $\mE_k=\mE$. 

In order to perform the $i$th iteration of the above-mentioned procedure, that is, to find an expansion $G_i\twoheadrightarrow G_{i-1}$ that is critical for \[(M\times G_1\times \cdots\times G_{i-1},\mC_n),\] it should be noticed that the group $G_{i-1}\cong G_1\times\cdots\times G_{i-1}$ itself extends the inverse monoid $M\times G_1\times\cdots\times G_{i-1}$, hence is the starting group for the application of the $n$-cycle procedure, which for the inverse monoid \[M\times G_1\times \cdots\times G_{i-1}\] outputs a group critical  for $(M\times G_1\times \cdots\times G_{i-1},\mC_n)$. Altogether, 
Ash's procedure for the construction of a group suitable for a given finite inverse monoid $M$ and a given finite connected graph $\mE$ having $n$ vertices and cyclomatic number $k$, boils down to the $k$-fold repetition of the $n$-cycle procedure.

From~\cite{geometryofprofinite} it even follows that a group $H$ suitable for $M$ and any (fixed) finite connected graph $\mE$ is contained in any series of groups $G_1\twoheadleftarrow G_2\twoheadleftarrow \cdots$ starting with a group $G_1$ extending $M$ and such that every group $G_k$ in the series avoids all $2$-cyclic configurations of the Cayley graph of its predecessor $G_{k-1}$. However, the proof is by a non-constructive topological compactness argument. In particular, it does not deliver a bound for a number $t$ for which $H$ can be chosen to be $G_t$. A concrete instance of such a series with an explicit bound is presented below.

The first author, partly in cooperation with B.~Steinberg~\cite{mytype2,AS,ASconstructiveRZ},  found a simpler construction for expansions of a group $G$ which avoid cyclic configurations of $\mG$ and therefore for the construction of groups which are suitable for a given finite inverse monoid $M$ and a given finite connected graph $\mE$. For an $A$-generated group $G$ and an integer $p\ge 2$ there exists an $A$-generated expansion ${}^{\mathbf{Ab}_p}G$ of $G$ that is universal among all $A$-generated expansions $K$ of $G$ for which the kernel of $K\twoheadrightarrow G$ is an Abelian group of exponent dividing $p$: for every such expansion $K$, the canonical morphism ${}^{\mathbf{Ab}_p}G\twoheadrightarrow G$ factors through $K$ and the kernel of ${}^{\mathbf{Ab}_p}G\twoheadrightarrow G$ itself is Abelian group of exponent dividing $p$. If $G\cong F/N$ then ${}^{\mathbf{Ab}_p}G\cong F/[N,N]N^p$ and the group ${}^{\mathbf{Ab}_p}G$ admits a transparent representation as a subgroup of a semidirect product of the free $|G\times A|$-generated Abelian group of exponent $p$  by $G$. It is shown in~\cite{AS}  that if $G$ extends the monoid $M$, then ${}^{\mathbf{Ab}_p}G$ is suitable for $M$ and all graphs with two vertices (in semigroup-theoretic terms this says that ${}^{\mathbf{Ab}_p}G$ is a witness for the group-pointlike sets of $M$). In particular, ${}^{\mathbf{Ab}_p}G$ is suitable for $M$ and the cycle graph $\mC_2$. We can iterate the expansion $G\mapsto {}^{\mathbf{Ab}_p}G\mapsto {}^{\mathbf{Ab}_p}({}^{\mathbf{Ab}_p}G)\mapsto \cdots$ and let $G^{(p,k)}$ denote the outcome of expanding $G$ $k$ times (with respect to $p$). It is shown in~\cite{mytype2} that $G^{(p,n-1)}$ is suitable for $M$ and $\mC_{n}$ for every $n$. We can further use the fact that loop edges (though they contribute to the cyclomatic number)  do not contribute to the complexity concerned in the construction of suitable groups (see the proof of Theorem~\ref{thm: ash} by use of Lemma~\ref{lem: key lemma RZ->Ash}), and that only one expansion of the form $G\mapsto {}^{\mathbf{Ab}_p}G$ is necessary to handle multiple edges. Combination of all this leads to the following 
\begin{Prop}\label{prop:mytype2} Let $M$ be a finite $A$-generated inverse monoid and let $G$ be a finite $A$-generated group extending $M$; let $\mE$ be a connected graph on $n$ vertices and let $\mE^\circ$ be the graph obtained from $\mE$ by removing all loop edges and replacing every multiple edge by a single edge. Let $k$ be the cyclomatic number of $\mE^\circ$. Then $G^{(p,(n-1)k+1)}$ is suitable for $M$ and $\mE$
\end{Prop}
This is Theorem 4.3 in~\cite{mytype2} except that for the number of expansions needed, instead of $(n-1)k+1$ 	the worse bound $(n-1)(2k-1)+1$ had been calculated there. In the sequence of iterations one may actually use different integers $p_i$ at every iteration.

\section{Theorems Strengthened}\label{sec:strengthened}
The groups provided by {the third author's
	construction}~\cite{ABO,otto0,otto3}  deserve special attention (not only) in
semigroup theory. Indeed, they provide
an affirmative solution to the Henckell--Rhodes Problem: every finite inverse monoid admits a
finite $F$-inverse cover --- see~\cite{ABO} for a detailed
presentation of the topic. However, this is not the end of the
story. By taking full advantage of the construction we get a
strengthening of Ash's Theorem; in Section~\ref{sec:strictly stronger} some evidence is given that the new theorem is substantially stronger than Ash's Theorem. 
In turn, this leads to a strengthening of the Herwig--Lascar
Theorem~\ref{thm: Herwig--Lascar}, see Theorem~\ref{thm: HL strengthened} below. 
We also get strengthenings of the group-theoretic versions of Section~\ref{sec: group-theoretic version}. The following is the intended strengthening of
Ash's Theorem.

\begin{Thm}[main theorem]\label{thm:otto-ash} Let $M$ be a finite $A$-generated inverse monoid and $n$ a positive integer. Then there exists a finite $A$-generated group $H$ such that, for every finite connected graph $\mE$ on at most $n$ vertices, every $H$-commuting labelling $v$ of $\mE$ 
admits a labelling $u$ such that
	\begin{enumerate}
		\item $u$ commutes over the free group $F$,
		\item $[v(e)]_M=[u(e)]_M$ for every edge $e\in \mE$,
		\item $[v(e)]_H=[u(e)]_H$ for every edge $e\in \mE$.
	\end{enumerate}
\end{Thm}

In short, every
$H$-commuting labelling of $\mE$ 
admits an $(M\times H)$-related relabelling that commutes over the free group $F$.
This theorem relates to Ash's Theorem just as  the construction
of a finite $A$-generated $F$-inverse cover of a finite $A$-generated inverse monoid $M$ relates to
the construction of a finite group witnessing 
the group-pointlike sets of  $M$
(see the discussion in~\cite[Section 2.7]{ABO}). Note that for a group $H$ satisfying
Theorem~\ref{thm:otto-ash} the inverse monoid $M\times H$ is an
$F$-inverse cover of $M$ {(provided that $\mE$ contains a cycle graph of length at least $2$}). Theorem~\ref{thm:otto-ash} is proved in
Section~\ref{sec:proofs}. 
In preparation for that poof we recall some key concepts and results from~\cite{ABO, otto0, otto3} in Section~\ref{sec:prep}.

\subsection{Summary of some results from~\cite{ABO,otto0,otto3}}\label{sec:prep} We use the notation of~\cite{ABO}. In this section, we shall be concerned with two different generating sets: $A$ and $E$. The set $E$ denotes the set of positive edges of some finite connected oriented graph $\mQ$ (to be specified later). In Section~\ref{sec:generalcase} we treat the general case where $E$ is the set of positive edges of a finite connected oriented graph $\mQ$ while in Section~\ref{sec:Cayleygraphcase}, $E=Q\times A$ is the set of positive edges of the Cayley graph $\mQ$ of some finite $A$-generated group $Q$.

\subsubsection{The general case}\label{sec:generalcase}
For $a\in E$ and $p\in \til{E}^*$ let $p_{a\to 1}$ be the word obtained
from $p$ by deletion of all occurrences of $a$ and $a^{-1}$ in
$p$. 
	An $E$-generated group $G$ is
	\emph{retractable} if, for all words $p,q\in \til{E}^*$ and every
	letter $a\in E$ the following holds:%
	\footnote{It suffices to restrict this postulate to the case $q=1$.}
	\[ [p]_G=[q]_G\Longrightarrow [p_{a\to 1}]_G=[q_{a\to 1}]_G.\]

\begin{Def}\label{def:word content}
For a word $p\in \til{E}^*$, the \emph{word content} $\co(p)$ is the set of all letters $a\in E$ for which $a$ or $a^{-1}$ occurs in $p$.
\end{Def}

Assume that $G$ is retractable: for $p,q\in \til{E}^*$ and $a\in E$ the equality $[p]_G=[q]_G$ implies $[p_{a\to 1}]_G=[q_{a\to 1}]_G$. Suppose now that $a\in \co(p)$ but $a\notin \co(q)$. Then
the words $q$ and $q_{a\to 1}$ are identical. Hence $[p]_G=[q]_G$ implies
\[
[p_{a\to 1}]_G=[q_{a\to 1}]_G=[q]_G=[p]_G.
\]
In this way, we may (without changing the value $[p]_G$) delete every letter in a word $p$  that does not occur in every other word representation $q$ of the group element $[p]_G$. 
This leads to the following definition.
\begin{Def}\label {def:G-content}
	Let $G$ be a retractable $E$-generated group and $g\in G$.
	The \emph{content} $\mathrm{C}(g)$ of $g$ is
	\[\mathrm{C}(g):=\bigcap\big\{\co(q)\colon q\in \til{E}^*, [q]_G=g\big\}.\]
	For a word $p\in \til{E}^*$ the $G$\emph{-content} of $p$ is the content $\mathrm{C}([p]_G)$.
\end{Def}
\begin{Rmk}
 For any retractable $E$-group $G$ and any word $w\in \til E^*$ the word content ${\rm co}(w)$ of $w$ always contains the $G$-content ${\rm C}([w]_G)$ of $w$.  
\end{Rmk}
For a retractable $E$-generated group $G$, the $G$-content of a word $p\in \til E^*$ is the unique $\subseteq$-minimal subset $B$ of $E$ such that $[p]_G$ can be represented as a word over $\til B$. For a given word $p$ denote by $p^\circ$ the word obtained from $p$ by deleting every occurrence of a letter (and its inverse) of $p$ not in the $G$-content of $p$. Then $p^\circ$ is a representation with minimal content of the group element $[p]_G$, that is, $[p^\circ]_G=[p]_G$ and $\mathrm{co}(p^\circ)=\mathrm{C}([p^\circ]_G)=\mathrm{C}([p]_G)$.

Now let $G$ be some retractable $E$-generated group and let $p\in \til E^*$ be a word forming a path $p\colon u\longrightarrow v$ in the graph $\mQ$ (recall that $E$ is the set of positive edges of the oriented graph $\mQ$). If we form the aforementioned minimal representation $p^\circ$ of $[p]_G$ as a word over $\til E$, deletion of letters in $p$ in general destroys the \textsl{path quality} of $p$ with respect to $\mQ$: $p^\circ$ is \textsl{some} minimal representation of $[p]_G$, but we cannot expect that $p^\circ$ forms a path in $\mQ$, let alone a path $u\longrightarrow v$.

The main result of~\cite[Sections 3--5]{ABO} describes the construction of a \textsl{finite} $E$-generated group $G$ having the following crucial property~\cite[Cor.~5.7]{ABO}. 
\begin{Lemma}\label{lem:crucial G}
For every word $p\in \til E^*$ forming a path $u\longrightarrow v$ in $\mQ$, there exists a word $q\in \til E^*$ such that
\begin{enumerate}
	\item $q$ also forms a path $u\longrightarrow v$ in $\mQ$,
	\item $[q]_G=[p]_G$,
	\item $\mathrm{co}(q)=\mathrm{co}(p^\circ)$.
\end{enumerate}
\end{Lemma}

That is, $q$ is $G$-related to $p$, uses only letters (i.e.~edges) from the common $G$-content and also forms a path $u\longrightarrow v$ in $\mQ$. Note that finiteness of $G$ is 
of the essence here, since the free $E$-generated group obviously 
satisfies these conditions.

The group $G$ is constructed in a series of expansions
\begin{equation}\label{eq: expanding groups}
G_1\twoheadleftarrow H_1\twoheadleftarrow G_2\twoheadleftarrow\cdots\twoheadleftarrow G_{|E|-1}\twoheadleftarrow H_{|E|-1}\twoheadleftarrow G_{|E|}=G,
\end{equation}
the members of which are the transition groups of complete $E$-graphs
\begin{equation}\label{eq: ascending graphs}
\mX_1\subseteq \mY_1\subseteq \mX_2\subseteq\cdots\subseteq \mX_{|E|-1}\subseteq \mY_{|E|-1}\subseteq \mX_{|E|},
\end{equation}
that is, $G_i=\cT(\mX_i)$ and $H_i=\cT(\mY_i)$. 
Every graph in the series \eqref{eq: ascending graphs} is obtained from its predecessor by adding certain complete components tailored from the Cayley graph of the corresponding preceding group; the (complicated) details are worked out in~\cite[Sections 3--5]{ABO}. The graph $\mX_1$ is a certain completion of the original graph $\mQ$ that is considered as an $E$-graph in which every edge  
is labelled by itself
(observe that the labelling alphabet $E$ coincides with the set of positive edges $E$ of $\mQ$). 
In addition, the completion $\mX_1$ of $\mQ$ is obtained by adding just edges 
so that $\mQ$ and $\mX_1$ have the same set of vertices. Let us choose and fix an arbitrary vertex of $\mQ$ (i.e.~of $\mX_1$) and denote it by $1_Q$. From the construction it follows that for every group $G_i$ and $H_i$ in the series \eqref{eq: expanding groups} there exists a unique canonical graph morphism $\varphi_i\colon \mG_i\twoheadrightarrow \mX_1$ (respectively $\psi_i\colon \mH_i\twoheadrightarrow \mX_1$) mapping $1_{G_i}$ (respectively $1_{H_i}$) to $1_Q$ where $\mG_i$ (respectively $\mH_i$) is the Cayley graph of $G_i$ (respectively $H_i$). To simplify notation put $H:=H_{|E|-1}$ and recall that $G:=G_{|E|}$, and denote their Cayley graphs as $\mG$ respectively $\mH$. We let 
\begin{equation}\label{eq:Q-cover}
\wh{\mQ}:=\mbox{ the connected component of } 1_G \mbox{ in }\varphi_{|E|}\inv(\mQ)
\end{equation} and likewise 
\[\wh{\mQ}_H:=\mbox{ the connected component of } 1_H\mbox{ in }\psi_{|E|-1}\inv(\mQ)\] 
and call these graphs the \emph{$\mQ$-cover at $1_G$ in $\mG$} and the \emph{$\mQ$-cover at $1_H$ in $\mH$}, respectively. 
For the canonical mapping $\phi\colon \mG\twoheadrightarrow \mH$ induced by the canonical morphism $G\twoheadrightarrow H$ we have $\phi(\wh{\mQ})=\wh{\mQ}_H$.

Restricting the canonical graph morphism $\varphi_{|E|}\colon \mG\twoheadrightarrow \mX_1$ to $\wh{\mQ}$, we get the canonical graph morphism $\wh{\mQ}\twoheadrightarrow \mQ$ whose restriction to the vertex set we denote as $\gamma\mapsto \ol{\gamma}$. For the next lemma,  a path $\Pi$ in an $E$-graph is a \emph{$B$-path} for some subset $B\subseteq E$ if the label 
of $\Pi$ is a word over $\til B$.
The crucial property  of the content function $\mathrm{C}$ of $G$, mentioned in Lemma~\ref{lem:crucial G}, implies the following result.

\begin{Lemma}\label{lem:contentpath} Let $\gamma,\delta$ be vertices of $\wh{\mQ}$ and $\Pi\colon \gamma\longrightarrow \delta$ be a $B$-path in $\mG$ for some $B\subseteq E$; then there exists a $B$-path $\Lambda\colon \gamma\longrightarrow \delta$ which runs in $\wh{\mQ}$.	
\end{Lemma}
\begin{proof}
Recall that the label $\ell(\underline{\phantom{\tau}})$ of every path in the following is a word over $\til E$, $E$ the generating set of $G$. 
Let $\mathrm{P}\colon \gamma\longrightarrow \delta$ be some path in $\wh{\mQ}$
that links the start and end points of the given path $\Pi \colon \gamma\longrightarrow \delta$.
    Then $\gamma[\ell(\Pi)]_G=\delta=\gamma[\ell(\mathrm{P})]_G$, whence $[\ell(\Pi)]_G=[\ell(\mathrm{P})]_G$. In particular,
	$\ell(\Pi)$ and $\ell(\mathrm{P})$ have the same $G$-content, 
    which is necessarily contained in $B$ 
    (since the word content of $\ell(\Pi)$ is contained in $B$). 
    Now $\ell(\mathrm{P})$ forms a path
	$\ol{\gamma}\longrightarrow\ol{\delta}$ in $\mQ$. By 
	Lemma~\ref{lem:crucial G}, there exists a word
	$\lambda$ over $\til E$ which forms a path
	$\ol{\gamma}\longrightarrow\ol{\delta}$ in $\mQ$ such that
	$[\ell(\mathrm{P})]_G=[\lambda]_G$ and such that $\lambda$
	uses only letters (i.e.\ edges) from the $G$-content of
	$\ell(\mathrm{P})$, which is contained in $B$. The path
	$\lambda\colon \ol{\gamma}\longrightarrow\ol{\delta}$ admits a
	unique lift $\wh{\lambda}\colon \gamma\longrightarrow\delta'$
	to $\wh{\mQ}$ for some $\delta'\in \wh\mQ$ for which
	$\ol{\delta}=\ol{\delta'}$. Note that the label
	$\ell(\wh{\lambda})$ of $\wh{\lambda}$ is just
	$\lambda$. Since $[\ell(\mathrm{P})]_G=[\lambda]_G$ we also
	have
	$\delta=\gamma[\ell(\mathrm{P})]_G=\gamma[\lambda]_G=\gamma[\ell(\wh{\lambda})]_G$;
	hence $\Lambda:=\wh{\lambda}$ is a path $\gamma\longrightarrow \delta$ and
	this path runs in $\wh{\mQ}$. Since this path uses only letter
	(i.e.~edges) from its $G$-content, it is a $B$-path as
	required.
\end{proof}

We fix some notation. For an $E$-graph $\mE$, a subset $B\subseteq E$ and a vertex $q$ of
$\mE$, $q\mE[B]$ denotes the 
\emph{$B$-component of $\mE$ containing $q$}:
the subgraph of $\mE$ spanned by all $B$-paths in $\mE$
starting at $q$. If $\varphi\colon \mE\twoheadrightarrow\mD$ is a canonical morphism of complete $E$-graphs, $q$ a vertex of $\mE$ and $B\subsetneq E$ then $\varphi(q\mE[B])=\varphi(q)\mD[B]$. For the Cayley graph $\mG$ of $G$
(being an $E$-graph), some $\gamma\in G$ and subset $K\subseteq E$, 
the graph $\gamma\mG[K]$ is the $K$-component of $\mG$
at $\gamma$; we shall refer to such graphs as \emph{$K$-coset graphs}
or just \emph{coset graphs} if the set $K$ is clear from the
context. As a $K$-graph, this is clearly isomorphic to the Cayley
graph of the subgroup $G[K]$ of $G$ generated by $K$.
Of crucial use will be Lemma~\ref{lem:freeness} below
(being an extension of Lemma~\ref{lem:contentpath}): it follows from the construction of $G$ in~\cite{ABO} but it is neither mentioned nor used in that paper (but it is a key ingredient in~\cite{otto3}: see Definition~6.13 and Lemma~6.17 there). To this end we present some more details of the construction  of $G$ in~\cite[Section 5]{ABO}. This concerns a specific component $\mU$ of $\mX_{|E|}$ (as earlier, we let $H=H_{|E|-1}$ and let $\mH$ be the Cayley graph of $H$).
\begin{Lemma}\label{lem:crucial G2}
	The graph $\mX_{|E|}$ contains a connected component $\mU$ satisfying the following conditions:
	\begin{enumerate}
		\item $\mU$ contains $\wh{\mQ}_H$ as a subgraph,
		\item for all proper subsets $B,C\subsetneq E$ and all vertices $\xi, \eta \in \wh{\mQ}_H\subseteq \mU$ we have: 
		\begin{enumerate}
			\item 	$\xi\mU[B]\cap\eta\mU[C]\ne\varnothing\Longrightarrow \xi\mU[B]\cap\eta\mU[C]\cap \wh{\mQ}_H\ne\varnothing$,
			\item 
			$\xi\mU[B]\cap\eta\mU[C]=\kappa\mU[B\cap C]$ for any $\kappa\in\xi\mU[B]\cap\eta\mU[C]$,
            \item the graphs $\xi\mU[B]\cap \wh{\mQ}_H$ and $\eta\mU[C]\cap \wh{\mQ}_H$ are connected.
		\end{enumerate}
		\end{enumerate}
\end{Lemma}
\begin{Lemma}\label{lem:crucial G3} Let $\mU$ be the graph of Lemma~\ref{lem:crucial G2} and $\vartheta\in \wh{\mQ}\subseteq \mG$; then for any $D\subsetneq E$, the canonical morphism $\varphi\colon \mG\twoheadrightarrow\mU$ (mapping $1_G$ to $1_H\in \wh{\mQ}_H$) is injective 
in restriction 
to the coset subgraph $\vartheta\mG[D]$.
\end{Lemma}
For proofs of these assertions, the reader is referred to Sections~3.3.2, 4 and 5.2, in particular, Definition~3.16, (3.8),  Theorem~4.7 and the proof of Condition~5.2, (ii) and (iii) (for $k=|E|$) in~\cite{ABO}. {In the notation of~\cite{ABO}, according to item (2) of Definition 5.3 
there, the graph $\mU$ occurs  in $\ol{\mZ_k}\subseteq \mX_{k+1}$ as the trivial completion%
\footnote{{This means that, for every letter $a$, a loop edge labelled $a$ is attached at every vertex 
of $\mathcal{V}$ that is not caintained in an $a$-cycle.}}
$\mU=\ol{\mathcal V}$  of $\mathcal{V}=\sfCE(H_k,\mC_{\mH_k};\mathbb{P}_A)$ for $k=|E|-1$, $H_k=H$, $A=E$, $\mC=\mA=\langle E\rangle=\mQ$ and $\mC_{\mH_k}=\wh{\mQ}_H$.} 
\begin{Lemma}\label{lem:freeness}
	Let $K,L\subsetneq E$ and let $\gamma,\delta$ be vertices of
	$\wh{\mQ}$.
	If $\gamma\mG[K]\cap \delta\mG[L]\ne \varnothing$ in $\mG$,
	then $\gamma\mG[K]\cap \delta\mG[L]\cap\wh{\mQ}\ne \varnothing$. Moreover, if $\varepsilon$ is a vertex in that intersection, then there exist a $K$-path $\gamma\longrightarrow \varepsilon$ and an $L$-path $\delta\longrightarrow \varepsilon$, both of which run in $\wh{\mQ}$.
\end{Lemma}

\begin{proof} 
Let $\theta\in \gamma\mG[K]\cap \delta\mG[L]$ and 
consider the canonical graph morphism $\varphi\colon
\mG\twoheadrightarrow\mU$ (mapping
$1_G$ to $1_H\in \wh{\mQ}_H$). According to Lemma~\ref{lem:crucial G3} the restrictions $\varphi_\gamma:=\varphi\upharpoonright \gamma\mG[K]$ and $\varphi_\delta:=\varphi\upharpoonright \delta\mG[L]$ are injective. Hence $\varphi_\gamma\colon \gamma\mG[K]\twoheadrightarrow \varphi(\gamma)\mU[K]$ and $\varphi_\delta\colon \delta\mG[L]\twoheadrightarrow \varphi(\delta)\mU[L]$ are graph isomorphisms. Moreover, $\varphi(\theta)\in \varphi(\gamma)\mU[K]\cap \varphi(\delta)\mU[L]$. 
In particular,
$\varphi(\gamma)\mU[K]\cap\varphi(\delta)\mU[L]\ne\varnothing$. Since
$\varphi(\gamma),\varphi(\delta)\in \varphi(\wh{\mQ})=\wh{\mQ}_H$ it follows from Lemma~\ref{lem:crucial G2} (2a) 
that $\varphi(\gamma)\mU[K]\cap\varphi(\delta)\mU[L]\cap \wh{\mQ}_H\ne
\varnothing$. Let $\xi$ be a vertex in this intersection; by Lemma~\ref{lem:crucial G2} (2b) there
exists a $(K\cap L)$-path $\Upsilon\colon \varphi(\theta)\longrightarrow
\xi$ (running in
$\varphi(\theta)\mU[K\cap L]$). Since $\varphi_\gamma$ and $\varphi_\delta$ agree on $\gamma\mG[K]\cap\delta\mG[L]=\theta\mG[K\cap L]$ (this equality holds by retractability of $G$),  their inverses $\varphi_\gamma\inv$ and $\varphi_\delta\inv$ agree on $\varphi(\theta)\mU[K\cap L]$. We then have $\varphi_\gamma\inv(\Upsilon)=\varphi_\delta\inv(\Upsilon)$ which is a $(K\cap L)$-path $\theta\longrightarrow\hat\xi:= \varphi_\gamma\inv(\xi)=\varphi_\delta\inv(\xi)$ in $\theta\mG[K\cap L]=\gamma\mG[K]\cap\delta\mG[L]$.
Since $\varphi(\gamma)\mU[K]\cap\wh{\mQ}_H$ is connected 
by Lemma~\ref{lem:crucial G2}~(2c), 
so is $\varphi_\gamma\inv(\varphi(\gamma)\mU[K]\cap\wh{\mQ}_H)$. It follows that $\gamma=\varphi_\gamma\inv(\varphi(\gamma))$ and $\hat\xi=\varphi_\gamma\inv(\xi)$ are in the same connected component of $\varphi\inv(\wh{\mQ}_H)$, that is $\hat\xi\in \wh{\mQ}$.
 Consequently, $\gamma\mG[K]\cap \delta\mG[L]\cap\wh{\mQ}\ne \varnothing$.
The second assertion of the lemma follows immediately from Lemma~\ref{lem:contentpath}.
\end{proof}

In addition to the previous facts we require a further property of $G$, namely  
\emph{$n$-acyclicity} for some positive integer $n$ (to be specified
later) which we are going to discuss now. The notions of $n$-coset cycles and $n$-acyclic groups are at the heart of~\cite{otto0,otto3}.
In the following we use again the cyclic index set $\mathbb{Z}_n$, that is, indices $i\in \mathbb{Z}_n$ are understood $\bmod\ n$.
\begin{Def}\label{def: cosetcycle} \rm Let $G$ be an $E$-generated retractable group, $B\subseteq E$, $E_1,\dots, E_n\subsetneq B$ and $\eta_1,\dots,\eta_n\in G$. The subgraph
	\[\eta_1\mG[E_1]\cup \eta_2\mG[E_2]\cup\cdots\cup \eta_n\mG[E_n]\]
	of $\mG$ is an $n$\emph{-coset cycle} provided that:
	\begin{enumerate}
		\item $\eta_{i+1}\in \eta_i\mG[E_i]$ for all $i\in\mathbb{Z}_n$,
		\item $\eta_i\mG[E_{i-1}\cap E_i]\cap \eta_{i+1}\mG[E_i\cap E_{i+1}]=\varnothing$ for all $i\in\mathbb{Z}_n$.
	\end{enumerate}
\end{Def}
This notion is reminiscent of the notion of an $n$-cyclic configuration
in Section~\ref{sec:groups for Ash}. In fact, $n$-coset cycles are
special types of $n$-cyclic configurations. The context, however,
is somehow different. In Section~\ref{sec:groups for Ash} we had a
group $G$ and we wanted to construct an {expansion} $H$ of $G$ that
avoids a given $n$-cyclic configuration (or even all of them)
\textsl{of its quotient} $G$.
In the present context we are looking for a group $H$ which avoids all $n$-coset cycles \textsl{in its own} Cayley graph, which is a more ambitious task.%
\footnote{{This is somehow analogous to the fact, discussed in Section~\ref{sec:strictly stronger}, that the series of expansions used to prove Proposition~\ref{prop:ash-pure} is always ``one step away'' from delivering a proof of Proposition~\ref{prop:otto-ash-pure}.}}

{The $n$-coset cycle of $\mG$ in Definition~\ref{def: cosetcycle} is determined by the $n$-tuples
\[(\eta_1,\dots,\eta_n)\mbox{ and }(E_1,\dots, E_n).\]
For every $i$, $\delta_i:=\eta_i\inv \eta_{i+1}\in G[E_i]$ (the subgroup of $G$ generated by $E_i$) and Definition~\ref{def: cosetcycle} implies that $\delta_1\cdots \delta_n=1_G$.}

\begin{Def}\label{def:N-acyclicity} Let $n>1$ be a positive integer; an $E$-generated retractable group $G$ is $n$\emph{-acyclic} if the Cayley graph $\mG$ does not contain any coset cycles of length up to $n$.
\end{Def}
Note that $2$-acyclicity of $G$ is the same as $G[B\cap C]=G[B]\cap G[C]$ for all subsets $B,C\subsetneq E$. This just says that two coset subgraphs $g\mG[B]$ and $h\mG[C]$ of $\mG$ either have empty intersection or otherwise $g\mG[B]\cap h\mG[C]=k\mG[B\cap C]$ for any $k\in g\mG[B]\cap h\mG[C]$; this condition is guaranteed by retractability of $G$.
We briefly sketch how we can achieve that the group $G$ in question be $n$-acyclic for any given (but fixed) positive integer $n$. 

Let $B\subseteq E$, $B_1,\dots,B_n\subsetneq B$ and let $H$ be an $E$-generated group with Cayley graph $\mH$ {and assume that $H[B]$ is retractable}; for every $i\in \mathbb{Z}_n$ choose a vertex $\gamma_i\in \mH[B_i]$ such that
\begin{equation}\label{eq:amalgamchain}
\mH[B_{i-1}\cap B_i]\cap\gamma_i\mH[B_i\cap B_{i+1}]=\varnothing.    
\end{equation}
Note that $\mH[B_{i-1}\cap B_i]\cup\gamma_i\mH[B_i\cap B_{i+1}]$ is a subgraph of $\mH[B_i]$.
Now take the disjoint union 
$\bigsqcup_{i=1}^n \mH[B_i]$
and factor by the smallest $E$-graph congruence $\Theta$ which
identifies $\gamma_i\in\mH[B_i]$ with $1\in \mH[B_{i+1}]$ for
$1\le i\le n-1$. {Condition~\eqref{eq:amalgamchain} guarantees that $\gamma_n\in \mH[B_n]$ and $1\in \mH[B_1]$ are not identified under $\Theta$.} Call the resulting graph a 
$\pss B$-\emph{amalgamation chain of length $n$ for the group $H$}.
For $n\ge 6$,  $n$-acyclicity  of $G$ is achieved by adding in Definition 5.3 of the graph $\mZ_k$ in~\cite{ABO} as a third item (3) all $\pss A$-amalgamation chains of length up to $n-2$ for the group $H_k$ in the series \eqref{eq: expanding groups} for all subsets $A\subseteq E$ of size $|A|=k+1$ (for more details see~\cite{otto0} and~\cite{otto3}).\footnote{The group $G$ as it is defined in~\cite{ABO} is already $5$-acyclic since the augmented clusters in item (1) of Definition 5.3 contain all amalgamation chains of length up to $3$.} 

One can check that for $C\subsetneq B$, every $C$-component of a 
$\pss B$-amalgamation chain of length $n$ is either a $C$-coset $\gamma\mH[C]$ or 
a $\pss\,C$-amalgamation chain of length at most $n$ (see~\cite[Section 2]{otto0} or~\cite[Lemma 4.8]{otto3}). 
{From this it follows that for every $k$, the expansion $H_k\twoheadleftarrow G_{k+1}$ is \emph{$k$-stable}, that is, for every $k$-element subset $B\subseteq E$, $G_{k+1}[B]\cong H_k[B]$ via restriction of the canonical morphism $G_{k+1}\twoheadrightarrow H_k$; this is an essential ingredient of the series~\eqref{eq: expanding groups}.} 

{It should be mentioned that the procedure 
in~\cite{otto0} was to add, at every step, amalgamation chains of length up to $n$ in order to get an $n$-acyclic group, while in~\cite[Lemma 5.1]{otto3} it is pointed out that amalgamation chains of length up to $n$ are sufficient to gain an $(n+2)$-acyclic group at the end.} 
{We give a sketch how this is proved
                     in~\cite{otto3}. Suppose the
sets 
$E_0, E_1,\dots, E_{n+1}\subsetneq E$ and elements $\gamma_0,\gamma_1,\dots, \gamma_{n+1}$ 
would give rise to 
an $(n+2)$-coset cycle in the final group $G=G_{|E|}$. For $i\in\mathbb{Z}_{n+2}$ let $\delta_i:=\gamma_i\inv \gamma_{i+1}$ and $\eta_i:=\varphi(\delta_i)\in H$ where $\varphi$ denotes the canonical morphism $G\twoheadrightarrow H$ and  $H:=H_{|E|-1}$ is the penultimate group in the series~\eqref{eq: expanding groups}. Recall that 
\[\gamma_0\mG[E_0]\cup \gamma_1\mG[E_1]\cup\cdots\cup \gamma_{n+1}\mG[E_{n+1}\]
being an $(n+2)$-coset cycle of $G$ implies that $\delta_0\cdots \delta_{n+1}=1_G$.
Denote by $\mH$ and $\mH[E_i]$ the Cayley graphs of $H$ and $H[E_i]$, respectively. The data $((\mH[E_i],\eta_i)\colon i=1,\dots n)$ define an $\pss E$-amalgamation chain $\mC$ of length $n$. According to the definition, $\mC=\bigsqcup_{i=1}^n \mH[E_i]\,/\,\Theta$ for the congruence $\Theta$ defined above. For every $i=1,\dots, n$ let $\mH_i:=\mH[E_i]\Theta$, (essentially the realisation of $\mH[E_i] $ within $\mC$). Note that $\mC=\mH_1\cup\cdots\cup \mH_n$ and $\mH_i\cong \mH[E_i]$ for all $i$. Let $\varepsilon:=1_{\mH[E_1]}\Theta$  (the $1$-element of $\mH[E_1]$ within $\mC$) and $\omega:=\eta_n\Theta\in \mH_n$ (the element $\eta_n$ of $\mH[E_n]$ within $\mC$). }

{By our assumption, $\mC$ belongs to $\mZ_{|E|-1}$, hence there is a canonical graph morphism $\mG\twoheadrightarrow \ol\mC$ with $\ol\mC$ the trivial completion of $\mC$. Now, for every $i\in\mathbb{Z}_{n+2}$, take a word $w_i\in \til {E_i}^*$ such that $[w_i]_G=\delta_{i}$. We look at actions of the words $w_i$ in the complete graph $\ol\mC$. Firstly, let $\tau:=\varepsilon w_0\inv$; then, in $\mC$,  $\tau=\varepsilon u\inv$ where $u$ is the word obtained from $w_0$ by deleting all letters not in $E_1$. Hence
$\tau$ belongs to the  $(E_0\cap E_1)$-component of $\varepsilon$ in $\mH_1$, that is,  $\tau\in \varepsilon\mH_1[E_0\cap E_1]$.  In $\mC$, and therefore also in $\ol\mC$, we further have $\varepsilon w_1\cdots w_n=\omega \in \mH_n$.  Let $v$ be the word obtained from $w_{n+1}$ by deleting all letters not in $E_n$; then $\omega w_{n+1}=\omega v$ (action to be taken in $\ol\mC$, 
although
$\omega v$ is defined even in $\mC$) belongs to the $(E_n\cap E_{n+1})$-component of $\omega$ in $\mH_n$, that is, $\omega w_{n+1}\in \omega\mH_n[E_n\cap E_{n+1}]$. Since $\varepsilon\mH_1[E_0\cap E_1]\cap \omega\mH_n[E_n\cap E_{n+1}]=\varnothing$ we get $\tau\ne \omega w_{n+1}$. Altogether this leads to 
\[\tau w_0w_1\cdots w_{n+1}\ne \tau\]
which holds in $\ol\mC$. Since $\mG\twoheadrightarrow \ol\mC$ we get $[w_0\cdots w_{n+1}]_G\ne 1_G$, 
contradicting 
$\delta_0\cdots\delta_{n+1}=1_G$. Consequently, 
$E_0,\dots, E_{n+1}$ and $\gamma_0,\dots, \gamma_{n+1}$ cannot give rise to an $(n+2)$-coset cycle in $G$. The situation is depicted in Figure~\ref{fig:amalgamationC}.}
\begin{figure}[ht]
\begin{tikzpicture}[scale=1]
	\draw plot [smooth cycle] coordinates {(-6,1.3) (-3,1.3) (-3,-1.3) (-6,-1.3)};
	\draw plot [smooth cycle] coordinates {(-4.5,1) (-1.5,1) (-1.5,-1) (-4.5,-1)};
	\draw(0,0) node{$\cdots$};
	\draw plot [smooth cycle] coordinates {(6,1.2) (3.2,1.2) (3.2,-1.2) (6,-1.2)};
	\draw plot [smooth cycle] coordinates {(4.5,1) (1.5,1) (1.5,-1) (4.5,-1)};
	\filldraw (-5.2,0) circle (1pt);
	\filldraw (-3.2,0) circle (1pt);
	\filldraw (-1.5,0) circle (1pt);
	\filldraw (5.2,0) circle (1pt);
	\filldraw (3.2,0) circle (1pt);
	\filldraw (1.5,0) circle (1pt);
	\filldraw(-5.7,1)circle(1pt);
	\filldraw(5.7,1)circle(1pt);
	\draw[latex-](-5.2,0) to (-5.7,1);
	\draw[-latex](-5.2,0) to (-3.2,0);
	\draw[-latex](-3.2,0) to (-1.5,0);
	\draw[dotted](-1.5,0) to(-0.5,0);
	\draw[dotted,-latex](0.5,0)to(1.5,0);
	\draw[-latex](1.5,0)to(3.2,0);
	\draw[-latex](3.2,0)to(5.2,0);
	\draw[-latex](5.2,0)to(5.7,1);
	\draw(-5.2,0)node[below]{$\varepsilon$};
	\draw(5.2,0)node[below]{$\omega$};
	\draw(-5.65,0.95)node[above]{$\tau$};
	\draw(-5.3,0.6)node{$u$};
	\draw(5.25,0.5)node{$v$};
	\draw(-4.2,-0.1)node[above]{$w_1$};
	\draw(4.2,-0.1)node[above]{$w_n$};
	\draw(2.2,-0.1)node[above]{$w_{n-1}$};
	\draw(-2.3,-0.1)node[above]{$w_2$};
	\draw(-5.7,-1)node{$\mathcal{H}_1$};
	\draw(5.6,-1)node{$\mathcal{H}_n$};
	\draw(-1.75,-0.7)node{$\mathcal{H}_2$};
	\draw(1.95,-0.7)node{$\mathcal{H}_{n-1}$};
\end{tikzpicture}
\caption{The amalgamation chain $\mC$}\label{fig:amalgamationC}
\end{figure}

\subsubsection{The case of $\mQ$ being the Cayley graph of an $A$-generated group $Q$}\label{sec:group G}\label{sec:Cayleygraphcase}
We now assume that the oriented graph $\mQ$ is the Cayley graph of a finite $A$-generated group $Q$ for some finite alphabet $A$. Recall that the set of vertices  of $\mQ$ is $Q$ and the set of positive edges $E$ of $\mQ$ is $Q\times A$ where 
\begin{equation}\label{eq:edges of Q}
\alpha(q,a)=q,\ \omega(q,a)=q[a]_Q\mbox{ and }(q,a)\inv=(q[a]_Q,a\inv).
\end{equation} 
The last equality in \eqref{eq:edges of Q} implies the equality of  sets $(Q\times A)\inv = Q\times A\inv$ and hence $\til E=Q\times \til A$ (this is strictly true and not just a manner of speaking).

Let $G$ be the $E$-generated group as discussed in the previous subsection, but for $\mQ$ the Cayley graph of a finite $A$-generated group $Q$. 
We shall use the Cayley graph $\mG$ of $G$, which has vertex set $G$ and edge set $G\times \til E$. The edges of $\mG$ therefore are of the form $(\gamma, (q,a))$ with $\gamma\in G$,  $(q,a)\in \til E$ and
\[(\gamma, (q,a))\inv=(\gamma[(q,a)]_G,(q[a]_Q,a\inv)).\]
Recall~\eqref{eq:Q-cover}, the definition of $\wh{\mQ}$,  the $\mQ$-cover 
at $1_G$ in $\mG$ (in the previous subsection $1_Q$ was a vertex of $\mQ$ chosen arbitrarily, now we let $1_Q$ be the identity element of the group $Q$). 
The vertex set of $\wh \mQ$ is denoted $\wh
Q$. The surjective mapping $\wh{Q}\twoheadrightarrow Q$ induced by the
canonical morphism $\wh{\mQ}\twoheadrightarrow \mQ$ will be denoted
$\gamma\mapsto\ol{\gamma}$ (as in the previous subsection).

\begin{Rmk} The (canonical) mapping $\wh{Q}\twoheadrightarrow Q$ can be extended to a surjective mapping $G\to Q$ by restricting the canonical graph morphism $\mG\twoheadrightarrow\mX_1$ to its vertex set and taking into account that the vertex set of $\mX_1$ is $Q$. This mapping $G\to Q$  neither is a group homomorphism nor does it give rise to a graph homomorphism $\mG\to \mQ$, and in this paper we will use only its restriction to $\wh{Q}$.
\end{Rmk}

In the following we shall deal with $A$-graphs (related to subgraphs of $\mQ$) as well as $E$-graphs (related to subgraphs of $\mG$). 
The $E$-graph $\wh\mQ$ is the subgraph of $\mG$  spanned by all paths in $\mG$ starting at $1_G$ and whose labels (being words over $\til E$) form paths in $\mQ$ starting at $1_Q$. 
Every such path in $\mQ$ is \textsl{induced} by a word over $\til A$: for every word $x\in \til A^*$ and every $q\in Q$ there exists a unique path in $\mQ$ with initial vertex $q$ whose label is $x$ --- we denote this path by $\pi_q^\mQ (x)$ and say this path is \emph{induced} by the word $x$.  We give a short summary of the situation.
\begin{Sum}\rm Let
	\[\Pi=(\gamma_1,(q_1,a_1))\cdots(\gamma_n,(q_n,a_n))\]
	be a word over $G\times \til E$. Then $\Pi$ forms a path $\gamma\longrightarrow\delta$ in $\wh{\mQ}$ if and only if its label 
	\[\ell(\Pi)=(q_1,a_1)\cdots(q_n,a_n)\] (a word over $\til E$) forms a path $\ol{\gamma}\longrightarrow\ol{\delta}$ in $\mQ$ induced by the word $x=a_1\dots a_n$, that is 
	\[\ell(\Pi)=\pi_{\ol{\gamma}}^\mQ(x).\]
	Moreover, in this case, $\gamma=\gamma_1$ and $\ol{\gamma}=\ol{\gamma_1}=q_1$. The canonical graph morphism $\wh{\mQ}\twoheadrightarrow\mQ$ maps the path $\Pi$ onto the path $\ell(\Pi)=\pi_{\ol{\gamma}}^\mQ(x)$. The original path $\Pi$  can be seen as the path in $\mG$ starting at $\gamma$ induced by the word $\ell(\Pi)$, and this path actually runs in $\wh{\mQ}\subseteq \mG$. 	
\end{Sum}

By~\cite[Proposition 5.5]{ABO} the construction of $G$ guarantees that every bijection $\varphi$ of $E$ that is induced by an automorphism of the oriented graph $\mQ$ extends to an automorphism $\wh{\varphi}$ of $G$. In particular, this is the case for every automorphism $\varphi_q$ of $\mQ$ induced by  multiplication by the element $q$ on the right, for any $q\in Q$. Hence,  for $E=Q\times A$, the permutation $\varphi_q\colon E\to E$, $(g,a)\mapsto {}^q(g,a):=(qg,a)$ extends to  an automorphism of $G$. 
It follows that the group $Q$ acts on the group $G$ by automorphisms on the left via $\gamma\mapsto {}^{q}\gamma$ for $\gamma\in G, q\in Q$. So we can form the semidirect product $G\rtimes Q$.

In~\cite{ABO} the group $H$ (depending on $G$ and $Q$) is defined as: 
\begin{equation}\label{eq:group H}
	H:=\langle([(1_Q,a)]_G,[a]_Q)\colon a\in A \rangle\le G\rtimes Q,
\end{equation}
which is an $A$-generated group. 
(Compare Section~2.5 in~\cite{ABO} and especially equation~(2.3) there.)
For a word $x$ over $\til A$, we have $[x]_H=([\pi_1^\mQ(x)]_G,[x]_Q)$ (see (2.4) in~\cite{ABO}).

In a sense, the group $H$ is sitting inside $G$. More precisely, the
vertices of $\wh{\mQ}$ are in bijective correspondence with the
elements of $H$ via the map $\gamma\mapsto
(\gamma,\ol\gamma)$. Indeed, take any word $x\in \til A^*$; then the
terminal vertex of the path $\pi_1^\mQ(x)$ is $[x]_Q$. Being  a
word over $\til E$, $\pi_1^\mQ(x)$ induces a path
$\Pi$ in $\mG$ starting at $1_G$ running in $\wh{\mQ}$. Its terminal
vertex is $[\ell(\Pi)]_G=[\pi_1^\mQ(x)]_G$. The canonical graph
morphism $\wh{\mQ}\twoheadrightarrow \mQ$ maps $\Pi$ onto
$\ell(\Pi)=\pi_1^\mQ(x)$, and hence maps the terminal vertex $[\ell(\Pi)]_G=[\pi_1^\mQ(x)]_G$ of $\Pi$ onto the terminal vertex $[x]_Q$ of $\pi_1^\mQ(x)$. Consequently, the element $([\pi_1^\mQ(x)]_G,[x]_Q)$ is of the form $(\gamma,\ol\gamma)$ for $\gamma=[\pi_1^\mQ(x)]_G$, that is, the mapping $\wh Q\to H$, $\gamma\mapsto (\gamma,\ol{\gamma})$ is surjective. Clearly,  for two distinct vertices $\gamma,\eta\in \wh{\mQ}$ we have $(\gamma,\ol\gamma)\ne(\eta,\ol\eta)$. 
As a consequence we have the following.
\begin{Lemma}\label{lem:wordproblemH} For two words $x,y\in \til A^*$, $[x]_H=[y]_H$ if any only if $[\pi_q^\mQ(x)]_G=[\pi_q^\mQ(y)]_G$ for some (any) $q\in Q$.
\end{Lemma}
This result also shows that the group $H$ could equally well be defined as the quotient $H:=F/N$ where $F$ is the $A$-generated free group and \[N=\{[w]_F\colon w\in \til A^*,
[\pi_1^\mQ(w)]_G=1_G\}.\] 

Without going too much into the details, we mention that the
``equivalence'' of $\wh{Q}$ and $H$ admits a groupoidal 
interpretation. Since the generating set of the group $G$ is the set
$E$ of positive edges of the graph $\mQ$ we can as well consider the
\emph{groupoid $\mathrsfs{Q}$ generated by the graph $\mQ$ subject to
the relations of $G$} defined as follows: the vertices
(i.e.\ objects) of $\mathrsfs{Q}$ are the vertices of $\mQ$, that is
the elements of $Q$. For $p,q\in Q$ the arrows (i.e.\ morphisms)
$p\longrightarrow q$ are all $G$-values of paths $p\longrightarrow q$
in $\mQ$  (recall that every such path is a word over $\til E$):
\[\mathrsfs{Q}(p,q)=\{[\pi]_G\colon \pi\mbox{ is a path } p\longrightarrow q\mbox{ in }\mQ\}.\] 
The $G$-value $ [\pi]_G$ of every path $\pi$ in $\mQ$ starting at $1_Q$ is the terminal vertex $\omega\wh{\pi}$ of the lift $\wh{\pi}$ of $\pi$ in $\mG$ starting at $1_G$. Hence the set
\begin{equation}\label{eq:rep of H}
	\mathrsfs{H}:=\bigsqcup_{q\in Q}\mathrsfs{Q}(1_Q,q)
\end{equation}
of arrows of $\mathrsfs{Q}$ with initial vertex $1_Q$ is in bijective correspondence with the set of vertices of $\wh{\mQ}$ that is, with the elements of $\wh{Q}$. We can define a group structure on this set  of arrows of $\mathrsfs{Q}$ (and hence on $\wh{Q}$) as follows. For two such arrows  $\gamma\in\mathrsfs{Q}(1_Q,\ol{\gamma})$, $\delta\in \mathrsfs{Q}(1_Q,\ol{\delta})$ we define \[\gamma\delta:=\gamma\cdot{}^{\ol\gamma}\delta\]
where the multiplication $\cdot$ on the right hand side is the multiplication in the groupoid $\mathrsfs{Q}$ and
${}^{\ol\gamma}\delta$ is the arrow  $\ol\gamma\longrightarrow \ol
\gamma \ol \delta$ in $\mathrsfs{Q}$ obtained by ``shifting'' the
arrow $\delta$ from $1_Q\longrightarrow \ol\delta$ to $\ol
\gamma\longrightarrow\ol\gamma \ol\delta$. To make this more precise, recall that $Q$ acts on its Cayley graph $\mQ$ by left  multiplication: $g\mapsto {}^rg:=rg$ and $(g,a)\mapsto {}^r(g,a):=(rg,a)$ for $r\in Q$, $g$ a vertex of $\mQ$ and $(g,a)$ and egde of $\mQ$. This action extends to an action by groupoid automorphisms of $\mathrsfs{Q}$ in a natural way: on the set of vertices (i.e.~objects) the action is the same as for $\mQ$; on the set of arrows $[\pi]_G\colon p\longrightarrow q$ it becomes ${}^r[\pi]_G\colon rp\longrightarrow rq$ where ${}^r[\pi]_G=[{}^r\pi]_G$ and ${}^r\pi={}^re_1\cdots {}^re_n$ for $\pi=e_1\cdots e_n$ with $e_i\in \til E$. So the arrow ${}^{\ol{\gamma}}\delta\colon \ol{\gamma}\longrightarrow \ol{\gamma}\ol{\delta}$ is just the result of the application of the automorphism ${{x}}\mapsto {}^{\ol{\gamma}}{{x}}$  to the arrow $\delta\colon 1_Q\longrightarrow \ol{\delta}$. 
The resulting group is $A$-generated, namely by (the $G$-values of) the edges $1_Q\overset{a}{\longrightarrow}[a]_Q$, more precisely by the set  $\{[(1_Q,a)]_G\colon a\in A\}$, and  is isomorphic to the group $H$. 

In the notation of  Margolis and Pin~\cite{MarPin}, the group so obtained is exactly  $\mathrsfs{Q}/Q$, see Proposition 3.11 in~\cite{MarPin}. Here $\mathrsfs{Q}/Q$ denotes the set of orbits of the action of $Q$ on $\mathrsfs{Q}$ on which the groupoid structure of $\mathrsfs{Q}$ induces a (quotient) groupoid structure which turns out to be a group since the transitivity of the action of $Q$ on the set of vertices (i.e.~objects) of $\mathrsfs{Q}$ forces $\mathrsfs{Q}/Q$ to have only one vertex (i.e.~one object). The set $\mathrsfs{H}$ in~\ref{eq:rep of H} is a complete set of representatives of the orbits,  hence may be used  to define the quotient structure. Altogether, the groupoidal interpretation just discussed may be summarised in one sentence: \textsl{ignore the (redundant) second components $\ol\gamma$ in the elements $(\gamma,\ol{\gamma})\in H$, but compose the first components $\gamma$ correctly}.

\subsection{Proof of Theorem~\ref{thm:otto-ash}}\label{sec:proofs}
{The task of this section is to prove Theorem~\ref{thm:otto-ash}, with Lemma~\ref{lem:main lemma} being the main achievement.} While inverse monoids are seemingly absent in the lemma, the connection with inverse monoids will become clear in the discussion after its proof.

We let $Q$ be a finite $A$-generated group with Cayley graph $\mQ$ and
consider the $A$-generated group $H$ defined in \eqref{eq:group H}, where $G$ is the $E$-generated group, for $E=Q\times A$, as discussed in Sections~\ref{sec:generalcase} and~\ref{sec:group G}. In particular $G$ satisfies the conditions of Lemmas~\ref{lem:crucial G} and~\ref{lem:freeness} and we  also 
assume that $G$ is $n$-acyclic in the sense of Definition~\ref{def:N-acyclicity}
for some positive integer $n$ --- the existence of such a group has been sketched in the last two paragraphs of Section~\ref{sec:generalcase}. Note 
that $G$ and $H$ are uniquely determined by the input data $Q$ and $n$ according to the construction outlined in Sections~\ref{sec:generalcase} and~\ref{sec:Cayleygraphcase}
and we set $H=H(Q,n)$.
Recall the word content $\co(\underline{\phantom{v}})$ of a word over $\til E$ according to Definition~\ref{def:word content} and that, for a word $w$ over $\til A$ and $g\in Q$, $\pi_1^\mQ(w)$ 
and $\pi_g^\mQ(w)$ denote the paths in $\mQ$ labelled $w$, starting at $1_Q$ and at $g$ respectively, considered as words over $\til E$. 

\begin{Lemma}[main lemma]\label{lem:main lemma}
Let $Q$ be a finite $A$-generated group, $n$ a positive integer and $H=H(Q,n)$. Then for any finite connected graph 
$\mE=(V\cup K;\alpha,\omega,\inv)$  on at most $n$ vertices and without loop edges  
and for any $H$-commuting word labelling $v\colon K\to \til A^*$  there exist a finite
tree $\mT=(W\cup L;\alpha,\omega,\inv)$ with vertex set $W\supseteq
V$ and a word labelling $w\colon L\to \til A^*$ such that: for every edge
$k\in  K$, if $\pi_k$ is the  reduced path $\alpha k\longrightarrow
\omega k$ in $\mT$ and $w_k=w(\pi_k)$, then 
\begin{enumerate}
    \item $[v(k)]_H=[w_k]_H$ and
    \item $\co(\pi_1^\mQ(v(k))\supseteq \co(\pi_1^\mQ(w_k))$.
\end{enumerate} 
In addition, $\mT$ can be chosen such that all vertices of $W\setminus V$ have degree at least $3$.
\end{Lemma} 

\begin{proof}
The proof of this lemma is split into several cases,
reflecting increasing levels of complexity with respect to the graph structure, 
in \S~\ref{subsubsection:2vertices}--~\ref{subsubsection:gengraphs}. 
In the following,  elements of $Q$ will always be denoted by $g,h, g_i, h_i$, etc., those of $G$ by $\gamma, \eta, \gamma_i,\eta_i$, etc., and those of $H$ by $(\gamma,g)$, $(\gamma_i,g_i)$, etc.

\subsubsection{Graphs with two vertices} 
\label{subsubsection:2vertices}
Let us consider first all
graphs  having only two vertices, say $o_1,o_2$ and $t>1$ geometric
edges between these vertices, and let  $k_1,\dots, k_t\colon
o_1\longrightarrow o_2$. This case includes cycle graphs of length
$2$. Let $v_i := v(k_i)$; $H$-commutativity of the labelling just says
that $[v_1]_H=\cdots= [v_t]_H$, that is,
\[([\pi_1^\mQ(v_i)]_G,[v_i]_Q)=[v_i]_H=[v_j]_H=([\pi_1^\mQ(v_j)]_G,[v_j]_Q)\]
for all $i,j$.
In particular, the $Q$-values of all words $v_i$ coincide, say $q:=[v_i]_Q$ for all $i$. Also, the $G$-values of all sequences $\pi_1^\mQ(v_i)$ coincide, say $\gamma:=[\pi_1^\mQ(v_i)]_G$ for all $i$.  A fortiori, all sequences $\pi_1^\mQ(v_i)$ have the same
$G$-content, say $E_0:=\mathrm{C}(\gamma)$; also note that $E_0\subseteq \mathrm{co}(\pi_1^\mQ(v_i))$ for all $i$ (the word content always contains the $G$-content).
According to Lemma~\ref{lem:crucial G} there exists a word $\pi$ over $ \til{E_0}$ forming a path $\pi\colon 1_Q\longrightarrow q$ in $\mQ$ such that $[\pi]_G=\gamma$.
Since $\pi$ forms a path in $\mQ$ it is induced by a word $w$ over $\til A$, that is, $\pi=\pi_1^\mQ(w)$ for some word $w\in \til A^*$. Then $\mathrm{co}(\pi_1^\mQ(w))=E_0\subseteq \mathrm{co}(\pi_1^\mQ(v_i))$ for all $i$ and $[\pi_1^\mQ(w)]_G=\gamma$. 
From Lemma~\ref{lem:wordproblemH} we get $[w]_H=[v_1]_H=\cdots=[v_t]_H$. 
A suitable tree $\mT$ for the given graph labelling consists of the two vertices $o_1,o_2$, an edge $f\colon o_1\longrightarrow o_2$ and its inverse $f\inv $, and the labelling $w(f):= w$.

\subsubsection{Cycle graphs of length $3$} 
\label{subsubsection:3cycles}
We continue with cycle graphs of length $3$. Let $\mC=(V\cup {K};\alpha,\omega,\inv)$ be a cycle graph of length $3$; let $V=\{o_1,o_2,o_3\}$ and $K=\{k_1^{\pm 1},k_2^{\pm 1},k_3^{\pm 1}\}$ with $k_i\colon o_i\longrightarrow o_{i+1}$ (indices to be taken $\bmod\ 3$). Let $v\colon K\to \til A^*$
be an $H$-commuting labelling, given by $v(k_i)=v_i$. 
For $i=1,2,3$, let $[v_i]_H =:(\gamma_i,g_i)$ so that  $g_i=[v_i]_Q$ and
$\gamma_i=[\pi_1^\mQ(v_i)]_G$.
We set $h_i$ and $\eta_i$ according to 
\[
 \left.\begin{array}{r@{\;:=\;}l@{\qquad\quad}r@{\;:=\;}l}
 h_1 & 1_Q
 &
 \eta_1 & 1_G
 \\
 h_2 & g_1
 &
 \eta_2 & \gamma_1
\\
 h_3 & g_1g_2
 &
 \eta_3 & \gamma_1\, {}^{g_1}\!\gamma_2
 \end{array}\right\}
 \mbox{ so that }
 \left\{
 \begin{array}{r}
 (\eta_i,h_i)=[v_1\cdots v_{i-1}]_H \\ \mbox{and } h_i=\ol{\eta_i}
 \mbox{ for } 1\leq i \leq 3. 
 \end{array}
 \right.
 \]

For $i=1,2,3$ let $E_i\subseteq E$ be the content of ${}^{h_i}\gamma_i$, which in fact is the $G$-content of $\pi_{h_i}^\mQ(v_i)$. Since $G$ is $3$-acyclic, the subgraph
\[\eta_1\mG[E_1]\cup\eta_2\mG[E_2]\cup\eta_3\mG[E_3]\] 
of $\mG$ is not a coset cycle. Since $\eta_{i+1}=\eta_i{}^{h_i}\gamma_i\in \eta_i\mG[E_i]$ (indices to be taken $\bmod\ 3$, that is, $1_G=\eta_1=\eta_3{}^{h_3}\gamma_3\in \eta_3\mG[E_3]$), $3$-acyclicity implies that
\begin{equation}\label{eq:3-acyclic}
\eta_1\mG[E_1]\cap\eta_2\mG[E_2]\cap \eta_3\mG[E_3]\ne\varnothing.
\end{equation}

Since $\eta_1\in \eta_3\mG[E_3]$ we have
$\eta_3\mG[E_3]=\eta_1\mG[E_3]$, and
\eqref{eq:3-acyclic} can be written as
\begin{equation}\label{eq:3acyclic2}\eta_1\mG[E_1\cap
	E_3]\cap\eta_2\mG[E_2]=
\eta_1\mG[E_1]\cap \eta_2\mG[E_2]\cap\eta_1\mG[E_3]
\ne \varnothing.
\end{equation}
Here we use that retractability of $G$ implies $G[E_1]\cap
G[E_3]=G[E_1\cap E_3]$. Let $\wh{\mQ}$ be the cover of $\mQ$ at $1_G$ in $\mG$. Since $1_G=\eta_1$ and $\eta_2$ belong to $\wh{\mQ}$
we have by Lemma~\ref{lem:freeness} that 
\[\eta_1\mG[E_1\cap E_3]\cap \eta_2 \mG[E_2]\cap\wh{\mQ}\ne\varnothing.\]
\begin{figure}
\begin{tikzpicture}[xscale=2,yscale=2.3]
	\filldraw(0,0) circle(1pt);
	\filldraw(0,-1)circle(1pt);
	\filldraw(0.866025,0.5) circle(1pt);
	\filldraw(-0.866025,0.5)circle(1pt);
	\draw plot [smooth cycle] coordinates {(-0.2,0.2)(-0.2,-1.5)(1.2,-0.8)(1.2,0.7)};
	\draw plot [smooth cycle] coordinates {(0.2,0.2)(0.2,-1.5)(-1.2,-0.8)(-1.2,0.7)};
	\draw plot [smooth cycle] coordinates{(1.3,0.8)(0,-0.2)(-1.3,0.8)(0,1.5)};
	\filldraw[ gray, nearly transparent,rounded corners] plotcoordinates {(-0.2,-1.5)(-1,-0.5)(-2.1,0.75)(-0.6,1)(-0.4,1.7)(0.4,1.7)(0.6,1)(2.1,0.75)(1,-0.5)(0.1,-1.5)};
	\draw[below](0,-1.05)node{$\eta_1$};
	\draw[above](0,0)node{$\delta$};
	\draw[above left](-0.806025,0.5)node{$\eta_2$};
	\draw[above right](0.806025,0.5)node{$\eta_3$};
	\draw(0,0.8)node{$\widehat{\mathcal Q}$};
	\draw[thick,-latex](0,-1) to (0,-0.05);
	\draw[thick,-latex](-0.866025,0.5) to (-0.05,0.05);
	\draw[thick,-latex](0.866025,0.5) to (0.05,0.05);
	\draw[thick,dotted,-latex](0,-1).. controls (-0.5,-1) and (-3,1.6)..  (-0.9,0.5) ;
	\draw[thick, dotted,rounded corners,-latex](-0.866025,0.5)--(-0.4,1)--(-0.3,1.6)--(0.3,1.6)--(0.4,1)-- (0.866025,0.53);
	\draw[thick, dotted,-latex](0.866025,0.5)..controls (3,1.6) and (0.5,-1).. (0.03,-1);
	\draw[right](0,-0.4)node{$\Pi_1$};
	\draw(-0.35,0.35)node{$\Pi_2$};
	\draw(0.3,0.35)node{$\Pi_3$};
	\draw(0,1.1)node{$\eta_2\mathcal{G}[E_2]=\eta_3\mathcal{G}[E_2]$};
	\draw(-0.9,0)node{$\eta_2\mathcal{G}[E_1]=$};
	\draw(-0.6,-0.3)node{$\eta_1\mathcal{G}[E_1]$};
	\draw(0.88,0)node{$\eta_3\mathcal{G}[E_3]=$};
	\draw(0.6,-0.3)node{$\eta_1\mathcal{G}[E_3]$};
\end{tikzpicture}
\caption{Transforming a $3$-cycle into a tree}\label{fig:3cycle}	
\end{figure}%
Let $\delta$ be
some vertex in this intersection, see Figure~\ref{fig:3cycle}.

By Lemma~\ref{lem:freeness} there are paths $\Pi_i$ running in $\wh{\mQ}$, 
$\Pi_i\colon\eta_i\longrightarrow \delta$, such that $\Pi_1$ runs in $\eta_1\mG[E_1\cap E_3]$, $\Pi_2$ runs in $\eta_2\mG[E_1\cap E_2]$ and $\Pi_3$ runs in $\eta_3\mG[E_2\cap E_3]$. Since the paths $\Pi_i$ run in $\wh{\mQ}$ they are induced by words over $\til A$: there are $x_1,x_2,x_3\in \til A^*$ such that $\ell(\Pi_i)=\pi_{h_i}^\mQ(x_i)$. From the definition of the paths $\Pi_i$ it follows that
\begin{align*}
{}^{h_1}\gamma_1&=[\ell(\Pi_1\Pi_2\inv)]_G=[\pi_{h_1}^\mQ(x_1x_2\inv)]_G\\
{}^{h_2}\gamma_2&=[\ell(\Pi_2\Pi_3\inv)]_G=[\pi_{h_2}^\mQ(x_2x_3\inv)]_G\\
{}^{h_3}\gamma_3&=[\ell(\Pi_3\Pi_1\inv)]_G=[\pi_{h_3}^\mQ(x_3x_1\inv)]_G
\end{align*}
where $h_1=1_Q$. Figure~\ref{fig:3cycle} depicts the situation: the
shaded area represents the $\mQ$-cover $\wh{\mQ}$ and the dotted
arrows indicate the paths with labels $\pi_{h_i}^\mQ(v_i)$.
From Lemma~\ref{lem:wordproblemH} it follows that 
\[[v_1]_H=[x_1x_2\inv]_H,\ [v_2]_H=[x_2x_3\inv]_H,\ [v_3]_H=[x_3x_1\inv]_H.\]
Finally, from the definition of $E_i$,
\begin{align*}\co(\pi_{h_1}^\mQ(x_1x_2\inv))&=E_1\subseteq \co(\pi_{h_1}^\mQ(v_1))\\
\co(\pi_{h_2}^\mQ(x_2x_3\inv))&=E_2\subseteq \co(\pi_{h_3}^\mQ(v_2))\\
\co(\pi_{h_3}^\mQ(x_3x_1\inv))&=E_3\subseteq \co(\pi_{h_3}^\mQ(v_3))
\end{align*}
where $\mathrm{co(\underline{\phantom{x}})}$ 
denotes the word content of the sequence in question
(being a word over $\til E = Q \times \til A$). Also note that, for any words $y,z$ over $\til A$ and any $g\in Q$, $\pi_1^\mQ(y)\subseteq \pi_1^\mQ(z)$ if any only if $\pi_g^\mQ(y)\subseteq \pi_g^\mQ(z)$.

We construct the required tree $\mT=(W\cup  L;\alpha,\omega,\inv)$ as follows: we let $W=V\cup\{p\}=\{o_1,o_2,o_3,p\}$ with $p\ne o_i$ for  $i=1,2,3$ and edges $L=\{f_1^{\pm 1},f_2^{\pm 1},f_3^{\pm 1}\}$ where $f_i\colon o_i\longrightarrow p$ and labelling $w\colon  L\to \til A^*$ given by $w(f_i)=x_i$. Then, for all edges $k\in \mE$, if $\pi_k\colon  \alpha k\longrightarrow \omega k$ is the unique  reduced path in $\mT$ and $w_k=w(\pi_k)$, then $[v(k)]_H=[w_k]_H$ and $\co(\pi_1^\mQ(w_k))\subseteq\co(\pi_1^\mQ(v(k))$, as required. 
 Figure~\ref{fig:3cycle} indicates how the cycle (dotted arrows) is converted into a tree (solid arrows with additional vertex $\delta$), though this happens in the $G$-component of the corresponding $H$-values.

From Figure~\ref{fig:3cycle} one might get the impression that we have
tacitly assumed that $\delta\ne \eta_i$ and $[v_i]_H\ne 1_H$ for all
$i$. However it may happen that $\delta=\eta_i$ for some $i$; in this
case the path $\Pi_i$ and the word $x_i$ reduce to the empty path and 
empty word, respectively. The label $w(f_i)$ of the arrow $f_i$
will in turn be the empty word $1$. Since we allow the labelling of some edges
of $\mT$ to be the empty word $1$, this is not in conflict with our
construction. Similarly, if for some $i$,  $[v_i]_H=1_H$, then two of
the vertices $\eta_j$ in Figure~\ref{fig:3cycle} would coincide.
This would have the effect that two of the words $x_j$
coincide, which again would not be in conflict with our construction. 

\subsubsection{Cycle graphs of arbitrary length}
\label{subsubsection:gencycles}
\label{subsec: n-1 to n} We proceed by induction to prove the claim for all cycle graphs. Let $n>3$ and assume the claim to be true
for all cycle graphs of length at most $n-1$. In the following we use the cyclic index set $\mathbb{Z}_n$. Let $\mC=(V\cup  K;\alpha,\omega,\inv)$ be a cycle graph of length $n$ with vertex set $V=\{o_1,\dots,o_n\}$ and let $k_i\colon o_i\longrightarrow o_{i+1}$, $i\in \mathbb{Z}_n$, be edges.

Let an $H$-commuting labelling $v\colon K\to \til A^*$ be given by $v(k_i^{\pm 1})=v_i^{\pm 1}$ and let $[v_i]_H=(\gamma_i,g_i)$, that is, $g_i=[v_i]_Q$ and
$\gamma_i=[\pi_1^\mQ(v_i)]_G$. For $i=0,\dots,n-1$ we set
$h_{i+1}=g_1\cdots g_i$ and
$\eta_{i+1}={}^{h_1}\gamma_1\cdots{}^{h_i}\gamma_i$ (note that
$h_1=1_Q$ and $\eta_1=1_G$). Then $[v_1\cdots v_{i-1}]_H=(\eta_i,h_i)$
and in particular $\ol{\eta_i}=h_i$ (for $i=1,\dots,n$). Let $E_i$ be
the content of ${}^{h_i}\gamma_i=[\pi_{h_i}^\mQ(v_i)]_G$, that is, the $G$-content of
$\pi_{h_i}^\mQ(v_i)$ (see Definitions~\ref{def:word content} 
and~\ref{def:G-content}). Since $G$ is $n$-acyclic, the subgraph
\[\eta_1\mG[E_1]\cup \eta_2\mG[E_2]\cup\cdots\cup \eta_n\mG[E_n]\]
of $\mG$ is not a coset cycle. Since $\eta_{i+1}=\eta_i{}^{h_i}\gamma_i\in \eta_i\mG[E_i]$ and $1_G=\eta_1=\eta_n{}^{h_n}\gamma_n\in \eta_n\mG[E_n]$, condition (2) in Definition~\ref{def: cosetcycle} of coset cycles must be violated: 
for some~$i$ 
\[\eta_i\mG[E_{i-1}\cap E_i]\cap \eta_{i+1}\mG[E_i\cap E_{i+1}]\ne\varnothing.\]

\begin{figure}[ht]
\begin{tikzpicture}[xscale=0.9,yscale=0.9]
	\filldraw (-3,-3) circle (2pt);
	\filldraw(0,0) circle(2pt);
	\filldraw(-2.2,0)circle(2pt);
	\filldraw(4.5,-2)circle(2pt);
	\filldraw(1.8,1.8)circle(2pt);
	\draw(-1.1,-3.5)node{$\eta_{i-1}\mathcal{G}[E_{i-1}]=\eta_{i}\mathcal{G}[E_{i-1}]$};
	\draw(-.85,2.6)node{$\eta_{i}\mathcal{G}[E_{i}]=\eta_{i+1}\mathcal{G}[E_{i}]$};
	\draw(1.6,-2.4)node{$\eta_{i+1}\mathcal{G}[E_{i+1}]=\eta_{i+2}\mathcal{G}[E_{i+1}]$};
	\draw(-2.15,0.25)node[left]{$\eta_{i}$};
	\draw(0,0)node[below]{$\delta$};
	\draw(1.8,1.9)node[right]{$\eta_{i+1}$};
	\draw(4.5,-2)node[below]{$\eta_{i+2}$};
	\draw(-3,-3)node[right]{$\eta_{i-1}$};
	\draw plot [smooth cycle] coordinates {(-3.2,-4)(1,-4)(1,1)(-3.2,1)};
	\draw plot [smooth cycle] coordinates {(-3,-1)(2,-1)(2,2.7)(-3,2.7)};
	\draw plot [smooth cycle] coordinates {(-1,-2.6)(4.65,-2.6)(4.65,2)(-1,2)};
	\filldraw[ gray, nearly transparent] plotcoordinates {(-3.5,-4.7)(-4.5,-2)(-4,0.5)(-1.5,1.7)(-0.8,3.9)(0.6,3.9)(1.3,3.5)(2.1,3)(5,1)(6.2,0)(6,-3.5)(6,-4.7)(5.3,-4.7)(1.8,0)(0.5,-0.5)(-0.8,-0.5)(-1.8,-1)(-2.5,-4.7)};
	\draw(5.2,-3.5)node{$\widehat{\mathcal Q}$};
	\draw[thick,-latex](-2.2,0) to [out=20, in =160](-0.1,0.05);
	\draw[thick,-latex]((-3,-3) to[out=75, in=190] (-0.05,-0.05);
	\draw[thick, latex-] (0.05,0.05) to[out=20, in =120] (4.4,-1.95);
	\draw[thick,-latex] (1.75,1.75) to[out=190,in=80](0.0,0.1);
	\draw[dotted,thick,rounded corners,-latex] (-3.0,-3)--(-4,-2.5)--(-4,-0.5)--(-2.25,0);
	\draw[dotted,thick,rounded corners,-latex](-2.15,0.1)--(-1.3,0.9)--(-0.5,3.8)--(0.7,3.8)-- (1.75,1.85);
	\draw[dotted,thick,rounded corners,-latex](1.75,1.85)--(6,-0.3)--(6,-1)--(4.55,-2);
	\draw(-2.1,-0.75) node {$\Pi_1$};
	\draw(1.8,0.45) node{$\Pi_2$};
	\draw(-1.1,0.4) node{$\Lambda_1$};
	\draw(0.45,1.5) node{$\Lambda_2$};
\end{tikzpicture}
\caption{Transforming an $n$-cycle into a tree, inductively}\label{fig:n-1ton}
\end{figure}

All elements $1_G=\eta_1,\eta_2,\dots,\eta_n$ belong to $\wh{\mQ}$. From Lemma~\ref{lem:freeness} we get
\[\eta_i\mG[E_{i-1}\cap E_i]\cap\eta_{i+1}\mG[E_i\cap
E_{i+1}]\cap\wh{\mQ}\ne\varnothing.\]
Let $\delta$ be a vertex in this intersection; then, in particular, 
\[\delta\in \eta_{i-1}\mG[E_{i-1}]\cap \eta_{i+1}\mG[E_{i+1}]\cap\wh{\mQ}.\]
Hence there is an $E_{i-1}$-path $\Pi_1\colon
\eta_{i-1}\longrightarrow \delta$ and an $E_{i+1}$-path $\Pi_2\colon
\eta_{i+2}\longrightarrow \delta$, both running in
$\wh{\mQ}$. Likewise there is an $(E_{i-1}\cap E_i)$-path
$\Lambda_1\colon\eta_i\longrightarrow \delta$ and an $(E_i\cap
E_{i+1})$-path $\Lambda_2\colon \eta_{i+1}\longrightarrow \delta$,
both running in $\wh{\mQ}$.
The situation is depicted in Figure~\ref{fig:n-1ton}, where the dotted arrows indicate the paths with labels $\pi_{h_j}^\mQ(v_j)$ for $j=i-1,i,i+1$, the shaded area the $\mQ$-cover $\wh{\mQ}$. 
Note that $\eta_j\mG[E_{j}]= \eta_{j+1}\mG[E_{j}]$
     since $\eta_{j+1} \in \eta_{j} \mG[E_{j}]$ for $j=i-1,i,i+1$.

Since these four paths run in $\wh{\mQ}$ they are induced by words in $\til A^*$. So, there are $x_1,x_2,y_1,y_2\in \til A^*$ such that
\[\ell(\Pi_1)=\pi_{h_{i-1}}^{\mQ}(x_1),\ell(\Pi_2)=\pi_{h_{i+2}}^{\mQ}(x_2),\ell(\Lambda_1)=\pi_{h_{i}}^{\mQ}(y_1),\ell(\Lambda_2)=\pi_{h_{i+1}}^{\mQ}(y_2).\]
This leads to the following equalities
\[
\begin{aligned}
{[\pi_{h_{i-1}}^\mQ(v_{i-1})]_G}&=[\ell(\Pi_1\Lambda_1\inv)]_G=[\pi_{h_{i-1}}^\mQ(x_1y_1\inv)]_G\\
[\pi_{h_{i}}^\mQ(v_{i})]_G&=[\ell(\Lambda_1\Lambda_2\inv)]_G=[\pi_{h_{i}}^\mQ(y_1y_2\inv)]_G\\
[\pi_{h_{i+1}}^\mQ(v_{i+1})]_G&=[\ell(\Lambda_2\Pi_2\inv)]_G=[\pi_{h_{i+1}}^\mQ(y_2x_2\inv)]_G\\
\end{aligned}
\]	
so that, by Lemma~\ref{lem:wordproblemH},
\begin{equation}\label{eq:Hequivalence}		
[v_{i-1}]_H=[x_1y_1\inv]_H,\ [v_i]_H=[y_1y_2\inv]_H 	\mbox{ and }[v_{i+1}]_H=[y_2x_2\inv]_H,
\end{equation}
which implies 
\[[v_{i-1}v_iv_{i+1}]_H= [x_1x_2\inv]_H\]
and therefore also 
\begin{equation}\label{eq:C0compatible}
[v_1\cdots v_{i-2}x_1x_2\inv v_{i+1}\cdots v_n]_H=1_H.
\end{equation}

Furthermore, by construction, the paths $\Pi_1\Lambda_1\inv$, $\Lambda_1\Lambda_2\inv$ and $\Lambda_2\Pi_2\inv$ run in $\eta_{i-1}\mG[E_{i-1}]$, $\eta_i\mG[E_i]$ and $\eta_{i+1}\mG[E_{i+1}]$, respectively. For the respective contents of their labels this implies
\begin{align}
\co(\pi_{h_{i-1}}^\mQ(v_{i-1}))&\supseteq E_{i-1}=\co(\pi_{h_{i-1}}^\mQ(x_1y_1\inv)),\label{eq:xy}\\
\co(\pi_{h_i}^\mQ(v_i))&\supseteq E_{i\phantom{+1}}=\co(\pi_{h_i}^\mQ(y_1y_2\inv)),\label{eq:yy}\\
\co(\pi_{h_{i+1}}^\mQ(v_{i+1}))&\supseteq E_{i+1}= \co(\pi_{h_{i+1}}^\mQ(y_2x_2\inv))\label{eq:yx}
\end{align} 
since the $G$-content of a word is always contained in its word content.

In order to reduce the issue to the inductive hypothesis, we consider 
a cycle graph $\mC_0$ of length $n-1$ with sets of vertices
\[V_0=\{o_1,\dots,o_{i-1},p,o_{i+2},\dots,o_n\}\] 
and of edges 
\[\{k_t\colon t=1,\dots,i-2,i+2,\dots,n\}\cup\{f_1,f_2\}\]
(and their respective inverses),
where $f_1\colon o_{i-1}\longrightarrow
p$ and $f_2\colon p\longrightarrow o_{i+2}$ and the edges $k_t$ are as
in $\mC$, namely $k_t\colon o_t\longrightarrow o_{t+1}$
for all $t\ne i-1,i,i+1$. We define a labelling $v_0$ of $\mC_0$ by setting
$v_0(k_t)=v_t$ ($v_t$ as for $\mC$, for all $t\ne i-1,i,i+1$) and
$v_0(f_1)=x_1$ and $v_0(f_2)=x_2\inv$. In other words,  in order to
get $\mC_0$ from $\mC$, the labelled arc spanned by
\[o_{i-1}\overset{v_{i-1}}{\longrightarrow}o_i\overset{v_i}{\longrightarrow}o_{i+1}\overset{v_{i+1}}{\longrightarrow}o_{i+2}\]
in $\mC$ is replaced by the labelled arc spanned by
\[o_{i-1}\overset{x_1}{\longrightarrow}p\overset{x_2\inv}{\longrightarrow}o_{i+2}.\]

According to~\eqref{eq:C0compatible}, the labelling $v_0$ 
commutes over $H$. Hence by the inductive
assumption there exists a finite tree $\mT_0=(W_0\cup
{L_0},\alpha,\omega,\inv)$ with $V_0\subseteq W_0$ such that all
vertices in $W_0\setminus V_0$ have degree at least $3$.
We may assume that $o_i,o_{i+1}\notin W_0$. By the inductive
hypothesis, there is a labelling $w_0\colon {L_0}\to \til A^*$ such
that, for every edge $k$ of $\mC_0$, if $\pi_k\colon \alpha
k\longrightarrow \omega k$ is the unique  reduced path in $\mT_0$,
then $[v_0(k)]_H=[w_0(\pi_k)]_H$ and $\co(\pi_1^\mQ(w_0(\pi_k))\subseteq \co(\pi_1^\mQ(v_0(k))$.

\begin{figure}[ht]
\begin{tikzpicture}[scale=2]
	\filldraw(0,0) circle (1.5pt);
	\filldraw(-0.5,0.75) circle(1.5pt);
	\filldraw(0.5,0.75)circle(1.5pt);
	\draw[thick](0,0)--(-0.5,-0.3)--(-1,-0.9)--(-0.75,-1.3);
	\draw[thick](-0.5,-0.3)--(-0.2,-0.8)--(-0.6,-1.2);
	\draw[thick](-0.2,-0.8)--(0.1,-1.2);
	\draw[thick](0,0)--(0.5,-0.4)--(1,-0.9);
	\draw[thick](0.5,-0.4)--(0,-0.9);
	\draw[thick](0.5,-0.4)--(0.3,-1);
	\draw(0.7,-1.2) node{$\dots$};
	\draw[below](0,-0.5)node{$\mathcal{T}_0$};
	\draw[above](-0.475,0.8)node{$o_i$};
	\draw[above](.525,0.8)node{$o_{i+1}$};
	\draw[above](0,0.05)node{$p$};
	\draw[thick,-latex](-0.5,0.75)--(-0.05,0.05);
	\draw[thick,-latex](0.5,0.755)--(0.05,0.05);
	\draw(-0.35,0.35)node{$l_1$};
	\draw(0.375,0.35)node{$l_2$};
\end{tikzpicture}
\caption{The tree $\mT$}\label{fig:n-1ton_2}	
\end{figure}
We  define the tree $\mT=(W\cup {L},\alpha,\omega,\inv)$ according to
the picture in Figure~\ref{fig:n-1ton_2}:  let $W=W_0\cup
\{o_i,o_{i+1}\}$ and $L=L_0\cup\{l_1^{\pm 1},l_2^{\pm 1}\}$ where
$l_1\colon o_i\longrightarrow p$ and $l_2\colon o_{i+1}\longrightarrow
p$. We have $V\subseteq W$; the vertex $p$ is not in $V$
and has degree at least $3$ in $\mT$ by construction.
All other vertices not in $V$ are in $W_0\setminus V_0$ and 
have degree at least $3$ by the inductive assumption.
We define the labelling $w\colon L\to \til A^*$ by $w(l_t^{\pm 1})=y_t^{\pm 1}$ ($t=1,2$) and $w(l)=w_0(l)$ for every edge $l\in {L_0}$. Let $\pi_{f_1}\colon o_{i-1}\longrightarrow p$  and $\pi_{f_2}\colon p\longrightarrow o_{i+2}$ be the appropriate  reduced paths in $\mT_0$ corresponding to the edges $f_t$ in $\mC_0$ and denote their labels by $w_0(\pi_{f_t})=:w_t$ ($t=1,2$). Then $\pi_{f_1}l_1\inv\colon o_{i-1}\longrightarrow o_i$, $l_1l_2\inv\colon o_i\longrightarrow o_{i+1}$ and $l_2\pi_{f_2}\colon o_{i+1}\longrightarrow o_{i+2}$ are the  reduced paths in $\mT$ corresponding to the edges $k_{i-1}, k_i, k_{i+1}$ in $\mC$.  The labels of these paths are $w_1y_1\inv$, $y_1y_2\inv$ and $y_2w_2$, respectively.
By the inductive assumption on $\mC_0$ and $\mT_0$ we have $[w_1]_H=[x_1]_H$, $[w_2]_H=[x_2\inv]_H$ and 
$\co(\pi_1^\mQ(w_1))\subseteq \co(\pi_1^\mQ(x_1))$, $\co(\pi_1^\mQ(w_2))\subseteq \co(\pi_1^\mQ(x_2\inv))$, which implies 
\[[w_1y_1\inv]_H=[x_1y_1\inv]_H\mbox{ and }[y_2w_2]_H=[y_2x_2\inv]_H
\]
as well as
\[\co(\pi_1^\mQ(w_1y_1\inv))\subseteq\co(\pi_1^\mQ(x_1y_1\inv))\mbox{ and }
\co(\pi_1^\mQ(y_2w_2))\subseteq\co(\pi_1^\mQ(y_2x_2\inv)).\]

In combination with~\eqref{eq:Hequivalence} and~\eqref{eq:xy}, \eqref{eq:yy}, \eqref{eq:yx} we get for each one of the edges $k = k_{i-1},k_i,k_{i+1}$ that
\[[w_{k}]_H=[v(k)]_H\mbox{ and }
\co(\pi_1^\mQ(w_k))\subseteq \co(\pi_1^\mQ(v(k))\]
where $\pi_{k}$ is the  reduced path $\alpha k\longrightarrow \omega
k$ in $\mT$ and $w_k=w(\pi_k)$. For every edge $k$ different from $k_{i-1}^{\pm 1},k_i^{\pm
1}, k_{i+1}^{\pm 1}$ the claim follows from the inductive
hypothesis on $\mC_0$ and $\mT_0$,
since the corresponding path $\pi_k$ runs entirely in $\mT_0$.

Figure~\ref{fig:n-1ton} suggests that we have tacitly assumed that
$\eta_i\ne\delta\ne\eta_{i+1}$. 
However, as for the case of 
cycles of length~$3$, the discussion above
also covers the special cases $\eta_i=\delta$ or $\delta=\eta_{i+1}$: one of the paths $\Lambda_i$ and the corresponding word $y_i$ reduces to the trivial path and empty word, respectively. This is not  in conflict with the construction since we allow the tree $\mT$ to have empty labels.

\subsubsection{Arbitrary graphs} 
\label{subsubsection:gengraphs}
The general case will be handled by
induction on the number  $l$ of (geometric) edges of the  connected
graph $\mE$ on $n$ vertices. Base of the induction is $l=n-1$, in
which case $\mE$ is a tree. We may take $\mT=\mE$ and  $w=v$.
The idea is the same as for the proof of Lemma~\ref{lem: key lemma
RZ->Ash} and is borrowed from Section~7 of~\cite{Ash}. The proof is
actually simpler: we fix the group $H$ in advance, it works for all
connected graphs on at most $n$ vertices and need not be expanded in the
inductive process when we add edges.  We assume that $\mE$ has no loop
edges.

So, let $l>n-1$ and assume that the claim of Lemma~\ref{lem:main
lemma} is true for all connected graphs $\mD$ (without loop edges) on
$n$ vertices with fewer than $l$ geometric edges; let $\mE=(V\cup K,\alpha,\omega,\inv)$
be connected (without loop edges) with $|K^\bullet|=l$, $|V|=n$ and
an $H$-commuting labelling $v\colon  K\to \til A^*$ be given.

Recall that the $E$-generated group $G$ involved in the definition of
$H$ is $n$-acyclic. Let $k^{\pm 1}$ be a geometric edge of $\mE$
that is contained in a cycle subgraph of $\mE$ and let
$\mD=\mE\setminus\{k,k\inv\}$. The labelling $v\restriction \mathrm{K}(\mD)$
also commutes over $H$, and so the inductive assumption applies. So there exist a tree $\mX$  with $\mathrm{V}(\mD)\subseteq\mathrm{V}(\mX)$ 
and a labelling $w_\mX\colon \mathrm{K}(\mX)\to \til A^*$
such that:
\begin{enumerate}
\item[--]
all vertices of $\mX$ not in $\mD$ have degree at least $3$, 
\item[--]
for every edge $e$ of $\mD$, if $\pi_e$ is the  reduced path $\alpha e\longrightarrow \omega e$ in $\mX$, then 
\\
$[v(e)]_H=[w_\mX(\pi_e)]_H$ and $\co(\pi_1^\mQ(w_\mX(\pi_e)))\subseteq \co(\pi_1^\mQ(v(e)))$. 
\end{enumerate}
Note that the condition on the degree of the vertices guarantees that the diameter of $\mX$ is at most $n-1$.

Choose edges $k_1,\dots, k_t$ in $\mD$ such that $k_1\cdots k_t\colon \omega k\longrightarrow \alpha k$ is a simple path (that is, there are no repetitions of edges), then $kk_1\cdots k_t$ is a simple closed path in $\mE$. Since $v$ commutes over $H$, $[v(k)v(k_1\cdots k_t)]_H=1_H$.
Denote by $\pi\colon \omega k\longrightarrow \alpha k$ and
$\pi_j\colon \alpha k_j\longrightarrow \omega k_j$ ($j=1,\dots, t$)
the respective  reduced paths  in $\mX$. Then $\pi$ is the reduced
path of $\pi_1\cdots \pi_t$, whence
$[w_\mX(\pi)]_F=[w_\mX(\pi_1\cdots \pi_t)]_F$,
and therefore
$[w_\mX(\pi)]_H=[w_\mX(\pi_1\cdots \pi_t)]_H$.
Since $[w_\mX(\pi_j)]_H=[v(k_j)]_H$,
we get $[v(k)w_\mX(\pi)]_H=[v(k)v(k_1\cdots k_t)]_H=1_H$.

We look more carefully at the path $\pi$. Assume that $\pi=f_1\cdots
f_s$, where the $f_j$ are edges in $\mX$:
\[\omega k= p_0\underset{f_1}{\longrightarrow}p_1\underset{f_2}{\longrightarrow}\cdots\underset{f_{s-1}}{\longrightarrow} p_{s-1}\underset{f_s}{\longrightarrow}p_s=\alpha k;\] 
since the diameter of $\mX$ is at most $n-1$,
the length of $\pi$ is also at most $n-1$.
We add to $\mX$ a new edge $f_0\colon \alpha k=\omega \pi\longrightarrow \alpha \pi=\omega k$ and its inverse $f_0\inv$ and consider the closed path  $f_0\pi=f_0f_1\cdots f_s$ of length $s+1\le n$.  
Let $\mC$ be the cycle graph spanned by this closed path. We define a labelling $w_\mC$ of $\mC$ by $w_\mC(f_0):=v(k)$ and $w_\mC(f_j)=w_\mX(f_j)$ for $j\ne 0$ so that $w_\mC$ provides a labelling of $\mC$ which commutes over $H$. According to Section~\ref{subsec: n-1 to n}, there exists a tree $\mY$ on a vertex set $\mathrm{V}(\mY)\supseteq\mathrm{V}(\mC)=\{p_0,\dots,p_s\}$ containing the vertices of $\mC$, and a labelling $w_{\mY}$ of $\mY$ such that, for each of the edges $f_j$, if $\rho_j\colon \alpha f_j\longrightarrow \omega f_j$ denotes the  reduced path in $\mY$, then 
\begin{equation}\label{eq:assumption on D}
[w_{\mY}(\rho_j)]_H= [w_\mC(f_j)]_H
\mbox{ and }\co(\pi_1^\mQ(w_{\mY}(\rho_j)))\subseteq \co(\pi_1^\mQ(w_\mC(f_j))).
\end{equation}
We note that $w_\mC(f_j)=w_\mX(f_j)$ if $j\ne 0$ and $w_\mC(f_0)=v(k)$.
We may assume that $\mY$ and $\mX$ have no edges in common
and share no vertices other than those of $\mC$, viz.\ $p_0,\dots,p_s$.
Now define a tree $\mZ$ as
\[\mZ=\mY\cup \mX\setminus\{f_1^{\pm 1},\dots, f_s^{\pm1}\}\]
and a labelling $w_\mZ$ of $\mZ$ by 
\[w_\mZ(f)=\begin{cases} w_\mX(f)&\mbox{ if }f\in \mX \\ 
w_{\mY}(f)&\mbox{ if }f\in \mY. \end{cases}\]

The situation is completely the same as in the proof of
Lemma~\ref{lem: key lemma RZ->Ash}, cf.~Figure~\ref{fig:trees}.
We get a picture of how the trees $\mY$ and $\mX_0,\dots,
\mX_s$ are merged to form $\mZ$; the subtrees $\mX_{0}, \mX_{1},\dots,
\mX_{s}$ indicate the connected components of the forest $\mX\setminus
\{f_1^{\pm 1},\dots,f_s^{\pm 1}\}$,
$\mX_i$ being the component of $p_i$ in this forest. We note that $\mathrm{V}(\mE)=\mathrm{V}(\mD)\subseteq\mathrm{V}(\mX)\subseteq\mathrm{V}(\mZ)$.
Let $e$ be an edge of $\mE$ and let $\varsigma_e\colon \alpha e\longrightarrow \omega e$ be the  reduced path in $\mZ$. If $e=k$ then $\varsigma_e=\varsigma_k$ runs completely in $\mY$ and we have \[[v(k)]_H=[w_\mC(f_0)]_H=[w_{\mY}(\rho_0)]_H=[w_\mZ(\varsigma_k)]_H\] and 
\[\co(\pi_1^\mQ(v(k)))=\co(\pi_1^\mQ(w_\mC(f_0)))\supseteq \co(\pi_1^\mQ(w_{\mY}(\rho_0)))=\co(\pi_1^\mQ(w_\mZ(\varsigma_k)))\]
by the construction of $\mY$ from $\mC$, see \eqref{eq:assumption on D}. Suppose that $e\ne k^{\pm 1}$, that is $e\in \mD$. Then 
\begin{equation}\label{eq:equalitiforpi'}
[v(e)]_H=[w_\mX(\pi_e)]_H\mbox{ and } 
\co(\pi_1^\mQ(v(e)))\supseteq \co(\pi_1^\mQ(w_\mX(\pi_e)))
\end{equation}
by the inductive hypothesis on $\mX$, where $\pi_e$ is the  reduced path $\alpha e\longrightarrow \omega e$ in $\mX$.
The path $\pi_e$ may use some segment of the arc spanned by
$f_1\cdots f_s$, so 
\begin{equation}\label{eq:definitionofpii}
\pi_e=\pi_1f_l\cdots f_m\pi_2\mbox{ or }\pi_e=\pi_1f_m\inv \cdots f_l\inv\pi_2
\end{equation}
for some  reduced paths $\pi_1,\pi_2$ in $\mX\setminus\{f_1^{\pm 1},\dots,f_s^{\pm 1}\}$ and $l\le m$. We handle the first case, the second one is completely analogous. Let $\rho\colon \alpha f_l\longrightarrow \omega f_m$ be the  reduced path in $\mY$. Then the  reduced path $\varsigma_e\colon \alpha e\longrightarrow \omega e$ in $\mZ$ is $\varsigma_e=\pi_1\rho\pi_2$ and from the definition of $w_\mZ$:
\begin{equation}\label{eq:i to prime}
w_\mZ(\varsigma_e)=w_\mX(\pi_1)w_{\mY}(\rho)w_\mX(\pi_2).
\end{equation}
The path $\rho$ is the reduced form of $\rho_l\cdots \rho_m$; on the one hand, this implies  $[w_{\mY}(\rho)]_F=[w_{\mY}(\rho_l\cdots\rho_m)]_F$ so that $[w_{\mY}(\rho)]_H=[w_{\mY}(\rho_l\cdots\rho_m)]_H$.  On the other hand we get $\co(\pi_1^\mQ(w_\mY(\rho)))\subseteq \co(\pi_1^\mQ(w_{\mY}(\rho_l\cdots\rho_m)))$.
From~\eqref{eq:assumption on D} we have 
\[[w_{\mY}(\rho_j)]_H=[w_\mX(f_j)]_H\mbox{ and }
\co(\pi_1^\mQ(w_{\mY}(\rho_j)))\subseteq\co(\pi_1^\mQ(w_\mX(f_j)))\]
for every $j\ne 0$. So from~\eqref{eq:i to prime} and the structure
of $\pi_e$ in~\eqref{eq:definitionofpii} we find
$[w_\mX(\pi_e)]_H=[w_\mZ(\varsigma_e)]_H$ and
$\co(\pi_1^\mQ(w_\mX(\pi_e)))\supseteq \co(\pi_1^\mQ(w_\mZ(\varsigma_e)))$. In combination
with~\eqref{eq:equalitiforpi'} we get for $e\ne k^{\pm 1}$ (what we had above already for $e=k^{\pm 1}$):
\begin{equation}\label{eq:main assertin of w'}
[v(e)]_H=[w_\mZ(\varsigma_e)]_H\mbox{ and }
\co(\pi_1^\mQ(v(e)))\supseteq \co(\pi_1^\mQ(w_\mZ(\varsigma_e))).	
\end{equation}

We are \textsl{almost} done, but not yet completely: it may
happen that the tree $\mZ$  contains vertices not in $\mathrm{V}(\mE)$
of degree smaller than $3$. This may occur as follows: some
intermediary vertices in the path $f_1\cdots f_s$ may be not in
$\mathrm{V}(\mE)=\mathrm{V}(\mD)$ and have degree $3$ in $\mX$; 
these vertices may have degree $1$ in $\mY$ (in fact, they do have degree
$1$ as can be checked in the construction
presented in Section~\ref{subsec: n-1 to n}). Every such vertex has degree $2$
in $\mZ$. To rectify this, consider all maximal arcs spanned by paths
\begin{equation}\label{eq:bad path}
q_0\underset{e_1}{\longrightarrow}q_1\underset{e_2}{\longrightarrow}q_2\cdots q_{t-1}
\underset{e_t}{\longrightarrow}q_t 
\end{equation} of consecutive edges $e_1, \dots, e_t$ in $\mZ$ such
that the vertices $q_1,\dots, q_{t-1}$ have degree $2$ and are not in
$\mathrm{V}(\mE)$  while each of $q_0$ and $q_t$ is either in
$\mathrm{V}(\mE)$ or has degree at least $3$. 
In order to get the required tree $\mT$ from $\mZ$,
we may proceed as we did earlier (proofs of Lemma~\ref{lem:1stpart} and 
Lemma~\ref{lem: key lemma RZ->Ash}): just contract the path~\eqref{eq:bad path} to a single edge, that is, replace every  arc spanned by a path of the form~\eqref{eq:bad path} by the single edge
$e\colon q_0\longrightarrow q_t$ and label $e$ in  $\mT$
by $w_\mT(e):=w_\mZ(e_1)w_\mZ(e_2)\cdots w_\mZ(e_t)$ (and provide inverse edges accordingly);  edges $e$ of $\mZ$
that are not on such a path remain unchanged in $\mT$ and we set
$w_\mT(e):=w_\mZ(e)$. In the new tree $\mT$ just vertices from
$\mathrm{V}(\mE)$ may have degree less than $3$.
Every reduced path in $\mZ$
that traverses one of the edges $e_i$ in~\eqref{eq:bad
path} needs to traverse the entire arc~\eqref{eq:bad path}. 
So, if $\varsigma_e\colon \alpha
e\longrightarrow \omega e$ denotes the  reduced path in $\mZ$ and
$\tau_e\colon \alpha e\longrightarrow \omega e$ is the one in $\mT$
for an edge edge $e\in \mE$, then 
$w_\mZ(\varsigma_e)=w_\mT(\tau_e)$ (i.e.\ the labellings of these paths coincide).
Finally, from~\eqref{eq:main assertin of w'} we get
\[[v(e)]_H=[w_\mT(\tau_e)]_H\mbox{ and } 
\co(\pi_1^\mQ(v(e)))\supseteq \co(\pi_1^\mQ(w_\mT(\tau_e)))\]
for every edge $e\in \mE$, where $\tau_e\colon \alpha e\longrightarrow \omega e$  denotes the  reduced path in the new tree $\mT$.
\end{proof}

While no inverse monoid occurs explicitly in the preconditions and the formulation of Lemma~\ref{lem:main lemma} there is one ``lurking from behind''. The precise meaning of this will become clear from Lemma~\ref{lem:meakmarg} below, but intuitively it can be seen as follows. Let $Q$ be an $A$-generated group (not necessarily finite) and consider the equivalence relation $\sim_Q$ on the free involutory monoid $\til A^*$ defined by
\begin{equation}\label{eq:meakmarg}
 u\mathrel{\sim_Q} v\Longleftrightarrow [u]_Q=[v]_Q\mbox{ and }\co(\pi_1^\mQ(u))=\co(\pi_1^\mQ(v)).    
\end{equation}
This means that two words $u$ and $v$ relate under $\sim_Q$ if and only if
\begin{enumerate}
    \item $u$ and $v$ evaluate equally in $Q$: $[u]_Q=[v]_Q=:g$ and
    \item in the Cayley graph $\mQ$ of $Q$, the paths $1_Q\longrightarrow g$ labelled $u$ and $v$, respectively,  traverse the same geometric edges. 
\end{enumerate}
It turns out that $\sim_Q$ is a congruence on the free involutory monoid $\til A^*$ and the quotient  $M(Q):=\til A^*/{\sim_Q}$ is an inverse monoid. The $\sim_Q$-class of a word $u\in \til A^*$ is determined by the pair $(\co(\pi_1^\mQ(u)),[u]_Q)=:(R,g)$ where $R=\co(\pi_1^\mQ(u))\subseteq E=Q\times A$ and $g=[u]_Q\in Q$. Moreover, given two words, $u$ and $v$, say, then, for $g=[u]_Q$, 
\[\co(\pi_1^\mQ(uv))=\co(\pi_1^\mQ(u))\cup \co(\pi_g^\mQ(v))=\co(\pi_1^\mQ(u))\cup g\co(\pi_1^\mQ(v)).\]  It follows that the pairs $(R,g)\in 2^E\times Q$ describing congruence classes modulo $\sim_Q$ multiply by the rule $(R,g)(S,h)=(R\cup gS,gh)$. In other words, the inverse monoid $M(Q)$ can be realised as a submonoid of the semidirect product $2^E\rtimes Q$ where $2^E$ is the semilattice monoid (i.e.~the commutative and idempotent monoid) of all subsets of $E$ with $\cup$ as binary operation and the action of $Q$ on $2^E$ is the one induced by left multiplication of the elements of $E=Q\times A$ by elements of $Q$. One can check that in the inverse monoid $M(Q)$ the natural order is described by
\[(R,g)\le (S,h)\mbox{ if and only if } g=h\mbox{ and }R\supseteq S.\]
Looking back at Lemma~\ref{lem:main lemma} we see that (in terms of the notation of that lemma) we \textsl{actually} have proved that the labelling $w$ of the tree $\mT$ satisfies:
\begin{enumerate}
    \item $[v(k)]_H=[w_k]_H$ and
    \item $[v(k)]_{M(Q)}\le [w_k]_{M(Q)}$.
\end{enumerate}
Item (2) implies in particular that $[v(k)]_M\le [w_k]_M$ for every $A$-generated inverse monoid $M$ for which $M(Q)\twoheadrightarrow M$. The inverse monoid $M(Q)$ is called the \emph{Margolis--Meakin expansion} of $Q$~\cite{MM}. It plays an important r\^ole in the theory of inverse monoids and will be of significant use in Section~\ref{sec:strictly stronger}. Usually, the elements of $M(Q)$ are described by pairs $(\mR,g)$ where $\mR$ is a subgraph of $\mQ$ rather than pairs $(R,g)$ with $R\in 2^E$. The connection between the pairs $(\mR,g)$ and $(R,g)$ is clear, of course: instead of $R\subseteq E$ just take $\mR:=\langle R\rangle$, the subgraph of $\mQ$ spanned by $R$. This view allows one to characterise completely which pairs $(R,g)\in 2^E\times Q$ actually belong to $M(Q)$: $(R,g)$ belongs to $M(Q)$ if any only if $\langle R\rangle$ is a finite connected subgraph of $\mQ$ containing $1_Q$ and $g$ (with the natural convention that $\langle\varnothing\rangle=\{1_Q\}$).
The collection of all $A$-generated inverse monoids $M$ is 
partially ordered by $M\le N$ if and only if $N\twoheadrightarrow M$. The special r\^ole  of $M(Q)$ then lies in the fact that it is the maximum element of all $A$-generated inverse monoids $M$ that are extended by $Q$. Essentially, this follows from the next lemma. 

Let $M$ be an $A$-generated inverse monoid
and $Q$ an $A$-generated group extending $M$;
recall the word content $\mathrm{co}(w)$ of a word $w\in \til E^*$ according to Definition~\ref{def:word content}, for $E=Q\times A$. 

\begin{Lemma}\label{lem:meakmarg} For $u,v\in \til A^*$ such that $[u]_Q=[v]_Q$, 
the inclusion relationship $\langle\pi_1^\mQ(u)\rangle \subseteq  \langle \pi_1^\mQ(v)\rangle$, or equivalently, $\mathrm{co}(\pi_1^\mQ(u))\subseteq \mathrm{co}(\pi_1^\mQ(v))$, implies $[u]_M\ge[v]_M$.	
\end{Lemma}

\begin{proof} It is clear that  $\langle\pi_1^\mQ(u)\rangle \subseteq  \langle \pi_1^\mQ(v)\rangle$ if and only if $\mathrm{co}(\pi_1^\mQ(u))\subseteq \mathrm{co}(\pi_1^\mQ(v))$. 
Let $u,v\in \til A^*$ be such that $[u]_Q=[v]_Q$  and $\langle\pi_1^\mQ(u)\rangle \subseteq  \langle \pi_1^\mQ(v)\rangle$. 
We represent the $A$-generated inverse monoid $M$ as a
monoid of partial bijections on a finite set $X$; we
need to show that the partial bijection defined by $v$
is a restriction of the one defined by $u$.  Let $\mP$
be the $A$-graph on the vertex set $X$ encoding
the partial bijections $a\in \til A$. We consider
the \tilA{Q}-set $X_\mP Q$ generated by $X$ subject to $\mP$ defined in \eqref{eq:GAsubjectomP}. 
Our assumption on $Q$ implies that $Q$ extends $\mP$ in the sense of Definition~\ref{def: G extends mP}.
Take any connected component $\mP_i$ of $\mP$ and let $\mE_i$ be 
the component of $X_\mP Q$ containing $\mP_i$ as a subgraph. Then there exists,
for every vertex $y$ of $\mE_i$, a canonical graph morphism
$\varphi_y\colon\mQ\twoheadrightarrow \mE_i$ mapping $1_Q$ to $y$
($\mQ$ the Cayley graph of $Q$).\footnote{This follows from the fact that
$\mE_i$ is isomorphic to the Schreier graph $\bm{\mS}_{H}$ for some
subgroup $H$ of $F$ containing the normal subgroup $N$  of $F$ for which
$G=F/N$, and
clearly $\mQ\cong\bm{\mS}_N\twoheadrightarrow \bm{\mS}_{H}$.}
Let $x\in X$ be a vertex of $\mP_i$ and consider the graph morphism $\varphi_x$. Since $\langle \pi_1^\mQ(v)\rangle  \supseteq \langle \pi_1^\mQ(u)\rangle$ the same inclusion holds for the images of these graphs under $\varphi_x$. Hence, if the path in $\mE_i\subseteq X_\mP Q$ starting at $x$ and being labelled $v$ runs in $\mP_i$,  then so does the one labelled $u$.  Consequently, if ${x}\cdot[v]_M$ is defined then so is ${x}\cdot[u]_M$. If both are defined then they are equal since 
\[{x}\cdot[v]_M=x[v]_Q=x[u]_Q={x}\cdot[u]_M.\]
It follows that $[v]_M\le[u]_M$.
\end{proof}
The preceding lemma is essentially  folklore. 
Proofs of the lemma based on an analysis of the words $u$ and $v$ and not using the interpretation of $M$ as a
monoid of partial transformations can be found in
\cite{MM,steinberginverseautomata,bitterlichdiss,AKS, ABO}, the present proof is sketched in~\cite{ASconstructiveRZ} and perhaps elsewhere. 
From this lemma one can conclude that the canonical morphism $\til A^*\twoheadrightarrow M$ (for $M$ extended by $Q$) factors through $M(Q)$, and so $M(Q)$ is indeed the maximum among all $A$-generated inverse monoids extended by $Q$.
Therefore, Lemma~\ref{lem:meakmarg} implies that if $Q$ extends $M$ then statement (2)  $\co(\pi_1^\mQ(v(k))\supseteq \co(\pi_1^\mQ(w_k))$ in Lemma~\ref{lem:main lemma} may be replaced by $[v(k)]_M\le [w_k]_M$. 
We finally show that Lemma~\ref{lem:main lemma} implies
Theorem~\ref{thm:otto-ash}.

\begin{proof}[Proof of Theorem~\ref{thm:otto-ash} using
Lemma~\ref{lem:main lemma} and Lemma~\ref{lem:meakmarg}] Let $M$ and  $n$ be as in the statement of  Theorem~\ref{thm:otto-ash} and let $\mE$ be a connected graph on $n$ vertices.
Let $\mD$ be the graph obtained by removing all
loop edges of $\mE$. Let $Q$ be a finite $A$-generated group
extending $M$, $H$ the $A$-generated group defined by \eqref{eq:group H} with $G$ as in the statement of Lemma~\ref{lem:main lemma}.
By this lemma, there is a finite tree $\mT$ with
$\mathrm{V}(\mE)=\mathrm{V}(\mD)\subseteq \mathrm{V}(\mT)$ and a
labelling $w\colon \mathrm{K}(\mT)\to \til A^*$
as in  the statement of the lemma (for $\mD$
in the r\^ole
of $\mE$), but  ``$\co(\pi_1^\mQ(v(k))\supseteq\co(\pi_1^\mQ(w_k))$'' replaced by ``$[v(k)]_M\le [w_k]_M$''. 
Define the labelling $u\colon \mathrm{K}(\mE)\to \til A^*$ by $u(e):=v(e)v(e)\inv w_e$ for every edge $e$ of $\mE$ where $w_e=w(\pi_e)$ and $\pi_e$ is the  reduced path $\alpha e\longrightarrow \omega e$ in $\mT$ (which is empty if $\alpha e=\omega e$). 
For every loop edge $e$, $\pi_e$ is an empty path and we have $[u(e)]_H=[v(e)]_H=1_H$ since $v$ commutes over $H$, and $[u(e)]_M=[v(e)]_M$ since $[v(e)]_Q=1_Q$ implies that $[v(e)]_M$ is idempotent, hence $[v(e)]_M=[v(e)v(e)\inv]_M=[u(e)]_M$. 
For every non-loop edge the definition also implies $[u(e)]_H=[v(e)v(e)\inv w_e)]_H=[v(e)]_H$ since $[v(e)]_H=[w_e]_H$ and $[u(e)]_M=[v(e)v(e)\inv w_e]_M=[v(e)]_M$ since $[v(e)]_M\le[w_e]_M$.
Moreover,  $u$ commutes over $F$: for every closed path $e_1\cdots e_t$ in $\mE$ the reduced path of the corresponding product path $\pi_{e_1}\cdots\pi_{e_t}$ in $\mT$ is the empty path, hence the reduced form of the label of $e_1\cdots e_t$ is the empty word. 
\end{proof}

\subsection{Applications} We take Theorem~\ref{thm:otto-ash} and
see what we get if we apply it in the context of the different versions of the Ash/Herwig--Lascar/Ribes--Zalesskii Theorem.  {In the following, let $X$ be a finite set, $\tilde A$ a set of partial bijections of $X$, encoded in a finite $A$-graph $\mathcal{P}$; if $X$ is the base set of a $\sigma$-structure $\mathfrak{X}=(X;(R^\mathfrak{X})_{R\in\sigma})$, then the partial bijections are assumed to be partial automorphisms of $\mathfrak{X}$. For an $A$-generated group $G$ recall definition~\eqref{eq:GAsubjectomP} of the $\tilde{A}$-set $X_\mathcal{P}G$ and the $\sigma$-structure $\mathfrak{X}_\mathcal{P}G$; the free $A$-generated group is denoted by $F$.}

\subsubsection{Model-theoretic Herwig--Lascar}\label{sec:model-theoretic HL} Looking at the proof  in
Section~\ref{subsec: Ash===>HL} we see that  {Proposition~\ref{prop:abstract-extension} below has implicitly been proved there; this proposition seems to be interesting in its own right as it appears to be a statement that is equivalent to the theorems of Ash and Herwig--Lascar, and in some sense ``sitting between them''. The equivalence claimed follows from the fact that the proof in Section~\ref{subsec: Ash===>HL} essentially consists of showing the chain of implications
\begin{center}
\begin{tabular}[t]{c}Theorem~\ref{thm: ash} 
\\ Ash
\end{tabular}
$\;\Longrightarrow\;\;$ Proposition~\ref{prop:abstract-extension} $\;\;\Longrightarrow\;$ \begin{tabular}[t]{c} Theorem~\ref{thm: Herwig--Lascar}.
\\
Herwig--Lascar
\end{tabular}
\end{center}}

 {\begin{Prop}\label{prop:abstract-extension}
 For any finite set $X$ endowed with a set of partial bijections $\tilde{A}$ and every positive integer $n$, there exists a finite $A$-generated group $G$ such that for every choice of elements  $g_1,\dots, g_n\in G$ there exist $f_1,\dots,f_n\in F$ for which the set of pairs
\[\{(xg_i,xf_i)\colon x\in X, i=1,\dots, n\}\]
is a partial function $X_\mathcal{P}G\to X_\mathcal{P}F$.
\end{Prop}}
 {The condition formulated in Proposition~  \ref{prop:abstract-extension} just says that the assignment $xg_i\mapsto xf_i$ is a well-defined function $\bigcup_{i=1}^n Xg_i \to \bigcup_{i=1}^n Xf_i$; that is, for any $x,y\in X$ and any $i,j\in \{1,\dots,n\}$ the following implication holds:
\[xg_i=yg_j\Longrightarrow xf_i=yf_j.\]
An analysis of the proof in Section~\ref{subsec: Ash===>HL} then leads to the following reformulation of the Herwig--Lascar Theorem; recall the weight of a finite relational structure as defined in the first paragraph of Section~\ref{subsec: Ash===>HL}.}
 {\begin{Prop}[re-interpretation of Herwig--Lascar]
\label{prop:EPPAin4.3.1}
For every positive integer $n$, every finite extension problem $(\mathfrak{X},\tilde{A})$ has a finite solution of the form $\mathfrak{X}_\mathcal{P}G$ 
such that 
every weak substructure $\mathfrak{U} \subseteq \mathfrak{X}_\mathcal{P}G$ 
of weight at most $n$ admits a homomorphism $\phi\colon \mathfrak{U}\to \mathfrak{X}_\mathcal{P}F$.
\end{Prop}}

So, in a sense, the proof in Section~\ref{subsec: Ash===>HL} of the Herwig--Lascar Theorem provides, for every finite $\sigma$-structure $\bX$ 
with a set of partial automorphisms (encoded in an $A$-graph $\mP$) and any finite collection $\Sigma$ of finite $\sigma$-structures (for the specification of the forbidden homomorphisms), a finite extension $\bX_\mP G$, that --- from point of view of $\Sigma$ --- \textsl{looks like} the free extension $\bX_\mP F$. The group $G$, and hence the extension $\bX_\mP G$, just depends on the maximum of the weights of the structures in $\Sigma$, 
not on their isomorphism types.  {For the following, recall that a \emph{hypergraph} $\mathbb{V}=(V;\mathbf{E})$ consists of a set $V$ (of vertices) and a set $\mathbf{E}$ of finite subsets of $V$ (the \emph{hyperedges} of $\mathbb{V}$).} 
A similar result in~\cite{otto4} uses criteria of hypergraph acyclicity 
in application to the hypergraph structure induced by the $G$-orbit of
$X$ in $\bX_\mP G$. In the terminology of~\cite[Proposition~5.10]{otto4}
the hypergraph whose hyperedges are the images of the    
base set $X$ under the automorphisms induced by~$G$ is \emph{locally acyclic} 
in the sense that, up to a certain size bound, all its sub-hypergraphs 
are acyclic (i.e.\ chordal and conformal, or equivalently: tree-decomposable,
see~\cite{BeeriFaginetal,Berge} for the terminology).%
\footnote{And Corollary~5.12 there uses this local similarity 
to the free extension $\bX_\mP F$, which is globally tree-like in this sense, 
to rule out forbidden homomorphisms 
towards EPPA for classes $\mathbf{Excl}(\Sigma)$, similar to a use of Theorem~\ref{thm: HL strengthened} below.} 
The following strengthening of the Herwig--Lascar Theorem makes
it even clearer in which sense the resulting finite extension
$\bX_\mP G$ of the given relational structure $\bX$ provides an
approximation to the free extension $\bX_\mP F$.
In a sense, such
extensions are \textsl{locally free}, see Figure~\ref{fig: XGvsXF}
below. Other strengthenings of the Herwig--Lascar Theorem may be
found in~\cite{eppa classes} and in some
references cited therein. Recall the concepts used in
Section~\ref{subsec: Ash===>HL}: 
given a finite $\sigma$-structure $\bX=(X;(R^\bX)_{R\in\sigma})$ together
with a set $\til A$ of partial automorphisms of $\bX$ encoded in an
$A$-graph $\mP$ on the vertex set $X$,  $\bX_\mP F$ is the free
extension of $\bX$;  for any $A$-generated
group $G$ satisfying 
the conditions of Lemma~\ref{lem:X is induced}  we have the extension
$\bX_\mP G$. 
We denote the canonical morphism
$F\twoheadrightarrow G$ by $f\mapsto \ol{f}$.
Then the canonical mapping $\bX_\mP F\twoheadrightarrow
\bX_\mP G$, $xf\mapsto x\ol{f}$  is a morphism of $\til A$-sets
as well as of $\sigma$-structures.

 {\begin{Thm}[strengthened version of (re-interpreted) Herwig--Lascar]\label{thm: HL strengthened}
For every positive integer $n$, every finite extension problem $(\mathfrak{X};\tilde{A})$ has a finite solution of the form $\mathfrak{X}_\mathcal{P}G$ 
with the following property:
every weak substructure $\mathfrak{U} \subseteq \mathfrak{X}_\mathcal{P}G$ 
of weight at most $n$ admits a homomorphism $\phi\colon \mathfrak{U}\to \mathfrak{X}_\mathcal{P}F$ such that the canonical mapping $\mathfrak{X}_\mathcal{P}F\twoheadrightarrow\mathfrak{X}_\mathcal{P}G$ induces an isomorphism $\psi\colon \phi(\mathfrak{U})\to \mathfrak{U}$ for which $\psi\circ\phi=\mathrm{id}_\mathfrak{U}$.
\end{Thm} }
\begin{proof}
We proceed as in the proof of Theorem~\ref{thm: Herwig--Lascar}
presented in Section~\ref{subsec: Ash===>HL}.
Let $M=\cT(\mP)$ be the transition monoid of
$\mP$ and $G$ a finite $A$-generated group
according to Theorem~\ref{thm:otto-ash} for all connected graphs on at most $n$ vertices. Let $\bX_\mP G$  be the extension induced by $G$  {and let $\mathfrak{U}$ be a weak substructure of $\mathfrak{X}_\mathcal{P}G$ of weight at most $n$}.
There exists a subset $J\subseteq G$ of size at most $n$ such that
every
relational tuple $(t_1,\dots,t_r)$ of $\bT$ and every 
element of $\bT$ is contained in $Xg$ for some $g\in J$. 
Define a graph
structure $\mJ$ on the set $J$ exactly as in the proof
of Theorem~\ref{thm: Herwig--Lascar} given 
in Section~\ref{subsec: Ash===>HL}, 
with the same labelling as in that proof.  
That is, the set of vertices $\mathrm{V}(\mJ)$ is $J$ and  two vertices  $g$ and $h$ are connected by one or more edges if and only if $Xg\cap Xh\ne \varnothing$. More precisely, there will be $|Xg\cap Xh|$ many edges from $g$ to $h$ and also in the reverse direction. Suppose that $Xg\cap Xh=\{z_1,\dots, z_s\}$; then there are $x_1,\dots, x_s, y_1\dots,y_s\in X$ such that $x_ig=z_i=y_ih$ for all $i$, hence $x_igh\inv =y_i$ for all $i$. According to Lemma~\ref{lem:extending G} there exist words $v_1,\dots,v_s\in \til A^*$ such that $x_igh\inv =y_i={x_i}\cdot[v_i]_M$ 
and $[v_i]_G=gh\inv$. In this situation we put 
edges $e_1,\dots, e_s\colon g\longrightarrow h$ with labels $v_1,\dots, v_s$, respectively. The corresponding inverse edges $e_i\inv\colon h\longrightarrow g$ then have labels $v_i\inv$ and satisfy $y_ihg\inv = x_i = {y_i}\cdot[v_i\inv]_M$. 
We do this for all pairs $g,h\in J$ for which $Xg\cap Xh\ne\varnothing$.  We have thus defined a word labelling of $\mJ$. 
This  labelling is $G$-commuting since for every edge $e\colon g\longrightarrow h$ the $G$-value of its label is $gh\inv$. Hence the label of any closed path has $G$-value $1_G$. According to Theorem~\ref{thm:otto-ash} there exists a relabelling which commutes over the free group $F$.
Compared to the proof of Theorem~\ref{thm: Herwig--Lascar} given 
in Section~\ref{subsec: Ash===>HL}, the relabelling now satisfies an additional condition. Back there (that is, in the proof in Section~\ref{subsec: Ash===>HL})
we already had $[v_i]_F=[w_i]_F=[w^{(g,h)}]_F$ and $[v_i]_M=[w_i]_M\le[w^{(g,h)}]_M$ for the labels of
$e_1,\dots,e_s\colon g\longrightarrow h$ for specific $g,h\in G$ for
which $Xg\cap Xh\ne \varnothing$. Now we additionally have that
$[w^{(g,h)}]_G=[w_i]_G=[v_i]_G=gh\inv $.
We proceed as in the earlier proof. 
Take a connected component of $\mJ$ and choose a base vertex
$o$. Choose any $f_o\in F$ such that $\ol{f_o}=o$; for every $g$ in
that component of $\mJ$ let $f_g := [p]_Ff_o$ where $p$ is the label of some
path $g\longrightarrow o$ (the labelling is with respect to the latest relabelling which provided every edge $e\colon g\longrightarrow h$ with the label $w^{(g,h)}$ mentioned above) --- this does not depend on the choice of
that path because of $F$-commutativity of the labelling.
Taking into account that the $G$-value of the label of every path
$g\longrightarrow o$ is $go\inv$, 
it follows that
$\ol{f_g}=\overline{[p]_Ff_o}=[p]_G\overline{f_o}=go\inv o=g$. Hence
the choice of $f_g$ satisfies $\ol{f_g}=g$. This holds for all $g\in J$
since this can be done for every connected component of $\mJ$. As in the proof of Theorem~\ref{thm: Herwig--Lascar} in Section~\ref{subsec: Ash===>HL}, the mapping $ {\phi}\colon\bT\to \bX_\mP F$  defined by
\begin{equation}\label{eq:isoHL}
	\begin{array}{rcccl}   \bT\cap Xg & \longrightarrow &  X &\longrightarrow &Xf_g \subseteq X_\mP F\\
		z  & \longmapsto    & zg\inv& \longmapsto & (zg\inv)f_g \end{array}
\end{equation}
is well-defined (independent of the choice of $g$)
and a morphism of the $\sigma$-structure $\bT$. Since $\ol{f_g}=g$ for
all $g$,
it follows that the composition of
this morphism with the canonical mapping $\bX_\mP
F\twoheadrightarrow \bX_\mP G$
is the identity mapping on $\bT$. As a consequence, the
mapping~\eqref{eq:isoHL}   {$\bT\to\wh{\bT} :=\phi(\mathfrak{U})$} is bijective and
a homomorphism of $\sigma$-structures in both directions,
hence an isomorphism.
\end{proof}

\begin{figure}[ht]
\begin{tikzpicture}[scale=1.3]
	\draw plot [smooth cycle] coordinates  {(2,2.4)(4,2.6)(5,1.5)(3,-1)(1,-1.4)(-1,-1)};
	\draw plot [smooth cycle] coordinates   {(-1,1.2)(-2,1.3)(-2,-1)(-1,-1.5)(0.5,-0.5)};
	\filldraw[fill=gray, nearly transparent](-1.5,0)circle(0.5cm);
	\draw(-1.5,0)node{$\mathfrak{U}$};
	\draw(-0.3,-0.8)node{$\mathfrak{X}$};
	\draw[->>](3,2.5) .. controls (2,4) and (0,4) .. (-1.8,1);
	\draw(1.2,3.7)node{$xf\mapsto x\overline{f}$};
	\draw(-1.3,-1.8)node{$\mathfrak{X}_\mathcal{P}G$};
	\draw(3,-1.4)node{$\mathfrak{X}_\mathcal{P} F$};
	\filldraw[fill=gray, nearly transparent](2,1.5)circle(0.5cm);
	\draw(2,1.5)node{$\widehat{\mathfrak{U}}k$};
	\filldraw[fill=gray, nearly transparent](3,0)circle(0.5cm);
	\draw(3,0)node{$\widehat{\mathfrak{U}}$};
	\draw[>->>](1.9,1.8)--(-1.5,0.2);
	\draw[>->>](3,-0.3)--(-1.5,-0.2);
	\draw[>->>](2.85,0.25)--(2.1,1.25);
	\draw(2.65,0.9)node{$\cdot k$};
\end{tikzpicture}
\caption{Configuration as in the proof of Theorem~\ref{thm: HL strengthened}, as referred to in the last paragraph of Section~\ref{sec:model-theoretic HL}}\label{fig: XGvsXF}
\end{figure}	
 The mapping $\phi$ in~\eqref{eq:isoHL} is defined on the entire subset $Z:=\bigcup_{g\in J}Xg$. 
Its composition with the canonical morphism $\bX_\mP
F\twoheadrightarrow \bX_\mP G$ is the identity mapping on
$Z$. Hence $\phi$ is a bijection $Z\to \phi(Z)$. While the inverse
mapping $\phi\inv$ is an injective homomorphism on the substructure of
$\bX_\mP F$ induced on $\phi(Z)$, the mapping $\phi$ itself cannot be
guaranteed to be a homomorphism defined on every weak substructure of $\bX_\mP G$
on base set $Z$: it is just 
guaranteed to respect the relational tuples of $Z$ that are 
completely contained in some of the translates $Xg$ for some $g\in J$.  {However, this does not matter if we ignore the relational structure. Doing so, we get a strengthened version of Proposition~\ref{prop:abstract-extension}.}

 {\begin{Prop}[strengthening of Proposition~\ref{prop:abstract-extension}]\label{prop:strengthening-abstract-extension}
  For any finite set $X$ endowed with a set of partial bijections $\tilde{A}$ and every positive integer $n$ there exists a finite $A$-generated group $G=F/N$ such that for every choice of elements $g_1,\dots, g_n\in G$ there exist $f_1,\dots,f_n\in F$ such that $g_i=f_iN$ for all $i$ and the set of pairs
  \begin{equation}\label{eq:partial function=injective}
 \{(xg_i,xf_i)\colon x\in X, i=1,\dots, n\}     
  \end{equation}
is a partial function $X_\mathcal{P}G\to X_\mathcal{P}F$. In  particular,~\eqref{eq:partial function=injective} is a partial bijection whose inverse is induced by the canonical mapping $X_\mathcal{P}F\twoheadrightarrow X_\mathcal{P}G$.  
\end{Prop}}
 {One can interpret Propositions~\ref{prop:abstract-extension} and~\ref{prop:strengthening-abstract-extension} in terms of the hypergraphs \[\mathbb{U}=(U;\{Xg_i\colon i=1,\dots, n\})\mbox{ and }\mathbb{V}=(V;\{Xf_i\colon i=1,\dots,n\}),\] 
where $U=\bigcup_{i=1}^n Xg_i$ and $V=\bigcup_{i=1}^n Xf_i$, as follows. 
The condition formulated in Proposition~\ref{prop:abstract-extension} guarantees that the assignment $xg_i\mapsto xf_i$ ($x\in X, i=1,\dots,n)$ determines a (surjective) hypergraph homomorphism $\mathbb{U}\twoheadrightarrow\mathbb{V}$ while the same assignment in the context of Proposition~\ref{prop:strengthening-abstract-extension} provides a hypergraph isomorphism $\mathbb{U}\to \mathbb{V}$ whose inverse is induced by the canonical map $X_{\mathcal P}F\twoheadrightarrow X_{\mathcal P}G.$}

A final remark: in the proof of Theorem~\ref{thm: HL strengthened} we had to
choose, for the vertex $o$, an element $f_o$ in $F$ such that
$\ol{f_o}=o$.
There are many such choices and the set of all choices is given by
$f_oN$ where $N\unlhd F$ is the normal subgroup of $F$ for which
$F/N=G$. Two such choices $f_o$ and $f'_o$ are connected
via the equality $f'_o=f_ok$ for some $k\in N$. This leads to new choices for all
elements $f'_g=[p]_Ff'_o=[p]_Ff_ok=f_gk$. So, for every $k\in N$,
the substructure $\wh{\bT}k$ of $\bX_\mP F$
serves the same purpose as the original $\wh{\bT}$.
These substructures of $\bX_\mP F$
are related by the automorphisms $y\mapsto yk$ for $k\in N$. Figure~\ref{fig: XGvsXF} gives an overview over the situation.

\subsubsection{Group-theoretic Herwig--Lascar} Finally we formulate strengthened forms of the group-theoretic versions presented in Section~\ref{sec: group-theoretic version}. They are consequences of the following lemma.
\begin{Lemma}
    Let $\mE=(V\cup K;\alpha,\omega,\inv)$ be a finite connected graph and $\ell\colon K\to 2^F$, $k\mapsto H_kg_k$ a coset labelling. Then there exists a finite $A$-generated group $H$ such that, for every  $H$-commuting word labelling $l\colon K\to F$ with $l(k)\in \ell(k)$ for all $k\in K$, there exists an $F$-commuting word labelling $l'\colon K\to F$ with $l'(k)\in \ell(k)$ for all $k\in K$ such that, in addition, $[l(k)]_H=[l'(k)]_H$ for every $k\in K$.
\end{Lemma}
\begin{Rmk} \rm The statement without the last ``in addition'' condition follows from Ash's Theorem. The strengthening  consists entirely of the seemingly innocent additional condition saying that ``$[l(k)]_H=[l'(k)]_H$ for every $k\in K$'', that is, that the two labellings are $H$-related.    
\end{Rmk}
\begin{proof}
    For every geometric edge $e=k^\bullet\in K^\bullet$ let $\mS_e$ be the plain 
    Stallings graph of the coset $H_kg_k$ and let $\mS:=\bigsqcup_{e\in K^\bullet} \mS_e$ be the disjoint union of these graphs. Let $H$ be a finite $A$-generated group as in Theorem~\ref{thm:otto-ash} for $n=|V|$ and $M=\cT(\mS)$. Let $l\colon K\to F$ be an $H$-commuting word labelling of $\mE$ with $l(k)\in \ell(k)$ for every $k$. Let $l'\colon K\to F$ be an $M$- and $H$-related relabelling that commutes over $F$ according to Theorem~\ref{thm:otto-ash}. We only need to argue that $l'(k)\in \ell(k)$ for every $k$. Now $l(k)\in \ell(k)=H_kg_k$ implies that $\iota\cdot [l(k)]_M=\tau$ holds in $\mS_{k^\bullet}$ where $\iota$ and $\tau$ are the initial and terminal vertices of the Stallings graph $\mS_{H_kg_k}=(\iota, \mS_{k^\bullet},\tau)$. Since $[l(k)]_M=[l'(k)]_M$ it follows that $\iota\cdot [l'(k)]_M=\tau$ also holds in $\mS_{k^\bullet}$ which in turn implies that $l'(k)\in H_kg_k=\ell(k)$.
\end{proof}
We are now ready to proceed as in Section~\ref{sec: group-theoretic version}. We replace right cosets with left cosets and get strengthened versions of Propositions~\ref{prop:HLgroupformulationABO},~\ref{prop:ashalmeida} and~\ref{prop:HLgroupformulation}. We explicitly formulate the strengthened versions of the latter two. Note that, in both cases, the strengthening entirely consists of the seemingly innocent additional condition that the solutions $y_i$ belong to the cosets $x_iN$. 
\begin{Prop}[strengthened version of Proposition~\ref{prop:ashalmeida}]\label{prop:ashalmeidastrengthened}
    Let  $\mE=(V\cup K; \alpha,\omega,\inv)$ be a finite graph and $\ell\colon V\cup K\to 2^F$, $p\mapsto g_pH_p$ an extended coset labelling and let $(X_p)_{p\in V\cup K}$ be variables. Then there exists a finite index normal subgroup $N$ of $F$ such that: if the system 
    \[X_{\alpha k}X_k=X_{\omega k}\mbox{ for all }k\in K\]
    has a solution $(x_p)_{p\in V\cup K}$ in $F^{V\cup K}$ subject to the constraints $X_p\in g_pH_pN$ for all $p\in V\cup K$, then it even has a solution $(y_p)_{p\in V\cup K}$ subject to the constraints $X_p\in g_pH_p$ (all $p\in V\cup K$), such that, in addition, $y_p\in x_pN$ for all $p\in V\cup K$. The normal subgroup $N$ only depends on the size $|V|$ of $V$ and the collection of cosets $\{g_pH_p\colon p\in V\cup K\}$. 
\end{Prop}
\begin{Prop}[strengthened version of Proposition~\ref{prop:HLgroupformulation}]\label{prop:HLgroupformulationstrengthened} Let $I$ and $L$ be finite index sets, $(X_i)_{i\in I}$ variables, $(g_l)_{l\in L}$ elements of $F$ and $(H_l)_{l\in L}$ finitely generated subgroups of $F$; let $S\subseteq I\times L$ and $T\subseteq I^2\times L$. Then there exists a finite index normal subgroup $N$ of $F$ such that: if the system of constraints
\[X_i\in g_lH_lN\mbox{ for all } (i,l)\in S\mbox{ and }X_j\inv X_i\in g_lH_lN \mbox{ for all } (j,i,l)\in T\]
has a solution $(x_i)_{i\in I}\in F^I$, then even the system  
\[X_i\in g_lH_l\mbox{ for all } (i,l)\in S\mbox{ and }X_j\inv X_i\in g_lH_l \mbox{ for all } (j,i,l)\in T\]
has a solution $(y_i)_{i\in I}\in F^I$ such that, in addition, $y_i\in x_iN$ for every $i\in I$. The normal subgroup $N$ only depends on the size $|I|$ of $I$ and the collection of cosets $\{g_lH_l\colon l\in L\}$.
\end{Prop}

\section{{Comparing Section~\ref{sec: revisited} with Section~\ref{sec:strengthened}}}\label{sec:strictly stronger}
In view of Lemma~\ref{lem:main lemma} and the discussion after its proof, the essence of  Theorem~\ref{thm:otto-ash} 
 may be summarised as in Proposition~\ref{prop:otto-ash-pure} below, seemingly in a purely group-theoretic manner and without the use of inverse monoids. Recall that for a finite $A$-generated group $Q$ and a word $x$ over $\til A$, $\pi_1^\mQ(x)$ denotes the path starting at $1_Q$ in the Cayley graph $\mQ$ of $Q$, considered as a word over $Q\times \til A$, and also recall the definition of the word content $\co(\underline{\phantom{v}})$ of a word over $Q\times\til A$ according to Definition~\ref{def:word content}.
\begin{Prop}\label{prop:otto-ash-pure}
	For every finite $A$-generated group $Q$ and every positive integer $n$ there exists a finite $A$-generated group $H$ expanding $Q$ such that the following holds. For every finite connected graph $\mE$ on at most $n$ vertices and every $H$-commuting labelling $v\colon \mathrm{K}(\mE)\to \til A^*$ there exists a labelling $w\colon \mathrm{K}(\mE)\to \til A^*$ that commutes over $F$ and such that
		for every edge $k$ of $\mE$
	\begin{enumerate}
		\item $[v(k)]_H=[w(k)]_H$ and
		\item $\co(\pi_1^\mQ(v(k)))\supseteq \co(\pi_1^\mQ(w(k)))$.
	\end{enumerate}
\end{Prop}
As we have seen, Theorem~\ref{thm:otto-ash} is essentially 
the same as Proposition~\ref{prop:otto-ash-pure}. 
The condition occurring in Lemma~\ref{lem:main lemma} concerning the absence of loop edges is of proof-technical nature and not relevant for the validity of the claim. It is interesting that the essence of Theorem~\ref{thm: ash} 
can  be formulated in a similar inverse-monoid-free fashion. Here we refer to Definition~\ref{def: critical group}, Lemma~\ref{lem: key lemma RZ->Ash} and Proposition~\ref{prop:mytype2}.

\begin{Prop}\label{prop:ash-pure}
	For every finite $A$-generated group $Q$ and every positive integer $n$ there exists a finite $A$-generated group $H$ expanding $Q$ such that the following holds. For every finite connected graph $\mE$ on at most $n$ vertices and every $H$-commuting labelling $v\colon \mathrm{K}(\mE)\to \til A^*$ there exists a labelling $w\colon \mathrm{K}(\mE)\to \til A^*$ that commutes over $F$ and such that for every edge $k$ of $\mE$
	\begin{enumerate}
		\item $[v(k)]_Q=[w(k)]_Q$ and
		\item $\co(\pi_1^\mQ(v(k)))\supseteq \co(\pi_1^\mQ(w(k)))$.
	\end{enumerate}
\end{Prop}
Again, the conditions on the absence of loop edges in $\mE$ in Definition~\ref{def: critical group} 
is of proof-technical nature and not relevant for the formulation of the essence. According to Ash's proof in~\cite{Ash} (cf.~the discussion at the beginning of Section~\ref{sec:groups for Ash}), the group $H$ does not only depend on $Q$ and the number of vertices of $\mE$, but also on the cyclomatic number of $\mE$. 
The proof in~\cite{mytype2} of Proposition~\ref{prop:mytype2} by use of iterated expansions of the form $G\mapsto{}^{\mathbf{Ab}_m}G$ also proceeds by induction on the maximal length of a cycle subgraph of $\mE$  \textsl{and} the  cyclomatic number of $\mE$. However,  it also shows that \textsl{all} multiple edges can be handled by just one iteration. 
This implies that for any connected graph on $n$ vertices at most $(n-1)\binom{n-1}{2}+1$ iterations of expansions of the form 
$G\mapsto{}^{\mathbf{Ab}_m}G$ 
are needed to obtain the desired group $H$ from $Q$ (cf.~the formulation of Proposition~\ref{prop:mytype2}). {Hence the formulation of Proposition~\ref{prop:ash-pure} without mention of the cyclomatic number  is justified, although Ash's Theorem only states a slightly weaker assertion.}

{If we are given an $A$-generated group $Q$ and a positive integer $n$,  
then for any group $H$ as in Proposition~\ref{prop:ash-pure}, \textsl{every} expansion $K$ of $H$ fulfils the same r\^ole. 
In other words, the class of all groups $H$ serving for given $Q$ and $n$ as in Proposition~\ref{prop:ash-pure} is closed under  
expansions. In particular, the group $H$ 
can always be chosen to be \emph{relatively free} (that is, \textsl{every} mapping $A\to H$ extends to an endomorphism $H\to H$). In contrast, the class of all groups $H$ serving for  given $Q$ and $n$ as in Proposition~\ref{prop:otto-ash-pure} is \textsl{not} closed under expansions (because of item (1)).
Nevertheless,  it would be interesting whether the group $H$ in Proposition~\ref{prop:otto-ash-pure}  can be chosen to be relatively free in general.}\footnote{{This should be related to the discussion after the proof of Theorem~2.4 in~\cite{ABO}.}}

For the following discussion recall the definition of the Margolis--Meakin expansion $M(Q)$ of the $A$-generated group $Q$ presented after the proof of Lemma~\ref{lem:main lemma}. Proposition~\ref{prop:ash-pure} can be interpreted in the sense that $H$ is an \textsl{approximation of the free group $F$ for the inverse monoid $M(Q)$ with respect to the graph $\mE$}.\footnote{In fact, this and the following arguments concern uniformly all finite graphs whose connected components have at most $n$ vertices where $n$ is the parameter occurring in the proposition.} Indeed, conditions (1) and (2) of Proposition~\ref{prop:ash-pure} say that the inequality $[v(k)]_{M(Q)}\le[w(k)]_{M(G)}$ holds for all edges $k$ and the labellings $v$ and $w$ as mentioned.
	We then get the equality $[v(k)]_{M(Q)}=[w'(k)]_{M(Q)}$
	for the labelling $w'(k):=v(k)v(k)^{-1}w(k)$, which clearly still commutes over $F$.
In other words, any $H$-commuting labelling $v$ of $\mE$ can be replaced by an $F$-commuting labelling $w'$ of $\mE$ that cannot be distinguished by $M(Q)$ from the original labelling $v$. 
In this sense, based on just the graph $\mE$, the inverse monoid $M(Q)$ cannot recognise any difference between the finite group $H$ and the free group $F$. 
Or, as seen from the other side, 
the finite group $H$ \textsl{simulates} the behaviour of the free group $F$ on the graph $\mE$ for the inverse monoid $M(Q)$. 

Since $H$ is an expansion of $Q$, condition (1) of Proposition~\ref{prop:otto-ash-pure} implies $[v(k)]_Q=[w(k)]_Q$, so at first we have the same situation as in Proposition~\ref{prop:ash-pure}: $H$ still can be interpreted as an approximation of $F$ for  the inverse monoid $M(Q)$ with respect to the graph $\mE$. But the same condition (1) also encapsulates the essential  difference between Proposition~\ref{prop:otto-ash-pure} and Proposition~\ref{prop:ash-pure}: while the group $H$ of the stronger version guarantees that the old labelling $v$ and the new labelling $w$  are 
$H$-related, the weaker version only ensures that these labellings are $Q$-related.\footnote{{This should be compared with the discussion in Section~2.7 of \cite{ABO}.}} 
Actually, the series of expansions $Q\mapsto Q_1\mapsto \cdots \mapsto Q_{N-1}\mapsto Q_N=:H$, where for every $i$, $Q_{i+1}={}^{\mathbf{Ab}_{m_i}}Q_i$ for some $m_i\ge 2$, leading to the expansion $H$ of Proposition~\ref{prop:ash-pure}, even guarantees that the labellings $v$ and $w$ are $Q_{N-1}$-related. This is stronger than being $Q$-related, but still weaker than being $Q_N$-related,
which would be required for a proof of Proposition~\ref{prop:otto-ash-pure} (a similar observation can be made for Ash's iterated expansions discussed in Section~\ref{sec:groups for Ash}). In a sense, the series of expansions used to prove Proposition~\ref{prop:ash-pure} 
always remains ``one step away'' 
from delivering a proof of Proposition~\ref{prop:otto-ash-pure}. This behaviour has been observed by the second author in his doctoral dissertation~\cite{bitterlichdiss}. 
Nevertheless, it could a priori still be possible 
that $v$ and $w$ are $H$-related for $H:=Q_N$ 
for some sufficiently large $N$. 
The second author conjectured that this is not the case for any $N$. As a consequence of Theorem~\ref{thm: otto-ash vs ash} below, this conjecture turns out to be true.  This also provides some  evidence that Theorem~\ref{thm:otto-ash} is substantially stronger than Theorem~\ref{thm: ash}.

The remainder of the present section is devoted to prove Theorem~\ref{thm: otto-ash vs ash} below. We start with some definitions. 
For a finite $A$-generated group $Q$ let $H(Q,n)$ be the $A$-generated expansion of $Q$ defined in \eqref{eq:group H}, with the intermediate $(Q\times A)$-generated group $G$ described in Section~\ref{sec:group G} being $n$-acyclic. Moreover recall, for an integer $m\ge 2$, the expansion $Q\mapsto {}^{\mathbf{Ab}_m}Q$ introduced at the end of Section~\ref{sec:groups for Ash} in the context of Proposition~\ref{prop:mytype2}. For integers $m_1,\dots, m_t\ge 2$ denote by ${}^{\mathbf{Ab}_{m_t}\cdots\mathbf{Ab}_{m_1}}Q$ the expansion of $Q$ obtained by iterated expansions $Q\mapsto {}^{\mathbf{Ab}_{m_1}}Q\mapsto {}^{\mathbf{Ab}_{m_2}}({}^{\mathbf{Ab}_{m_1}}Q)\mapsto\cdots \mapsto {}^{\mathbf{Ab}_{m_t}}(\cdots({}^{\mathbf{Ab}_{m_1}}Q)\cdots)$.

An \emph{expansion series} $\mathbf{exps}$  is a mapping $Q\mapsto \mathbf{exps}(Q)$ that, for every $A$, assigns to every finite $A$-generated group $Q$ a set $\mathbf{exps}(Q)$ of finite $A$-generated expansions of $Q$. Examples of expansion series that are significant in the present context are $\mathbf{solk}$ and $\mathbf{conc}$ where:
\begin{itemize}
	\item[i)] $\mathbf{solk}(Q)=\big\{{}^{\mathbf{Ab}_{m_t}\cdots\mathbf{Ab}_{m_1}}Q\colon t\ge 1,m_1,\dots,m_t\ge 2\big\} $
	\item[ii)]$\mathbf{conc}(Q)=\{H(Q,n)\colon n\ge 5\}$. 
\end{itemize}
Here $\mathbf{solk}$ stands for \emph{solvable kernel} while $\mathbf{conc}$ stands for \emph{connected content}.
\begin{Def}\rm
	An expansion series $\mathbf{exps}$ is \emph{sufficient} for Theorem~\ref{thm: ash} 
if, for every finite alphabet $A$,  every $A$-generated group $Q$ that extends an $A$-generated inverse monoid $M$ admits, 
for every finite graph $\mE$, an expansion $H \in \mathbf{exps}(Q)$ 
satisfying the requirements of Theorem~\ref{thm: ash}.
Sufficiency for Theorem~\ref{thm:otto-ash} is analogously defined. 
\end{Def}
For example, the expansion series $\mathbf{conc}$ is sufficient for Theorem~\ref{thm:otto-ash} (the proof of this claim is the content of Section~\ref{sec:proofs}) 
and hence is sufficient also for Theorem~\ref{thm: ash}. Proposition~\ref{prop:mytype2} tells us that the expansion series $\mathbf{solk}$ is sufficient for Theorem~\ref{thm: ash}. The main point
of the present section is to show that $\mathbf{solk}$ is not sufficient for Theorem~\ref{thm:otto-ash}. This may be interpreted as saying that Theorem~\ref{thm:otto-ash} is strictly stronger (or harder) than Theorem~\ref{thm: ash}. Indeed, the expansion series $\mathbf{solk}$ provides a tool that enables one to prove Theorem~\ref{thm: ash} while it is not possible to prove Theorem~\ref{thm:otto-ash} by use of this instrument. 
\begin{Thm}\label{thm: otto-ash vs ash} The expansion series $\mathbf{solk}$ is sufficient for Theorem~\ref{thm: ash}, but is not sufficient for Theorem~\ref{thm:otto-ash}. 
\end{Thm}

We start to recall some concepts and results from the theory of group varieties; 
a standard reference for this topic is the monograph \cite{HNeumann}. For two varieties $\mathbf{U}$ and $\mathbf{V}$, their product $\mathbf{UV}$ consists of all groups $G$ having a normal subgroup $N\in \mathbf{U}$ for which the quotient $G/N$ belongs to $\mathbf{V}$. The class $\mathbf{UV}$ is again a variety and for any three varieties $\mathbf{U},\mathbf{V},\mathbf{W}$ we have $(\mathbf{UV})\mathbf{W}=\mathbf{U}(\mathbf{VW})$. For a variety $\mathbf{U}$ and a group $G$ the \emph{$\mathbf{U}$-verbal subgroup of $G$} denoted by $\mathbf{U}(G)$ is the (unique) smallest normal subgroup $N$ of $G$ for which $G/N$ belongs to $\mathbf{U}$. The term ``verbal'' in this context stems from the fact that $\mathbf{U}(G)$ is comprised of all values of words $w$ such that $w=1$ is an identity
of $\mathbf{U}$.\footnote{Here ``value of a word'' has a slightly different meaning than in Section~\ref{sec:prelims}.} For example, for the variety $\mathbf{Ab}$ of all Abelian groups, the $\mathbf{Ab}$-verbal subgroup of a group $G$ is the subgroup of $G$ generated by all \emph{commutators}, that is, all elements of $G$ of the form $ghg\inv h\inv$ for $g,h\in G$. In the context of 
formations of finite groups the term \emph{$\mathbf{U}$-residual} is used instead of $\mathbf{U}$-verbal subgroup. For two varieties $\mathbf{U}$ and $\mathbf{V}$ and a group $G$ we have that 
\begin{equation}\label{eq:composed verbal subgroup}
	\mathbf{U}(\mathbf{V}(G))=(\mathbf{UV})(G).
\end{equation}
Next, let $Q=F/N$ be a not necessarily finite $A$-generated group and $\mathbf{U}$ be a variety of groups. The \emph{$\mathbf{U}$-universal expansion} ${}^{\mathbf{U}}Q$ of $Q$ is defined to be 
\begin{equation}
{}^{\mathbf{U}}Q:=F/\mathbf{U}(N).	
\end{equation}
The $\mathbf{U}$-universal expansion of $Q$ is the largest $A$-generated expansion $H$ of $Q$ such that the kernel of the canonical morphism $H\twoheadrightarrow Q$ belongs to $\mathbf{U}$. That is, the canonical morphism ${}^{\mathbf{U}}Q\twoheadrightarrow Q$ factors through any such expansion $H$ of $Q$ with kernel in $\mathbf{U}$. For two varieties $\mathbf{U}$ and $\mathbf{V}$ and an $A$-generated group $Q$ we have ${}^{\mathbf{U}}\big({}^{\mathbf{V}}\!Q\big)={}^{\mathbf{UV}}Q$ --- this follows immediately from~\eqref{eq:composed verbal subgroup}. In particular, the notation of the iterated expansion  ${}^{\mathbf{Ab}_{m_t}\cdots\mathbf{Ab}_{m_1}}Q$ chosen at the beginning of the present section is consistent with this: ${}^{\mathbf{Ab}_{m_t}\cdots\mathbf{Ab}_{m_1}}Q$ denotes the universal expansion of $Q$ with respect to the product variety $\mathbf{Ab}_{m_t}\cdots\mathbf{Ab}_{m_1}$ where $\mathbf{Ab}_{m_i}$ denotes the variety defined by the identities $xy=yx$ and $x^{m_i}=1$, which consists of all Abelian groups of exponent dividing $m_i$. 
A variety $\mathbf{U}$ is \emph{locally finite} if every finitely generated member of $\mathbf{U}$ is finite; the product $\mathbf{UV}$ of two locally finite varieties $\mathbf{U}$ and $\mathbf{V}$ is again locally finite. In general, the $\mathbf{U}$-universal expansion ${}^{\mathbf{U}}Q$ need not be finite, even for finite $Q$. However, ${}^{\mathbf{U}}Q$ is finite for finite $Q$ provided that $\mathbf{U}$ is locally finite. We note that the varieties  $\mathbf{Ab}_{m_t}\cdots \mathbf{Ab}_{m_1}$ we deal with are all locally finite.

In the following we  give a description of groups of the form ${}^{\mathbf{U}}Q$, similar to the construction of the group $H$ in~\eqref{eq:group H}. In this context we note that if $Q$ happens to be the $A$-generated free object in the variety $\mathbf{V}$ then ${}^{\mathbf{U}}Q$ is the $A$-generated free object in the variety $\mathbf{UV}$. This fact will be exploited later. So let $Q$ be an $A$-generated group and $E=Q\times A$ the set of positive edges of the Cayley graph $\mQ$ of $Q$. Let $G$ be the $E$-generated free object in $\mathbf{U}$. Now $Q$ acts 
as a group of permutations
on $E$ via $(g,a)\mapsto {}^q(g,a)=(qg,a)$ for $q\in Q$ and $(g,a)\in E$. Since $G$ is free in $\mathbf{U}$, this action extends to an action of $Q$ on $G$ by automorphisms on the left, denoted $\gamma\mapsto {}^q\gamma$ for $q\in Q$ and $\gamma\in G$, and we may form the semidirect product $G\rtimes Q$. Let $H$ be the following $A$-generated subgroup of $G\rtimes Q$:
\begin{equation}\label{eq:group ^UQ}
H:=\langle([(1_Q,a)]_G,[a]_Q)\colon a\in A \rangle\le G\rtimes Q.
\end{equation}
Just as in the context of~\eqref{eq:group H}, for a word $x$ over $\til A$, we have $[x]_H=([\pi_1^\mQ(x)]_G,[x]_Q)$ (see (2.4) in~\cite{ABO}). The following is folklore (see e.g.~\cite{ASz}, or~\cite[Chapter 10]{almeida:book} for the discussion of a closely related topic).
\begin{Prop}
	The group $H$ defined in~\eqref{eq:group ^UQ} is a realisation of ${}^{\mathbf{U}}Q$. In particular, if $Q$ is the $A$-generated free object in the variety $\mathbf{V}$, then $H$ is a realisation of the $A$-generated free object in $\mathbf{UV}$.
\end{Prop}

In  Lemma~\ref{lem:crucialforsol} below we shall discuss two crucial identities, 
 the first of which holds in the variety $\mathbf{Ab}_{m_t}\cdots \mathbf{Ab}_{m_1}$ while the second one fails. These identities will be \textsl{the} essential ingredient for proving the desired result. For elements $w_1,\dots, w_n$ of an involutory monoid define $[w_1,\dots,w_n]$ inductively by $[w_1]:=w_1$, $[w_1,w_2]:=w_1w_2w_1\inv w_2\inv$ and $[w_1,\dots, w_n]:=[[w_1,\dots w_{n-1}], w_n]$. We first need the following lemma.

\begin{Lemma}\label{lem: commutators}
	Let $S,Q$ be groups, $S$ abelian and such that $Q$ acts on $S$ by automorphisms on the left. Let $a=(h,g), b=(x,1_Q)\in S\rtimes Q$. Then the following hold:
	\begin{enumerate}
		\item $[a,b]$ is independent of $h$,
		\item $[a,b]\ne 1_{S\rtimes Q}$ if and only if ${}^gx\ne x$.
	\end{enumerate} 
\end{Lemma} 
\begin{proof}
	By direct calculation we get
	\begin{align*}
		aba\inv b\inv &=(h,g)(x,1_Q)({}^{g\inv}h\inv,g\inv)(x\inv,1_Q) \\
	&=	(h\,{}^g\!x,g)({}^{g\inv}h\inv\cdot{}^{g\inv}\!x\inv,g\inv)\\
	&= (h\,{}^g\!xh\inv x\inv,1_Q)=({}^g\!xx\inv,1_Q)
	\end{align*}
since $S$ is abelian. We see that $({}^g\!xx\inv,1_Q)$ is independent of $h$ and coincides with the identity element $(1_S,1_Q)$ if and only if ${}^g\!x=x$.
\end{proof}

We come to the result about identities in $\mathbf{Ab}_{m_t}\cdots\mathbf{Ab}_{m_1}$ announced above. The first author is grateful to R.~Bryant and A.~Krasilnikov for having brought this crucial result to his attention back in the year 2002.
\begin{Lemma}\label{lem:crucialforsol}
	Let $t\ge 1$ and $m_1,\dots,m_t\ge 2$ and let $H$ be the free object on three letters $x,y,z$ in the variety $\mathbf{Ab}_{m_t}\cdots \mathbf{Ab}_{m_1}$. Then the following hold in $H$:
	\begin{enumerate}
		\item $[xyx\inv,z^{m_1},z^{m_1m_2},\dots,z^{m_1\cdots m_{t-1}}]= [y,z^{m_1},z^{m_1m_2},\dots,z^{m_1\cdots m_{t-1}}]$	
		\item $[y,z^{m_1},z^{m_1m_2},\dots,z^{m_1\cdots m_{t-1}}]\ne 1_H$.
	\end{enumerate}
\end{Lemma}
\begin{proof}
	 The proof is by induction on $t$. For $t=1$ this reduces to $xyx\inv = y$ and $y\ne 1_H$, which is trivial since $\mathbf{Ab}_{m_1}$ is Abelian and non-trivial (that is, $m_1\ge 2$). So let $t\ge 2$ and suppose that the claim is true for $t-1$. We realise $H$ as in~\eqref{eq:group ^UQ} for $Q$ the $A=\{x,y,z\}$-generated free object in $\mathbf{Ab}_{m_{t-1}}\cdots \mathbf{Ab}_{m_1}$ and $G$ the $(Q\times A)$-generated free object in $\mathbf{Ab}_{m_t}$. Let $w=[y,z^{m_1},z^{m_1m_2},\dots,z^{m_1\cdots m_{t-2}}]$, $w'= [xyx\inv,z^{m_1},z^{m_1m_2},\dots,z^{m_1\cdots m_{t-2}}]$ and $v=z^{m_1\cdots m_{t-1}}$  considered as words over $\til A$. We evaluate these words in $H$. By the inductive assumption,
	  $w$ and $w'$ 
evaluate identically 
in $Q$; hence $a:=[w]_H=(h,g)$ and $a':=[w']_H=(h',g)$ with $[w]_Q=g=[w']_Q$ for some $h,h'\in G$ and $g\in Q$. All members of the variety $\mathbf{Ab}_{m_{t-1}}\cdots \mathbf{Ab}_{m_1}$ have exponents dividing $m_1\cdots m_{t-1}$, hence $[v]_Q=1_Q$ so that $b:=[v]_H=(k,1_Q)$ for some $k\in G$. By Lemma~\ref{lem: commutators}~(1) we find that $[a,b]=[a',b]$ which proves assertion (1) of the Lemma.
	 
	 In order to prove item~(2) we need to show that ${}^gk\ne k$. For this we note that, by the inductive assumption, $[w]_Q=g\ne 1_Q$. 
     It follows 
     from the structure of the word $w$ 
     that $g$ is not in the subgroup of $Q$ generated by $z$, or, in other words, the content of $g$ is $\{y,z\}$. (Note that, since $Q$ is free in $\mathbf{Ab}_{m_{t-1}}\cdots \mathbf{Ab}_{m_1}$, it is in particular retractable, 
     and therefore has a content function.) On the other hand, $k$ belongs to the subgroup of $G$ generated by all elements of the form $(q,z)$ with $q$ in the subgroup of $Q$ generated by $z$. More precisely,
	 \[k=[(1_Q,z)([z]_Q,z)([z^2]_Q,z)\cdots ([z^{m_1\cdots m_{t-1}-1}]_Q,z)]_G.\]
	 Since the order of $[z]_Q$ is $m_1\cdots m_{t-1}$, the elements $([z^i]_Q,z)$ are pairwise distinct, hence none of them 
     cancel with each other, and in particular $k\ne 1_G$. Since $g$ has content $\{y,z\}$, $[z^i]_Q\ne g[z^j]_Q$ for all $i, j$ and the latter elements $g[z^j]_Q$ are also pairwise distinct. From this we get ${}^gk\ne k$. 
\end{proof}

 In order to continue we recall the
 \emph{Margolis--Meakin expansion} $M(Q)$ of an $A$-generated group  $Q$ 
 introduced in Section~\ref{sec:proofs}, which may be described as follows. 
 Recall that the group $Q$ acts on its Cayley graph $\mQ$
 by left multiplication. For a subgraph $\mK$ of $\mQ$ and $g\in Q$ we denote by ${}^g\mK$ the image of $\mK$ under the action of $g$.
 For a given $A$-generated group $Q$, the Margolis--Meakin expansion $M(Q)$
 consists of all pairs $(\mK,g)$ where $g\in Q$ and $\mK$ is a finite connected subgraph of the Cayley graph $\mQ$ of $Q$ containing the vertices $1_Q$ and $g$. Endowed with the multiplication
 \[(\mK,g)(\mL,h)=(\mK\cup{}^g\mL,gh)\] and involution \[(\mK,g)\inv=({}^{g\inv}\mK,g\inv)\] the set $M(Q)$ becomes an $A$-generated inverse monoid with identity element $(\{1_Q\},1_Q)$,
 and with the assignment function
 \[A\to M(Q),\quad a\mapsto (\langle (1_Q,a)\rangle,[a]_Q).\]
 The value of the word $p\in \til{A}^*$ in $M(Q)$ is
 \[ [p]_{M(Q)}=(\langle \pi_1^\mQ(p)\rangle,[p]_Q),\] where $\pi_1^\mQ(p)$ is the path in $\mQ$ starting at $1_Q$ and having label $p$; the natural partial order on $M(Q)$ is given by
 \[(\mK,g)\le(\mL,h)\mbox{ if and only if }\mK\supseteq \mL\mbox{ and }g=h.\]
 
 \begin{figure}[ht]
 	\begin{tikzpicture}[scale=1.4]
 		\filldraw(1,1) circle (1pt);
 		\filldraw(1,-1) circle (1pt);
 		\filldraw(-1,1) circle (1pt);
 		\filldraw(-1,-1) circle (1pt);
 		\draw[->, >=stealth](1.02,-0.95) .. controls (1.25,-0.5) and (1.25,0.5) .. (1.02,0.95); 
 		\draw[<-,>=stealth](0.98,-0.95) .. controls (0.75,-0.5) and (0.75,0.5) .. (0.98,0.95); 
 		\draw[<-,>=stealth] (-1.05,-0.95) .. controls (-1.25,-0.5) and (-1.25,0.5) .. (-1.05,0.95); 
 		\draw[->,>=stealth] (-0.98,-0.95) .. controls (-0.75,-0.5) and (-0.75,0.5) .. (-0.98,0.95); 
 		\draw[<-,>=stealth] (-0.98,1.05) .. controls (-0.5,1.25) and (0.5,1.25) .. (0.95,1.05); 
 		\draw[->, >=stealth] (-0.95,0.98) .. controls (-0.5,0.75) and (0.5,0.75) .. (0.95,0.95); 
 		\draw[->, >=stealth] (-0.95,-1.05) .. controls (-0.5,-1.25) and (0.5,-1.25) .. (0.98,-1.05); 
 		\draw[<-,>=stealth] (-0.95,-0.98) .. controls (-0.5,-0.75) and (0.5,-0.75) .. (0.95,-0.98); 
 		\draw[below right](1,-1)node{$1_Q$};
 		\draw[right](1.15,0)node{$a$};
 		\draw[left](0.85,0)node{$a$};
 		\draw[left](-1.15,0)node{$a$};
 		\draw[right](-0.85,0)node{$a$};
 		\draw[above] (0,1.15)node{$b$};
 		\draw[below] (0,0.85)node{$b$};
 		\draw[below] (0,-1.15)node{$b$};
 		\draw[above] (0,-0.85)node{$b$};
 	\end{tikzpicture}
 \qquad\quad
 \begin{tikzpicture}[scale=1.4]
 	\filldraw(1,1) circle (1pt);
 	\filldraw(1,-1) circle (1pt);
 	\filldraw(-1,1) circle (1pt);
 	\filldraw(-1,-1) circle (1pt);
 	\draw[->, >=stealth](1.02,-0.95) .. controls (1.25,-0.5) and (1.25,0.5) .. (1.02,0.95); 
 	\draw[<-,>=stealth](0.98,-0.95) .. controls (0.75,-0.5) and (0.75,0.5) .. (0.98,0.95); 
 	\draw[<-,>=stealth] (-1.05,-0.95) .. controls (-1.25,-0.5) and (-1.25,0.5) .. (-1.05,0.95); 
 	\draw[->,>=stealth] (-0.98,-0.95) .. controls (-0.75,-0.5) and (-0.75,0.5) .. (-0.98,0.95); 
 	\draw[<-,>=stealth] (-0.98,1.05) .. controls (-0.5,1.25) and (0.5,1.25) .. (0.95,1.05); 
 	\draw[->, >=stealth] (-0.95,0.98) .. controls (-0.5,0.75) and (0.5,0.75) .. (0.95,0.95); 
 	\draw[->, >=stealth] (-0.95,-1.05) .. controls (-0.5,-1.25) and (0.5,-1.25) .. (0.98,-1.05); 
 	\draw[<-,>=stealth] (-0.95,-0.98) .. controls (-0.5,-0.75) and (0.5,-0.75) .. (0.95,-0.98); 
 	\draw[below right](1,-1)node{$1_Q$};
 	\draw[right](1.15,0)node{$e$};
 	\draw[left](0.85,0)node{$f$};
 	\draw[left](-1.15,0)node{$i$};
 	\draw[right](-0.85,0)node{$j$};
 	\draw[above] (0,1.15)node{$g$};
 	\draw[below] (0,0.85)node{$h$};
 	\draw[below] (0,-1.15)node{$c$};
 	\draw[above] (0,-0.85)node{$d$};
 \end{tikzpicture}
 \caption{Cayley graph $\mQ$}\label{fig:cayleyV4}	
 \end{figure}
 
 In the following we let $Q:=\langle a,b\colon a^2=b^2=(ab)^2=1\rangle$ be the Klein four-group and consider its Margolis--Meakin expansion $M(Q)$. Let $\mE$ be the graph consisting of two vertices $o\ne p$, two  edges $e_1,e_2\colon o\longrightarrow p$ and their inverses. It is clear that $Q$ extends $M(Q)$. We show that no expansion $H\in \mathbf{solk}(Q)$ satisfies the requirements of Theorem~\ref{thm:otto-ash} for $M=M(Q)$ and $\mE$.
\begin{figure}[ht]
\begin{tikzpicture}[scale=1.4]
	\filldraw(1,1) circle (1pt);
	\filldraw(1,-1) circle (1pt);
	\filldraw(-1,1) circle (1pt);
	\filldraw(-1,-1) circle (1pt);
	\draw[->, >=stealth](1.02,-0.95) .. controls (1.25,-0.5) and (1.25,0.5) .. (1.02,0.95); 
		\draw[<-,>=stealth] (-0.98,1.05) .. controls (-0.5,1.25) and (0.5,1.25) .. (0.95,1.05); 
	\draw[->, >=stealth] (-0.95,0.98) .. controls (-0.5,0.75) and (0.5,0.75) .. (0.95,0.95); 
	\draw[->, >=stealth] (-0.95,-1.05) .. controls (-0.5,-1.25) and (0.5,-1.25) .. (0.98,-1.05); 
	\draw[<-,>=stealth] (-0.95,-0.98) .. controls (-0.5,-0.75) and (0.5,-0.75) .. (0.95,-0.98); 
	\draw[below right](1,-1)node{$1_Q$};
	\draw[left](1.2,0)node{$e$};
	\draw[above] (0,1.15)node{$g$};
	\draw[below] (0,0.85)node{$h$};
	\draw[below] (0,-1.15)node{$c$};
	\draw[above] (0,-0.85)node{$d$};
	\draw(0,0)node{$\mathcal{A}_1$};
\end{tikzpicture}
\qquad\quad
\begin{tikzpicture}[scale=1.4]
	\filldraw(1,1) circle (1pt);
	\filldraw(1,-1) circle (1pt);
	\filldraw(-1,1) circle (1pt);
	\filldraw(-1,-1) circle (1pt);
	\draw[<-,>=stealth](0.98,-0.95) .. controls (0.75,-0.5) and (0.75,0.5) .. (0.98,0.95); 
	\draw[<-,>=stealth] (-0.98,1.05) .. controls (-0.5,1.25) and (0.5,1.25) .. (0.95,1.05); 
	\draw[->, >=stealth] (-0.95,0.98) .. controls (-0.5,0.75) and (0.5,0.75) .. (0.95,0.95); 
	\draw[->, >=stealth] (-0.95,-1.05) .. controls (-0.5,-1.25) and (0.5,-1.25) .. (0.98,-1.05); 
	\draw[<-,>=stealth] (-0.95,-0.98) .. controls (-0.5,-0.75) and (0.5,-0.75) .. (0.95,-0.98); 
	\draw[below right](1,-1)node{$1_Q$};
	\draw[right](0.8,0)node{$f$};
		\draw[above] (0,1.15)node{$g$};
	\draw[below] (0,0.85)node{$h$};
	\draw[below] (0,-1.15)node{$c$};
	\draw[above] (0,-0.85)node{$d$};
	\draw(0,0)node{$\mathcal{A}_2$};
\end{tikzpicture}
	\caption{The graphs $\mA_1$ and $\mA_2$}\label{fig:A_1 and A_2}
\end{figure}
 
 Suppose to the contrary that  
 for some $t\ge 1$, $m_1,\dots, m_t\ge 2$ the expansion $H:={}^{\mathbf{Ab}_{m_t}\cdots\mathbf{Ab}_{m_1}}\!Q$ satisfies the requirements of Theorem~\ref{thm:otto-ash}. We denote the positive edges of the Cayley graph $\mQ$ by the letters in $B:=\{c,d,e,f,g,h,i,j\}$ as indicated in Figure~\ref{fig:cayleyV4}, on the right-hand side. For the letters in $C:=\{x,y,z\}$ set
 \[w(x,y,z):=[xyx\inv, z^{m_1}, z^{m_1m_2},\dots,z^{m_1\cdots m_{t-1}}],\]
 considered as a word over $\til C$. Then the word $w(e,gh,dc)$ over $\til B$ forms a path in $\mQ$ closed at $1_Q$ and so does the word $w(f\inv,gh,dc)$. Here $w(e,gh,dc)$ denotes the word obtained from $w(x,y,z)$ by substituting $x\mapsto e, y\mapsto gh, z\mapsto dc$ and analogously for $w(f\inv, gh,dc)$. Hence there exist words $u_1,u_2$ over $\til A$ such that $w(e,gh,dc)=\pi_1^\mQ(u_1)$ and $w(f\inv,gh,dc)=\pi_1^\mQ(u_2)$.  As a first consequence we have that $[u_1]_M=(\mA_1,1_Q)$ and $[u_2]_M=(\mA_2,1_Q)$ for $\mA_1$ and $\mA_2$ as indicated in Figure~\ref{fig:A_1 and A_2}. 
 
 {Recall that $\mE$ is the graph consisting of two vertices $o\ne p$, two edges $e_1,e_2\colon o\longrightarrow p$ and their inverses  and consider the word labelling of $\mE$ defined by $\ell(e_i)=u_i$}.
 We represent $H$ as in~\eqref{eq:group ^UQ} with $G$ the $B$-generated free object in $\mathbf{Ab}_{m_t}\cdots\mathbf{Ab}_{m_1}$. From Lemma~\ref{lem: commutators} (1) it follows that in $G$ we have
\begin{align*}
[\pi_1^\mQ(u_1)]_G&=[w(e,gh,dc)]_G=[w(1,gh,dc)]_G\\ &=[w(f\inv,gh,dc)]_G=[\pi_1^\mQ(u_2)]_G    
\end{align*}
 hence $[u_1]_H=[u_2]_H$ since clearly $[u_1]_Q=[u_2]_Q$. Consequently, $\ell$ commutes over $H$. 
 
Suppose $v_1,v_2$ are words over $\til A$ such that $\ell'(e_i)=v_i$ provides an $F$-commuting relabelling of $\ell$ that is $M$-related to $\ell$. The latter means that $[u_i]_M=[v_i]_M$ and thus $[u_i]_M\le [\mathrm{red}(v_i)]_M$. Since the labelling $\ell'$ commutes over $F$ we have $\mathrm{red}(v_1)=\mathrm{red}(v_2)=:v$  and therefore $[v]_M\ge [u_1]_M, [u_2]_M$. It follows that $\langle \pi_1^\mQ(v)\rangle\subseteq \mA_1\cap \mA_2$.
The latter inclusion is only possible if $v$ is a power of $b$ (positive or negative). If $[u_i]_H=[v_i]_H$ then also $[u_i]_H=[v]_H$ since $v_1$ and $v_2$ have the same reduced form. 
It follows from Lemma~\ref{lem: commutators} (2) that
 $[w(1,gh,dc)]_G\ne 1_G$
 and in particular \[[w(1,gh,dc)]_G\ne 1_G=[w(1,1,dc)]_G.\] In other words, $g$ and $h$ belong to the content of $[w(1,gh,dc)]_G$, which is the first entry of $[u_1]_H$ as well as of $[u_2]_H$. On the other hand, if $v$ is a power of $b$ (positive or negative) then $g$ and $h$ cannot belong to the content of $[\pi_1^\mQ(v)]_G$, which is the first entry of $[v]_H$. Consequently, $[u_1]_H\ne[v]_H\ne[u_2]_H$, a contradiction. {In fact, we have shown that for no expansion $H\in\mathbf{solk}(Q)$ can the $A$-generated direct product $M(Q)\times H$ be an $F$-inverse cover of $M(Q)$, see \cite[Section 2.7]{ABO}.}

 The above arguments apply to the Margolis--Meakin expansion $M(Q)$ of any finite $A$-generated group $Q$ whose Cayley graph admits two disjoint cycle subgraphs having two arcs $\mB_1$ and $\mB_2$ connecting them such that the intersection 
 $\mB_1\cap \mB_2$ does not connect the two cycles.

\end{document}